\documentclass[12pt]{article}

\usepackage[margin=1in]{geometry}
\usepackage{amsmath,amsthm,amssymb,amsfonts,float}
\usepackage{listings}
\usepackage{pdfpages}
\usepackage{mathtools,slashed}
\usepackage{tikz}
\usetikzlibrary{patterns}
\usepackage{cancel}

\definecolor{colorqqqq}{rgb}{0.1,0.7,0}

\def\aaa{\zeta}
\def\EE{\mathcal{E}}

\definecolor{coloroftheorems}{rgb}{0.45,0.0,0.0}

\def\les{\lesssim}
\makeatletter
\newcommand{\subjclass}[2][2010]{%
  \let\@oldtitle\@title%
  \gdef\@title{\@oldtitle\footnotetext{#1 \emph{Mathematics subject classification.} #2}}%
}
\makeatother

\usepackage[colorlinks=true, pdfstartview=FitV, linkcolor=blue,citecolor=blue, urlcolor=blue]{hyperref}

\usepackage[abbrev,lite,nobysame]{amsrefs}
\usepackage{color}
\usepackage{esint}

\usepackage{mathtools,enumitem,mathrsfs}

\usepackage[compact]{titlesec}

\usepackage{comment}

\mathtoolsset{showonlyrefs=true}

\newcommand{\ep}{\epsilon}
\newcommand{\eps}{\epsilon}

\newcommand{\na}{\nabla}
\newcommand{\grad}{\nabla}

\newcommand{\norm}[1]{\left\|  #1 \right\|}
\newcommand{\abs}[1]{\left| #1 \right|}

\newcommand{\brak}[1]{\left\langle #1 \right\rangle} 
\newcommand{\br}[1]{\left\langle #1 \right\rangle} 

\newcommand{\rr}{\mathbb{R}}

\newcommand{\dee}{\mathrm{d}}

\DeclareMathOperator{\curl}{\mathrm{curl}}

\usepackage{bbm}

\newtheorem{theorem}{Theorem}[section]
\newtheorem{proposition}[theorem]{Proposition}

\newtheorem{lemma}[theorem]{Lemma}
\newtheorem*{lemma*}{Lemma}

\theoremstyle{definition}

\newtheorem{remark}[theorem]{Remark}

\def\ZZ{\mathbb{Z}}

\def\paren#1{\left(#1\right)}
\def\abs#1{\left|#1\right|}
\newcommand{\n}{\ensuremath{\nonumber}}
\newcommand{\pa}{\ensuremath{\partial}}

\def\simhe#1{}% {\color{cyan} #1}

\newcommand{\enorm}[1]{\left\lVert #1\right\rVert_{L^2}}

\newcommand{\myr}[1]{{\color{red} #1 }}%

\newcommand{\wt}{\widetilde}
\newcommand{\al}{\alpha}

\newcommand{\de}{\Delta}
\newcommand{\lan}{\langle}
\newcommand{\ran}{\rangle}
\newcommand{\lf}{\left}
\newcommand{\rg}{\right}

\newcommand{\mf}{\mathfrak}
\newcommand{\mh}[1]{}

\newcommand{\jj}{{[j]}}
\newcommand{\jjj}{{[j']}}
\newcommand{\mc}{\mathcal}
\newcommand{\rma}{\mathrm{ap}}
\newcommand{\rme}{\mathrm{re}}
\newcommand{\rmi}{\mathrm{in}}
\newcommand{\U}{\mathcal{U}}
\newcommand{\exb}[1]{\exp\left\{ #1 \right\}}
\allowdisplaybreaks

\newcommand{\nv}{\mathrm{Na}}

\renewcommand{\Re}{\textup{Re}}

\allowdisplaybreaks

\numberwithin{equation}{section}
\begin{document}

\title{Stability threshold of close-to-Couette shear flows \\ with no-slip boundary conditions in 2D} 
\author{ Jacob Bedrossian\thanks{\footnotesize Department of Mathematics, University of California, Los Angeles, CA 90095, USA \href{mailto:jacob@math.ucla.edu}{\texttt{jacob@math.ucla.edu}}} \and Siming He\thanks{Department of Mathematics, University of South Carolina, Columbia, SC 29208, USA \href{mailto:siming@mailbox.sc.edu}{\texttt{siming@mailbox.sc.edu}}} \and Sameer Iyer\thanks{Department of Mathematics, University of California, Davis, Davis, CA 95616, USA \href{mailto:sameer@math.ucdavis.edu}{\texttt{sameer@math.ucdavis.edu}}}
\and Linfeng Li\thanks{\footnotesize Department of Mathematics, University of California, Los Angeles, CA 90095, USA \href{mailto:lli265@math.ucla.edu}{\texttt{lli265@math.ucla.edu}}}
 \and Fei Wang\thanks{School of Mathematical Sciences, CMA-Shanghai, Shanghai Jiao Tong University, 
		 Shanghai, China \href{mailto:fwang256@sjtu.edu.cn}{\texttt{fwang256@sjtu.edu.cn}}}}
		% jacob jacob
\maketitle

\begin{abstract}
In this paper, we develop a stability threshold theorem for the 2D incompressible Navier-Stokes equations on the channel, supplemented with the no-slip boundary condition. The initial datum is close to the Couette flow in the following sense: the shear component of the perturbation is small, but independent of the viscosity $\nu$. On the other hand, the $x$-dependent fluctuation is assumed small in a viscosity-dependent sense, namely, $O(\nu^{\frac12}|\log \nu|^{-2})$. Under this setup, we prove nonlinear enhanced dissipation of the vorticity and a time-integrated inviscid damping for the velocity. These stabilizing phenomena guarantee that the Navier-Stokes solution stays close to an evolving shear flow for all time. The analytical challenge stems from a time-dependent nonlocal term that appears in the associated linearized Navier-Stokes equations.
%In this article, we prove a threshold theorem for the 2D, incompressible Navier-Stokes equations set on the channel, $\mathbb{T} \times [-1,1]$. The equations are supplemented with the classical no-slip boundary condition. Initial datum is assumed to be a perturbation of the Couette flow in the following sense: the shear component of the perturbation is assumed small, but independent of viscosity, in an appropriate Sobolev space. On the other hand, the $x$-dependent component of the perturbation is assumed small in a viscosity-dependent sense, $O(\nu^{\frac12}|\log \nu|^{-2})$. Under this setup, we prove classical enhanced dissipation of nonzero $x$-modes, and a time integrated inviscid damping for the velocity. 

\end{abstract}

\setcounter{tocdepth}{2}
{\small\tableofcontents}

\section{Introduction} \label{sec:intro}
In this paper, we study the quantitative stability threshold for perturbations of two-dimensional nearly-Couette shear flows with no-slip boundary conditions on the channel $\mathbb{T} \times [-1,1]$. Perturbations of the nonzero tangential modes that are assumed small (measured by a power of viscosity, $\nu$) in an appropriate Sobolev space, are shown to undergo enhanced dissipation and inviscid damping. 
%In general, the study of stability thresholds for fluid equations is by now vast, and we discuss related results in more detail below. We establish that the solution undergoes enhanced dissipation and inviscid damping that is in line with the behavior predicted by the linearized purely-Couette flow, studied in the same setting in \cite{}.  
%In the case of Navier-type boundary conditions, this corresponding result was proved in \cite{...}; 
%however, the case of no-slip boundary conditions presents greater difficulties.
More precisely, we consider the incompressible 2D Navier-Stokes equations
	\begin{align}
	\begin{split}
		&\partial_t \tilde u
	+ 	\tilde u\cdot\nabla \tilde  u
	-	
	\nu\Delta  \tilde u+\nabla \tilde{p}=0, 
	\\ & 
	\nabla \cdot \tilde u= 0,  \quad 
	\tilde{u}(t,x,{y=\pm 1}) =(\pm 1, 0),
	\\&
	\tilde u(t=0) = \tilde u_{\rm in} (x,y), 
	\end{split}
	\label{eq:NS}
	\end{align} 
on the domain $ \Omega \times [0,\infty)=(\mathbb{T}\times[-1, 1])\times [0,\infty)$, where $\nu>0$ denotes the inverse Reynolds number, $\tilde{u}(x,y)$ is the velocity field, and $\tilde{p}(x,y)$ is the pressure. We denote $U_{\rm in}(y)$ to be the shear-flow component of the initial data and $u_{\rm in} (x,y)$ to be the $x$-dependent component:  
\begin{align}
	& \tilde u_{\rm in}(x,y)
	=
	(\tilde u^{(1)}_{\rm in}(x,y),\tilde u^{(2)}_{\rm in}(x,y))
	=
	(U_{\rm in} (y), 0)
	+ 
	u_{\rm in}(x,y),
	\label{EQ127a}
\end{align}
where 
\begin{align}
	U_{\rm in}(y) :=
	\frac{1}{2\pi} 
	\int_{-\pi}^\pi 
	\tilde{u}^{(1)}_{\rm in}(x,y) dx, \qquad 
	u_{\rm in}(x,y) = \tilde{u}_{\rm in}(x,y) 
	-
	\frac{1}{2\pi}
	 \int_{-\pi}^\pi 
	 \tilde{u}_{\rm in}(x,y) dx.
\end{align}
We assume that the shear component, $U_{\rm in}(y)$, is an $O(\eps)$ Sobolev perturbation of the Couette flow (which satisfies the boundary conditions $U_{\rm in}(y=\pm 1)=\pm 1$), whereas the $x$-dependent component, $u_{\rm in}(x, y)$ is $O(\eps \nu^{\frac12} |\log \nu|^{-2})$. 
Importantly, $\eps>0$ is a small universal number and the viscosity is much smaller than this universal number, namely, $0 < \nu \ll \eps$.

In general, the study of stability thresholds for fluid equations is by now vast, and we discuss related results in more detail below. For now, the reader should keep in mind two ``nearest neighbor" results. First, Q. Chen, T. Li, D. Wei and Z. Zhang studied the nonlinear stability threshold of the purely-Couette case, $U_{\rm in}(y) = y$ % \myr{and $u_{\rm in} (x,y)=0$},
in \cite{CLWZ20}.  Our result can be seen as a generalization of \cite{CLWZ20} to the close-to-Couette shear flow case, i.e., $U_{\rm in}(y) = y + O(\eps)$. %\myr{and nonzero $u_{\rm in} (x,y)$}
Such a generalization is important to obtain from the point of view of potential treatment of much larger, smoother data, such as those in \cites{BHIWnonlinear,bedrossian2025landau}. Indeed, in the uniform-in-$\nu$ argument of \cite{BHIWnonlinear}, a similar threshold theorem was applied after the enhanced viscous time scale $t \gg O(\nu^{-1/3})$. 
While nonzero modes have been heavily dissipated by this time, the zero mode has not, and is instead an $O(\eps)$ perturbation of the Couette flow. 
In this case, one needs to understand precisely the stability threshold of nearly-Couette shear flows as opposed to purely-Couette flows, which is the setting of the present paper.  The second ``nearest-neighbor" result is that of \cite{BHIW23I}, which treats nearly-Couette shear flows with the Navier boundary condition. 
Therefore, one can also think of the current paper as the no-slip version of \cite{BHIW23I}, which presents significant new difficulties.  

For the case when $ u_{\rm in} (x,y)\equiv 0$, the solution of the Navier-Stokes equations \eqref{eq:NS} is given by the heat extension of the initial data, that is
\begin{align}
	%\tilde{u} (t,x,y)=
	(U(t,y),0) 
	 := (y+e^{\nu t \partial_{y}^2} (U_{\rm in} (y)- y) , 0), \label{def:barU}
\end{align}
where $e^{\nu t \partial_{y}^2}$ is the semi-group generated by the $1$D heat equation subject to homogeneous Dirichlet boundary conditions. 
Naturally, if $U_{\rm in} (y)$ is regular, then the resulting shear flow $U(t,y)$ is similarly uniform-in-time regular 
and converges to the Couette flow like $O(e^{-\nu t})$.  
This is because $U(t,y)-y$ solves the $1$D heat equation with the Dirichlet boundary condition and the initial data $U_{\rm in} (y) -y$, thus all the higher order derivatives of $U(t,y)-y$ are uniformly bounded in time.

For the case when $ u_{\rm in} (x,y)$ is small and nonzero, we write the solution $\tilde u (t,x,y)$ as a perturbation of the time-dependent shear profile $(U(t,y),0)$
and hence define
\begin{align}
	u(t,x,y):=
	\tilde u (t,x,y) 
	-
	(U(t,y), 0).
	\label{EQ00}
\end{align}
From  \eqref{eq:NS} and \eqref{EQ00}, we infer that the exact perturbation $u(t,x,y)= (u^{(1)}, u^{(2)})$ solves
	\begin{align}
	\begin{split}
		&\pa_tu +U \partial_x  u+  u\cdot \nabla  u -\nu\Delta  u + (U' u^{(2)}, 0) +\nabla p =0, 
%	\label{EQ01}
	\\&
	\nabla \cdot  u =0
	,\quad u (t,x,y=\pm 1)=0
%	\label{EQ02} 
	,
	\\&
	u(t=0,x,y)= u_{\rm in} (x,y),
	\end{split}
	\label{EQ01}
	\end{align}
where we used the notation $U'=\partial_y U(t,y)$.
Taking the curl of \eqref{EQ01} and denoting 
\begin{align}
	 \omega=\curl  u=-\pa_y u^{(1)}+ \pa_x u^{(2)},
\end{align} 
we obtain the vorticity formulation of the perturbed Navier-Stokes equations
\begin{align}
	\label{vort-orig}
	\begin{split}
	&\pa_t \omega +U\partial_x  \omega+  u\cdot \nabla  \omega =U''\partial_x\psi+\nu\Delta  \omega,
	\\&
	u =\nabla^{\perp}\psi,\quad	\Delta \psi = \omega, \quad (\pa_x,\pa_y)\psi\big|_{y=\pm 1}= 0, 
	\\&	
	 \omega(t=0,x,y)=  \omega_{\rm in}(x,y),
	\end{split}
\end{align} 
where $U''=\partial_{y}^2 U(t,y)$ and $\grad^\perp = (-\partial_y,\partial_x)$.

%\jacob{I don't like $\brak{f}$ as it conflicts with $(1+x^2)^{1/2}$; I'd rather just use $f_0 := \hat{f}_0$ or just drop the hats entirely from $ \omega_k$.}

The \emph{quantitative stability threshold} is, given a norm on the initial data $\norm{\cdot}_X$, the smallest $\gamma= \gamma (X) \geq 0$ such that
$\norm{u_{\rm in}}_{X} \ll \nu^\gamma$
implies `stability' and $\norm{u_{\rm in}}_{X} \gg \nu^\gamma$ in general permits `instability', or at least fully nonlinear dynamics.  
The exact quantification of `stability' may vary, but generally involves at least the enhanced dissipation characteristic of the linearized problem which shall be discussed further below. 
See \cite{BGM_Bull19} for a review and e.g. \cite{Gallay18,CZEW20,CLWZ20,BVW16,MW20,MW19} and the references therein for 2D, and for 3D see  \cite{WZ21,BGM15I,BGM15II,BGM15III} and the references therein.
When studying high Reynolds number stability problems, the mixing induced by the shearing enhances the stability of the equilibrium,
leading to (A) \emph{inviscid damping} wherein the perturbation velocity decays in the linearized (or nonlinear) Euler equations and uniformly in Reynolds number and (B) \emph{enhanced dissipation}, wherein the mixing accelerates the viscous dissipation, leading to rapid decay of the $x$-dependence on time-scales such as e.g. $\approx \nu^{-1/3}$ (rather than $\approx \nu^{-1}$). Many works have studied these in the linearized Navier-Stokes and Euler equations; see \cite{Zillinger16,BCZV17,Jia2020,WZZ19,WZZ20,WZZ18,LinXu2017,ISJ22,J20,RWWZ23} and the references therein for inviscid damping results and e.g. \cite{Gallay18,J23,CZEW20,BH20,WZZ20,CWZ23} and the references therein for results on enhanced dissipation (some results study both, such as \cite{J23,CWZ23}). See for example \cite{GCZ21,ABN22,BCZ15,BW13} and the references therein on mixing and enhanced dissipation by laminar flows in passive scalars. 
The purpose of this paper is to obtain the stability threshold for all shear flows close to Couette, i.e. the shear profile in \eqref{EQ127a} satisfies
$\Vert U_{\rm in} -y\Vert_{H^4}\ll 1$.

To state our result,
we introduce the $x$-Fourier transform
\begin{align}
	f_k(y):=
	\frac{1}{2\pi}\int_{-\pi}^\pi f(x,y) \exp(-ikx)dx,
	\quad \text{for~}k\in \mathbb{Z}.
\end{align} 
In particular,
\begin{align}
	f_0(y):= \frac{1}{2\pi}\int_{-\pi}^\pi f(x,y)dx.
\end{align}
%\begin{align}
%	\label{x-Four}
%	& 
%	f_k(y):=
%	\frac{1}{2\pi}\int_{-\pi}^\pi f(x,y) \exp(-ikx)dx,
%	\quad f(x,y)=
%	\sum_{k=-\infty}^\infty \wh f_k(y)\exp(ikx),
%	\\
%	\label{x-avre}&
%	f_0(y):=\wh f_0(y)=\frac{1}{2\pi}\int_{-\pi}^\pi f(x,y)dx,\quad f_\nq(x,y): =f(x,y)-{f}_0(y)=\sum_{k\neq 0}\wh f_k(y)\exp(ikx).
%\end{align} 
For $T>0$,
we define the energy
\begin{align}
	E(T, \omega)
	 &
%	 =E_0[T, \omega_0]+E_\nq[T, \omega_\nq] 
%	\\&:
	=\| \omega_0\|_{L^\infty_t(0,T;L_y^2)}
	+
	\sum_{k\neq 0}\big(\| e^{\aaa\nu^{1/3}t}
	(1-y^2)^{1/2}
	 \omega_k\|_{L_t^\infty(0,T; L^2_y)}
	\\&\quad
	+	
	\nu^{1/4}|k|^{1/2}\| e^{\aaa\nu^{1/3}t}
	 \omega_k\|_{L_t^2(0,T; L_y^2)}
	+
	|k|\| e^{\aaa \nu^{1/3}t}
	u_k\|_{L_t^2(0,T;L_y^2)}
	\big),
\end{align}
where $0<\zeta\ll 1$ is a small universal constant.

The following theorem is the main result of the paper. 
\begin{theorem}
	\label{T01}
Consider the solution $ \omega$ to \eqref{vort-orig} with regular initial data $ \omega_{\rm{in}}\in H^2 (\Omega)$. 
%\footnote{\siming{HS: I raise the regularity level from $H^{2}$ to $H^{3}$. In the proof of the nonlinear stability, the \eqref{Initial_pt} cannot be derived from the $H^2$ data because it requires $\ell^1$ in $k$. To reduce this index to $H^2$ (as in CLWZ), we need to keep track of the $k$-dependence of the $\| \omega_{\mathrm{in};k}\|_{H_k^2}$-term. }
%}
%Consider the Navier-Stokes equation \eqref{vort-orig} subject to the initial condition $ \omega_{\rm{in}}\in H^2 ( \omega)$. 
There exist constants $\nu_0,\kappa, \epsilon>0$ such that, for any $\nu \in (0,\nu_0)$, if the following conditions are satisfied \begin{enumerate}
\item Smallness: 
\begin{align}
	&\|U_{\rm in}-y\|_{H^4([-1,1])}\leq 	
	\kappa,
		\label{smallness1}
	\\ 
	&\| \omega_{\rm in}\|_{H^2 (\Omega)}\leq \ep\nu^{1/2} |\log \nu|^{-2};
	\label{smallness2}
\end{align}
%$\|U_{\rm in}-y\|_{H^4([-1,1])}\leq \ep$ and $\| \omega_{\rm in}\|_{H^2}\leq \ep\nu^{1/2} |\log \nu|^{-2}$;
\item Compatibility: 
\begin{align}
	&
	\int_{-1}^1 e^{\pm ky} \omega_{{\rm in};k}(y)dy=0,
	\quad \forall k\in \mathbb{Z}, \label{eq:AverageZero}\\&
\label{eq:bc:Uin}
	\partial_{y}^j U_{\rm in} |_{y=\pm 1}
	=0,\quad
	j=2,4,
\end{align}
\end{enumerate}
then the solution of \eqref{vort-orig} satisfies the following global-in-time estimate 
\begin{align}
	\sup_{T \geq 0} E(T, \omega)
	\leq 
	C\ep \nu^{1/2} |\log \nu|^{-1}.
\end{align}
\end{theorem}

Throughout this paper, we always assume that the smallness assumptions \eqref{smallness1}--\eqref{smallness2} and compatibility conditions \eqref{eq:AverageZero}--\eqref{eq:bc:Uin} hold without mention.
%\jacob{What do we get for $L^\infty_t L^\infty$ estimates on the velocity? } \myb{HS: We might lose there.}
%\linfeng{it seems that we didn't prove exactly what CLWZ proved and we didn't have $L^\infty L^\infty$ of the velocity.}
\begin{remark}
The compatibility condition \eqref{eq:AverageZero} is equivalent to the no-slip boundary condition on the initial velocity. 
\end{remark}\simhe{
\begin{remark}
Suitable adjustments of our argument yield the same stability threshold result for strictly monotone, spectrally stable shear flows, i.e., the shear flows discussed in \cite{CWZ23}.
\end{remark}}

In \cite{BVW16}, the analogue of this theorem was proved without boundaries (i.e. $y \in \mathbb R$) with initial data in Sobolev spaces $X = H^s$ with $s > 1$ (also without the logarithmic correction in $\nu$) and for Navier boundary conditions in the channel in \cite{BHIW23I}.  
For the case of no-slip boundary conditions and $\norm{U_{\mathrm {in}} - y} \lesssim \nu^{1/2}$, this theorem was first proved in \cite{CLWZ20}. It was essentially proved in \cite{BGMZ25} that this result is sharp up to the logarithmic correction in $\nu$ even for $\norm{U_{\mathrm in} - y} \lesssim \nu^{1/2}$ in any finite regularity due to instabilities in the boundary layer.

In the case of $\mathbb T \times \mathbb R$, \cite{MW20} establishes stability of the Couette flow under the assumption  $\norm{U_{\mathrm {in}} - y} \lesssim \nu^{1/2}$ in the almost-critical space, $H^{\log}_xL^2_y$, whereas \cite{LiMasmoudiZhao22b} proves a complimentary instability result. 
However, at higher Sobolev regularities, this was later improved to $\gamma \leq 1/3$ for sufficiently high Sobolev regularity \cite{MW19, WeiZhang23} and $\gamma\in[0,1/3]$ for suitable Gevrey classes \cite{LiMasmoudiZhao25}.  These higher regularity results are also expected to be sharp, due to the results of \cite{DM18}.
When one has at least Gevrey-2 regularity on the other hand, the result was proved with $\gamma = 0$ in \cite{BMV14} (i.e. the results are uniform-in-$\nu$). 
The corresponding uniform-in-$\nu$ results with Navier boundary conditions were proved in \cite{BHIWLinear,BHIWnonlinear}.

\section{Outline} \label{sec:outline}
%
%In this section, we introduce some key ideas and propositions. 
\subsection{Setup} 
The main difficulty of this paper is to obtain a sufficiently precise understanding of the \emph{linearized} Navier-Stokes equations associated with \eqref{vort-orig}.
The relevant linearized Navier-Stokes system for our analysis is the following system
\begin{align} 
	\label{lnr_Ns}
	\begin{aligned}
	&\pa_t \omega  +U\pa_x \omega -U''\pa_x\psi -\nu \de \omega =\na\cdot(\mathbbm{f}_1,\mathbbm{f}_2),\\
	&\displaystyle u=\na^\perp\psi,\quad\de \psi = \omega, \quad (\pa_x,\pa_y)\psi |_{y=\pm 1}=0,\\
	&  \omega(t=0,x,y)= \omega_{\mathrm{in}} (x,y).
    \end{aligned}
\end{align}
We emphasize that $U=U(t,y)$ is time-dependent, as defined in \eqref{def:barU}.
Taking the Fourier transform in $x$, the equations \eqref{lnr_Ns} reduce to the following $k$-by-$k$ system:\begin{subequations}\label{lnr_Ns_k}%($k\neq0$)
\begin{align} 
	\label{lnr_Ns_k_main}
	\begin{aligned}
	&\pa_t \omega_k  +ikU  \omega_k-U''ik\psi_k -\nu \de_k \omega_k =\na_k\cdot(\mathbbm{f}_{1;k},\mathbbm{f}_{2;k}),\\
	&\displaystyle u_k=\na_k^\perp\psi_k,\quad\de_k \psi_k = \omega_k, \\
	&ik\psi_k |_{y=\pm 1}=\pa_y\psi_k |_{y=\pm 1}=0,\\
	&  \omega_k(t=0,y)= \omega_{\mathrm{in};k} (y),
    \end{aligned}
\end{align}
where $\nabla_k^\perp :=(-\partial_y, ik)$, $\na_k:= (ik,\pa_y)$ and $\de_k:=\pa_y^2-|k|^2$.
% \siming{We note that for the $k=0$-mode, the stream function $\psi_0$ solves the Poisson equation $\pa_{yy}\psi_0=\omega_0$ subject to the Neumann boundary condition $\pa_y\psi_0\Big|_{y=\pm 1}=0$. This naturally enforced the consistency condition that $\int_{-1}^1 \omega_{0}(y)dy=0$. This is sufficient to guarantee that $(u^1_0,u^2_0)\Big|_{y=\pm 1}=0$. However, the boundary value for the stream function $\psi_0$ is not fixed.}
We note that the boundary condition for $\psi_k$ is over-determined. To make the system solvable, we always assume the consistency condition on $\omega_k$:
\begin{align}
\label{lnr_Ns_bc}
\int_{-1}^1 e^{\pm ky} \omega_{k}(t,y)dy=0,
	\quad \forall t\geq 0,\quad k\in \mathbb{Z}.
\end{align}\scalebox{0.001}{\eqref{lnr_Ns_k_main}\eqref{lnr_Ns_bc}}
\end{subequations}
One can view this as an integral boundary condition for $\omega_k. $

{Next we specify the functional framework. Since the vorticity undergoes a transient growth near the boundary, the desirable norm should balance out this growth. We recall the standard setup here, see, e.g., \cite{CLWZ20}. The main focus of the estimates on the linearized Navier-Stokes equations \eqref{lnr_Ns_k} is to develop the uniform-in-time estimate for the weighted $L^2$-norm of the vorticity:
\begin{align}
\|(1-|y|)^{1/2}\omega_k\|_{L_y^2},
\end{align}
and the time-integral estimate of the following quantity:
\begin{align}
	\label{Frz_Z_norm}
	\|\omega_k\|_{\mathcal{Z}}^2
	:=
	\nu^{1/3}|k|^{2/3}\|\rho_k^{1/2}\omega_k\|_{L_y^2}^2
	+
	|k|^2\|\na_k\de_k^{-1}\omega_k\|_{L_y^2}^2
	+
	\nu^{1/2}|k|\|\omega_k\|_{L_y^2}^2. 
\end{align}
Here, $\de_k^{-1}$ is the inverse of the Dirichlet Laplacian and the boundary weight $\rho_k$ is defined as \begin{align}\label{Frz_rhok}
	\rho_k:=\min\{L(1-|y|),1\},\quad L:=\nu^{-1/3}|k|^{1/3}.
\end{align} We observe that in the high-frequency regime, i.e., $\nu |k|^2\gtrsim 1$, the $\mathcal{Z}$-norm is equivalent to the following simplified expression:
\begin{align}\label{Z-norm_h_mode}
\|\omega_k\|_{\mathcal{Z}}^2\approx	|k|^2\|\na_k\de_k^{-1}\omega_k\|_{L_y^2}^2
	+
	\nu^{1/2}|k|\|\omega_k\|_{L_y^2}^2.
\end{align}}
%We introduce the $\mathcal{Z}$ norm for convenience, as the $L^2 \mathcal{Z}$ norm of $ w_k$ ultimately controls the left-hand side of \eqref{Lnr_low_est}.

The main linear estimate we establish for solutions to \eqref{lnr_Ns_k} is stated in the following proposition, which will be used to prove Theorem~\ref{T01}.

\begin{proposition}
	\label{pro:lin_st}
	Let $ \omega_k$ be the solution of $\eqref{lnr_Ns_k}_{k\neq 0}$ and $U(t,y)$ be the heat extension defined in \eqref{def:barU}. 
	There exist constants $\zeta>0$ and $0<\nu_0,\kappa\ll 1$ such that
	if $\|U_{\rm in}-y\|_{H^4}\leq \kappa$ and $\nu \in (0,\nu_0]$, then
\begin{align}\label{Lnr_st_est}
	&\|
	e^{\zeta\nu^{1/3}t}(1-|y|)^{1/2} \omega_k\|_{L^\infty_t L_y^2}^2
	+ \|e^{\zeta\nu^{1/3}t} \omega_k\|_{L^2_t \mathcal{Z}}^2
	%+	|k|^2\|e^{\zeta\nu^{1/3}t}u_k\|_{L^2_t L_y^2}^2	+	\nu^{1/2}|k|\|e^{\zeta\nu^{1/3}t} \omega_k\|_{L^2_t L_y^2}^2
	\\&
	\lesssim 
	|\log \nu|^2
	\lf(
	|k|^{-2}	\| \omega_{{\rm in};k}\|_{H_k^{2}}^2 
%	\nu^{1/3}|k|^{-4/3}\|\pa_y \omega_{{\rm in};k}\|_{L^2}^2
%	+
%	\| \omega_{{\rm in};k}\|_{H_k^{2}}^2 
	+
	\min\{\nu^{-1/3}|k|^{4/3},\nu^{-1}\}
	\|e^{\zeta\nu^{1/3}t}\mathbbm{f}_{1;k}\|_{L^2_t L_y^2}^2
	+
	\nu^{-1}\|e^{\zeta\nu^{1/3}t}\mathbbm{f}_{2;k}\|_{L^2_t L_y^2}^2\rg).\nonumber
\end{align}
Here, the $\mathcal{Z}$-norm is defined in \eqref{Frz_Z_norm} and $\|f_k\|_{H_k^2}^2:=\sum_{\al+\beta=2}\||k|^\al \pa_y^\beta f_k\|_{L_y^2}^2$. 
%\myb{HS: Here, at the high mode, we use $\nu^{-1}$ instead of the $k$ dependent bound. This comes from the high mode estimate Proposition 2.3. }
\end{proposition}
Throughout this paper, we always assume that $k\neq 0$, unless stated otherwise.
\begin{remark}
The same proof also yields the space-time estimate \eqref{Lnr_st_est} for the system \eqref{lnr_Ns_k} with general time-dependent shear flow $\mathbb{U}(t,y)$, provided the flow satisfies the following conditions:
\begin{align}
	&\|\mathbb{U}-y\|_{L_t^\infty H_y^4}\leq \kappa,\quad  (\mathbb{U}-y)\big|_{y=\pm 1}=0,
	\\&
	\|\pa_y^2(\mathbb{U}(s,\cdot)-\mathbb{U}(t,\cdot))\|_{L_{y}^\infty}+\sup_{y\in [-1,1]}
	\Big|
	\frac{\mathbb{U}(s,y) -\mathbb{U}(t,y)}{1-|y|}
	\Big|\leq C\nu|s-t|,\quad \forall s, t\geq 0.
	\label{EQ140a}
\end{align}
\end{remark}

In proving Proposition~\ref{pro:lin_st}, we treat the high and low frequency regimes separately.

\noindent
{\bf a) High Frequency Regime:} 
We establish the following proposition for the high frequency regime. 
\begin{proposition}\label{pro:lin_st_hi}
Let $ \omega_k$ be the solution of \eqref{lnr_Ns_k} and $U(t,y)$ be the heat extension defined in \eqref{def:barU}.  
There exists a constant $\zeta>0$
such that for the high frequency regime $
	\nu|k|^2\ge \frac{1}{2}
	\|\pa_y U\|_{L_{t,y}^\infty},$ 
the following estimate holds:
	\begin{align} 
			&
	\|e^{\zeta\nu^{1/3}t}(1-|y|)^{1/2} \omega_k\|_{L_t^\infty L_y^2}^2
	+\|e^{\zeta\nu^{1/3}t} \omega_k\|_{L_t^2\mathcal Z}^2
	%$%+|k|^2\|e^{\zeta\nu^{1/3}t}u_k\|_{L_t^2L_y^2}^2+\nu^{1/2}|k|\|e^{\zeta\nu^{1/3}t} \omega_k\|_{L_t^2L_y^2}^2
	%\\&\quad
	\les
	\nu^{-1}
	(	
	\Vert
	e^{\aaa \nu^{1/3} t} \mathbbm{f}_{1;k} \Vert_{L_t^2 L_y^2}^2 
	+  
	\Vert
	e^{\aaa \nu^{1/3} t}
	\mathbbm{f}_{2;k}
	\Vert_{L_t^2 L_y^2}^2
	)
	+ 
	\norm{ \omega_{{\rm in}, k}}_{L_y^2}^2. 
\end{align}
\end{proposition}
The proof of this proposition is given in Section \ref{sec:high} and follows from standard energy estimates. Therefore, we omit further discussion and focus on the challenging low frequency regime.  
% \siming{I have adjust the proof in the previous high mode proposition and change all the $\U$ in the proof to $U$. The argument after \eqref{EQ127c} is similar so we can omit. } 

\noindent
{\bf b) Low Frequency Regime:}
To develop the linear estimates in the presence of the no-slip boundary condition \eqref{vort-orig}$_2$, we must use resolvent estimates similar to those studied, for example, in \cites{CLWZ20,CWZ23,ChenWeiZhang20} among many other papers.
Since the low frequency regime is more difficult, we will restrict our discussion to this case for the remainder of the section.
Our goal is to establish the following proposition.
\begin{proposition}
\label{pro:lin_low} 
Let $ \omega_k$ be the solution of \eqref{lnr_Ns_k} and $U(t,y)$ be the heat extension defined in \eqref{def:barU}. 
There exists a constant $\zeta>0$ such for the low frequency regime
$
\nu|k|^2\leq 
\inf_{t\in\rr_+}\|\pa_y U(t,\cdot)
\|_{L_{y}^\infty}, 
%\myr{\Rightarrow \nu|k|^2\leq \|\pa_y U\|_{L_{t,y}^\infty}}
$
the following estimate holds:
\begin{align}
\begin{split}
&\|e^{\zeta\nu^{1/3}t}(1-|y|)^{1/2} \omega_k\|_{L_t^\infty L_y^2}^2
+\|e^{\zeta\nu^{1/3}t} \omega_k\|_{L^2\mathcal{Z}}^2
%$%+|k|^2\|e^{\zeta\nu^{1/3}t}u_k\|_{L^2L^2}^2+\nu^{1/2}|k|\|e^{\zeta\nu^{1/3}t} \omega_k\|_{L^2L^2}^2
\\&\quad
\lesssim 
|\log \nu|^2
\lf(
%\nu^{1/3}|k|^{-4/3}
%\|\pa_y \omega_{{\rm in};k}\|_{L^2}^2
%+
|k|^{-2}
\| \omega_{{\rm in};k}\|_{H_k^{2}}^2\rg. 
%\\&\qquad
\lf.+\nu^{-1/3}|k|^{4/3}\|e^{\zeta\nu^{1/3}t}\mathbbm{f}_{1;k}\|_{L_t^2L_y^2}^2+
\nu^{-1}\|e^{\zeta\nu^{1/3}t}\mathbbm{f}_{2;k}\|_{L_t^2L_y^2}^2\rg).
\end{split}
\label{Lnr_low_est}
\end{align}
\end{proposition} 

We will sketch the proof of this proposition in the forthcoming subsections. The complete proof will be concluded in Section~\ref{sec:frozen}. Furthermore, we observe that combining Propositions \ref{pro:lin_st_hi} and~\ref{pro:lin_low} yields the full linear estimates in Proposition \ref{pro:lin_st}. This is because when the shear flow is sufficiently close to Couette flow, the high frequency regime $
	\nu|k|^2\ge \frac{1}{2}
	\|\pa_y U\|_{L_{t,y}^\infty}$ and the low frequency regime  $
\nu|k|^2\leq 
\inf_{t\in\rr_+}\|\pa_y U(t,\cdot)
\|_{L_{y}^\infty}$ cover all possible $k\neq 0$.

The main challenges in deriving reasonably sharp quantitative estimates for \eqref{lnr_Ns} in the limit $\nu \to 0$ are the following:
\begin{enumerate}[label=\Alph*.]
	\item The time-dependency of the background shear flow $U(t,y)$; 
	\item The nonlocal effect encoded by the term $U''\partial_x\psi$ in \eqref{vort-orig}$_1$;
	\item The boundary layer vorticity generation induced by the no-slip boundary condition \eqref{vort-orig}$_2$.
\end{enumerate} 
We elaborate on these three aspects here. 

\vspace{2 mm}

\noindent
{\bf Challenge A: }The first challenge arises from the fact that representing semigroups in terms of resolvents uses a Laplace transform in time, which is not applicable to time-dependent linear operators. 
To overcome this difficulty, we employ the frozen-time scheme introduced in \cite{ChenWeiZhang20}, wherein  the analysis is carried out for a family of Navier-Stokes equations linearized around the stationary shear flows $U(t=j\nu^{-1/3},y)$, where $j\in\mathbb{N}_0$.
After establishing suitable estimates for this family, a layer-cake construction (see Section~\ref{sec:Frozen} for details) enables us to reconstruct the solution of the original system \eqref{lnr_Ns_k} with a time-dependent shear flow. Consequently, the frozen-time reduction allows us to focus on the time-independent linearized Navier-Stokes system.

\vspace{2 mm}

\noindent
{\bf Challenge B:} The second difficulty arises from the linear nonlocal term $U''\pa_x\psi$. Although the closeness of the shear flow $U=U(t,y)$ to the Couette profile $(y,0)$ ensures $U''=O(\epsilon)$, this linear nonlocal term cannot be controlled directly by viscosity. 
A suitable inviscid damping mechanism is therefore necessary to absorb its effect. 
Under the Navier/Lions boundary condition $ \omega|_{y=\pm 1}=0$, one may use a hypocoercivity energy method (see, e.g., \cites{villani2009, BCZ17}), combined with the singular integral operator $\mathfrak{J}_k$ introduced in \cite{BHIW23I}, to derive the necessary estimates as carried out in \cite{BHIW23I}; in that setting, Challenge A is also avoided. 
In contrast, under the no-slip boundary condition, such energy-type estimates break down due to boundary contributions. 
In this work, we instead use the singular integral operator $\mathfrak{J}_k$ to provide the inviscid damping estimates on the resolvent which is one of the most crucial parts of the resolvent bounds (see Proposition \ref{pro_Nav}). The introduction of this operator also simplifies the compactness argument presented in \cite{CWZ23} in the close-to-Couette setting.
 
 \vspace{2 mm}
 
\noindent
{\bf Challenge C:} One of the most delicate aspects of the argument is the analysis of the boundary layers.
We decompose the solution of \eqref{lnr_Ns} into two parts: a component $ \omega^I$ which has zero initial data and accounts for the forcing on the right hand side, and a component $ \omega^H$, which incorporates the initial condition but is not subject to any external forcing.
The latter can be treated relatively easily if one has a good enough method for the former component $ \omega^I$. 
As established in \cites{CLWZ20,CWZ23}, it is convenient to further decompose $ \omega^I$ into $ \omega_{\rm Nav}$, which satisfies the Navier/Lions boundary condition $ \omega_{\rm Nav}|_{y=\pm1}=0$, and a boundary corrector $ \omega_{\rm B}$.  The corrector $ \omega_{\rm B}$ is determined by the velocity $u^{(1)}_{\rm Nav}$ generated by $ \omega_{\rm Nav}$ at the boundary points $y=\pm1$, and it will generally have a large $L^2$ norm but be very well localized near the boundary. 
The estimates in \cite{CWZ23} are insufficient to obtain a sharp enough bound on $ \omega_{\rm B}$. In particular, one cannot directly derive the estimate \eqref{cest_FL2}, which is essential for controlling $ \omega_{\rm B}$, from Proposition \ref{pro_Nav}. We therefore establish a new estimate in Lemma \ref{L02} to address this issue.

The remaining part of the section is organized as follows: In Section \ref{sec:tim_ind_outln}, we present the linear estimates for the time-independent shear flows, which serve as the building blocks for the frozen-time scheme. In Section \ref{sec:Frozen}, we present the frozen-time scheme and upgrade the linear estimates for the time-independent shear flows to the time-dependent case and conclude the proof of the key Proposition \ref{pro:lin_low}. Finally, in Section \ref{sec:nlnr_outln}, we provide the sketch of the proof of Theorem \ref{T01}.

\subsection{Linear Estimate: Time-independent Shear}
\label{sec:tim_ind_outln}
In this section, we consider the following frozen-time linearized Navier-Stokes equations
\begin{align}
	\begin{aligned}
	&\pa_t  \omega_{k} -\nu\de_k \omega_k+\mathcal{U}ik  \omega_{k}-\mathcal{U}''ik \psi_{k}
	=
	\nabla_k \cdot
	(\mathbbm{f}_{1} ,\mathbbm{f}_{2}),
	\\&
	\de_k\psi_{k}= \omega_{k},\quad(ik,\pa_y)\psi_{k}\big|_{y=\pm1}=0,\quad \int_{-1}^1 e^{\pm k y} \omega_{k} dy=0,
	\\& \omega_{k}(t=0,y)= \omega_{{\rm in};k}(y),
	\end{aligned}
\label{lnrzd_NS_1}
\end{align} 
where $\nabla_k= (ik, \partial_y)$ and $\mathcal{U}= \mathcal{U} (y)$ is a time-independent shear profile. 
We assume the boundary condition
\begin{align}
	(\U -y) |_{y=\pm 1}=0.
	\quad
	\label{EQcomp}
\end{align}
Note that the system \eqref{lnrzd_NS_1} is a general form of the frozen-time systems \eqref{Frz_om_0} and \eqref{Frz_om_j}, and its estimates apply to both systems. We focus on the following low frequency regime:
\begin{align}
\nu|k|^2\leq \|\U'\|_{L^\infty_y}. \label{low_reg}
\end{align}
%Here the constant ``$2$" is not essential. One can replace that by any universal constant. (HS: Check!  We extend the range of the low frequency regime, because we want to avoid the appearance of $k$ which is $\|\pa_y U\|_{L^\infty_{t,y}}> |k|> \|\pa_y U(j\nu^{-1/3},\cdot)\|_{L^\infty_y}$ for some $j$. That will violate both the high and low frequency condition (if we just have $\nu k^2\leq \|\U'\|_{L^\infty_y}$). 
\simhe{One can hope to extend the current result to the spectrally stable shear case. However, one has to be careful about this frequency cutoff. We might need to introduce the forcing to make sure that $U(t,y)\equiv \U(y)$ for all time. Then the definition of the high-$k$ and low-$k$ will not be ambiguous. The low frequency resolvent estimate is only proved for $\nu k^2\leq \|\U'\|_{L^\infty_y}$ in their setting. Even though I guess it is still true for $\nu k^2\leq 	C\|\U'\|_{L^\infty},$ I am not 100\% sure. }
Our goal is to establish the following key proposition.
\begin{proposition}
	\label{pro:lin_st:calU}
	Let $ \omega_k$ be the solution of \eqref{lnrzd_NS_1} in the low frequency regime \eqref{low_reg}.
	There exist constants $\delta>0$ and $0<\nu_0,\kappa\ll 1$ such that
	if $\|\mathcal{U}-y\|_{H^4}\leq \kappa$, the boundary condition \eqref{EQcomp} holds, and $\nu \in (0,\nu_0]$, then
	\begin{align}
		\begin{split}
			&\|e^{\delta\nu^{1/3}t}(1-|y|)^{1/2} \omega_k\|_{L^\infty_t L_y^2}^2
			+
			\Vert e^{\delta \nu^{1/3} t} \omega_k\Vert_{L^2_t \mathcal{Z}}^2
%			+
%			\nu^{1/3}|k|^{2/3}
%			\|e^{\delta\nu^{1/3}t}\rho_k^{1/2} \omega_k\|_{L^2_t L_y^2}^2
%			\\&\qquad
%			+
%			|k|^2\|e^{\delta\nu^{1/3}t}u_k\|_{L^2_t L_y^2}^2
%			+
%			\nu^{1/2}|k|\|e^{\delta\nu^{1/3}t} \omega_k\|_{L^2_t L_y^2}^2
			\\&\quad
			\lesssim 
			|\log \nu|^2\lf(
			|k|^{-2}
			\| \omega_{{\rm in};k}\|_{H_k^{2}}^2 
			%	\rg. 
			%%	\\&\qquad
			%	\lf.
			+
			\nu^{-1/3}|k|^{4/3}\|e^{\delta\nu^{1/3}t}\mathbbm{f}_1\|_{L^2_t L_y^2}^2
			+
			\nu^{-1}\|e^{\delta\nu^{1/3}t}\mathbbm{f}_2\|_{L^2_t L_y^2}^2\rg).
		\end{split}\label{Lnr_st_est_low}
	\end{align}
%\begin{align}
%\begin{split}
%	&\|e^{\delta\nu^{1/3}t}(1-|y|)^{1/2} \omega_k\|_{L^\infty_t L_y^2}^2
%	+
%	\nu^{1/3}|k|^{2/3}
%	\|e^{\delta\nu^{1/3}t}\rho_k^{1/2} \omega_k\|_{L^2_t L_y^2}^2
%		\\&\qquad
%		+
%	|k|^2\|e^{\delta\nu^{1/3}t}u_k\|_{L^2_t L_y^2}^2
%	+
%	\nu^{1/2}|k|\|e^{\delta\nu^{1/3}t} \omega_k\|_{L^2_t L_y^2}^2
%	\\&\quad
%	\lesssim 
%	|\log \nu|^2\lf(
%	|k|^{-2}
%	\| \omega_{{\rm in};k}\|_{H_k^{2}}^2 
%%	\rg. 
%%%	\\&\qquad
%%	\lf.
%	+
%	\nu^{-1/3}|k|^{4/3}\|e^{\delta\nu^{1/3}t}\mathbbm{f}_1\|_{L^2_t L_y^2}^2
%	+
%	\nu^{-1}\|e^{\delta\nu^{1/3}t}\mathbbm{f}_2\|_{L^2_t L_y^2}^2\rg).
%	\end{split}\label{Lnr_st_est_low}
%\end{align}
\end{proposition}
We present the proof of Proposition~\ref{pro:lin_st:calU}
in Section \ref{sec:space-time}.  
%To prove it, we consider two frequency regimes, the high-frequency case where
%\begin{align}
%	\nu|k|^2> \|\U'\|_{L^\infty_y},
%\end{align} 
%and low-frequency case, where
%\begin{align}
%	\nu|k|^2\leq \|\U'\|_{L^\infty_y}.
%\end{align} 
%
%
%\begin{itemize}
%\item High-frequency regime: 
%\begin{align}
%\nu|k|^2> \|\U'\|_{L^\infty_y}.
%\end{align}
%\item Low-frequency regime: 
%\begin{align}
%\nu|k|^2\leq \|\U'\|_{L^\infty_y}.
%\end{align}
%\end{itemize}
%In the high-frequency regime, the dissipation term $\nu|k|^2 w_k$ is sufficiently strong, so that a standard energy estimate yields the desired result (see Proposition \ref{pro:lin_st_hi}). 
%We therefore focus on the low-frequency regime in the remaining discussion. 
To prove it, we first decompose the solution of \eqref{lnrzd_NS_1} into its homogeneous component $ \omega_k^H$ and inhomogeneous component $ \omega_k^I$, namely,
\begin{align} \label{deco:1}
 \omega_k =  \omega_k^H +  \omega_k^I, 
\end{align}
where $ \omega_k^I$ solves
\begin{align}\label{om_I}
	\left\{
	\begin{aligned}
	&\pa_t  \omega^I_{k}
	+
	ik
	\mathcal{U}(y) \omega^I_{k}
	-
	ik
	\mathcal{U}''(y)
	\psi^I_{k}
	-\nu \de_k \omega^I_{k}
	=
	\mathbb{F}_k,
	\\&
	\de_k\psi_k^I= \omega_k^I,\quad\psi_k^I\big|_{y=\pm 1}=0,\quad\int_{-1}^1 e^{\pm k y} \omega^I_{k}(t,y)dy= 0, \\
	& \omega^I_{k}(t=0,y)=0,
	\end{aligned}
	\right.
\end{align}
and $ \omega_k^H$ solves 
\begin{align}
\left\{
	\begin{aligned}\label{om_H}   &\pa_t  \omega_{k}^H
		+
		ik
		\mathcal{U}(y)  \omega^H_{k}
		-ik
		\mathcal{U}''(y)\psi^H_{k}
		-
		\nu\de_k \omega^H_{k}=0,
		\\&
		\de_k\psi_k^H= \omega_k^H,\quad\psi_k^H\big|_{y=\pm 1}=0, \quad\int_{-1}^1 e^{\pm k y} \omega^H_{k}(t,y)dy= 0, \\
     & \omega^H_{k}(t=0,y)= \omega_{{\rm in};k}(y).
 \end{aligned}
 \right.
\end{align}
The inhomogeneous component $ \omega_k^I$ aims to address the external forcing, while the homogeneous component captures the linear evolution of the initial data. 
Importantly, obtaining suitable estimates for the inhomogeneous problem also provides control over the homogeneous part. 
We begin with the analysis of $ \omega_k^I$ in Section \ref{sec:tmnd_I}, followed by that of $ \omega_k^H$ in Section \ref{sec:tmnd_H}.

\subsubsection{Decomposition of the Inhomogeneous Solution $ \omega_k^I$}
\label{sec:tmnd_I}

%The main strategy to derive the estimate for the $ w_k^I$ is to consider the resolvent estimates. 
As in \cite{CLWZ20,CWZ23,BH20}, we use resolvent estimates. 
To motivate the key equations under consideration, we apply the Fourier transform in time to a suitably normalized solution.
% (here $\wt \lambda\in \rr$ is the Fourier variable with respect to $t$, and we define $\lambda_r=-\wt \lambda/k$)
% \begin{align}
%	&w^I_k^\circ(\lambda_r,y)=\int_0^\infty  w^I_{k}(t,y) e^{-it \wt \lambda}dt=\int_0^\infty  w^I_{k}(t,y) e^{it  k\lambda_r}dt;\\
%	&\psi^I_k^\circ(\lambda_r,y)=\int_0^\infty \phi^I_{k}(t,y) e^{-it \wt \lambda}dt=\int_0^\infty \phi^I_{k}(t,y) e^{it  k\lambda_r}dt;\\
%	&F_j^\circ(\lambda_r,k,y)=\int_0^\infty f_j(t,k,y)e^{itk\lambda_r }dt,\quad j=1,2.
%\end{align}
%These are real Fourier transforms of the solution in the $t$ variable. 
Since we expect the solutions to decay with rate $e^{-\zeta \nu^{1/3} k^{2/3}t}$ for some constant $\zeta>0$, we consider the normalized vorticity, stream function, and forcing with the enhanced dissipation time weight, i.e.,
\begin{align}
	 \omega_{k}^I(t,y)e^{\zeta\nu^{1/3}|k|^{2/3}t},\quad \psi_k^{I}(t,y)e^{\zeta\nu^{1/3} |k|^{2/3}t},\quad \mathbb{F}_{k}(t,y)e^{\zeta\nu^{1/3}|k|^{2/3}t}
	. 
\end{align}
%Without loss of generality, we set $k> 0$ with straightforward adjustments for negative $k$.  
We define the intermediate Fourier variable which corresponds to $t$, i.e., $\wt \lambda\in \rr$. 
To be consistent with the literature, we also recall the actual spectral parameter 
\begin{align}
	\lambda_r:=-\wt \lambda/k,\quad
	\lambda_i:=-\zeta(\nu k^{-1})^{1/3},\quad \lambda:=\lambda_r+i\lambda_i.
\end{align}
This is equivalent to shifting the contour in the complex plane accordingly. 
Furthermore, we extend the aforementioned functions to be zero for the time regime $t<0. $ This yields a transform which is essentially equivalent to the Laplace transform (though less flexible). The corresponding Fourier transforms hence read as follows
\begin{align}
	\label{phys_to_Fourier}	 
	\begin{split}
	w^I_k(\lambda,y) :=&\int_0^\infty \omega_{k}^I(t,y)e^{it  k \lambda}dt, \quad	\phi^I_k(\lambda,y):= \int_0^\infty \psi_{k}^I(t,y)e^{it  k\lambda}dt, 
	\\
	F_k(\lambda,y):=&\int_0^\infty \mathbb{F}_k (t,y) e^{it  k\lambda}dt.
	\end{split}
\end{align}
Note that $ w^I_{k}(t=0,y)=0$, and thus the following relation holds
\begin{align}
	\int_0^\infty\pa_t \omega_k^{I} e^{itk\lambda}d t=- \omega_{k}^I(t=0,y)-ik\lambda\int_0^\infty  \omega_{k}^I e^{itk\lambda}d t=-ik\lambda w^I_k.
\end{align}
Consequently, after the Fourier transform in the $t$-variable, the system \eqref{om_I} becomes
\begin{align}
	\label{eq:w_I}
	\left\{
	\begin{aligned}
	&-\nu\Delta_k
	w^I_k
	+
	ik\lf(\mathcal{U}-\lambda\rg)w^I_k-ik\;\mathcal{U}''\phi^I_k	
	=	F_k,
	\\&
	\Delta_{k}\phi^I_k=
	w^I_k,\quad (ik,\pa_y)\phi^I_k\Big|_{y=\pm1}=0,\quad 
	%	\psi^I_k\Big|_{y=\pm 1}=0,
	%	\\& 
	\int_{-1} ^1 e^{\pm ky}w^I_k dy
	=0.
	\end{aligned}
	\right.
\end{align}
%Here, we can clearly see that we needed $\mathcal{U}$ to be time independent for this.   
%\textcolor{red}{ Sometimes in the analysis, one might choose $\zeta=0$. In this case, the $\lambda=\lambda_r.$ Sometimes, you will see this in the Chen-Li-Wei-Zhang paper. However, if we get estimate which are independent of $\zeta$, we can apply them in both the $\zeta=0$ and $\zeta>0$ case. }
%\jacob{Do we actually do this anywhere in this paper?}
The above boundary condition on $w^I_k$ is more difficult than the Navier boundary condition case. 
A classical idea used in \cite{CLWZ20} is to decompose the solution into a contribution associated with the Navier boundary conditions and then additional boundary correctors: 
\begin{align}\label{w_I_decmp}
	 w^I_k(\lambda,y)
	=
	 w_{{\rm Na}; k}(\lambda,y)
	+
	\underbrace{c_+(\lambda; k) w_{+; k}(\lambda,y)%$%}_{=:  w_{\mathrm{BC}, +; k}(\lambda, y)}
	+
	c_-(\lambda; k)  w_{-; k}(\lambda,y)}_{w_{\mathrm{BC};k}}%\underbrace{}_{=:   w_{\mathrm{BC}, -; k}(\lambda, y)}
	,
\end{align}
where $w_{{\rm Na}; k}$ is the solution of the system with Navier boundary condition and $w_{\mathrm{BC}; k}$ is the boundary corrector that ensure that $w^I_k$ satisfies the no-slip boundary condition. More precisely, the Navier component $w_{{\rm Na};k}$ solves
\begin{align} 
\left\{
\begin{aligned}
	\label{Nav_eqn}
	&    -\nu(\pa_y^2-k^2) w_{{\rm Na}; k}
	+
	ik\lf(\mathcal{U}-\lambda\rg) w_{{\rm Na}; k}
	-
	ik\mathcal{U}''\phi_{{\rm Na}; k}	
	=
	F_k,
	%\mathbb {F}_k,
	\\&
	\Delta_{k}\phi_{{\rm Na}; k}
	= 
	 w_{{\rm Na}; k},\quad \phi_{{\rm Na}; k}\big|_{y=\pm1}=0,
%	\\&
%	\psi_{{\rm Na}; k}\big|_{y=\pm 1}=0,
	\quad
	 w_{{\rm Na}; k}\big|_{y=\pm 1}=0.
\end{aligned}
\right.
\end{align}
Note that the boundary conditions \eqref{Nav_eqn}$_2$ imply that $
\pa_y^2	 \phi_{{\rm Na}; k} \big|_{y=\pm 1}=0$. 
%The boundary corrector $w_{{\rm BC},\pm; k}(\lambda, y)$ in \eqref{w_I_decmp} has two components: the elementary parts $w_{\pm ;k}$ and the coefficients $c_\pm$. 
The elementary parts $w_{\pm;k}$ in the boundary corrector solve 
\begin{align} 
	\label{bdy_crt}
	\left\{
	\begin{aligned}
		&    -\nu(\pa_y^2-k^2) w_{\pm;k}
		+
		ik\lf(\mathcal{U}-\lambda\rg) w_{\pm;k}
		-ik
		\U''\phi_{\pm;k}	=0,
		\\&
		\phi_{\pm;k}=\Delta_{k}^{-1}
		 w_{\pm;k},
	\quad
		\pa_y\phi_{\pm ;k}\big|_{y=\pm 1}=1, \quad 
		\pa_y\phi_{\pm ;k}\big|_{y=\mp 1}=0.  
	\end{aligned}
	\right.
\end{align}Here, $\de_k^{-1}$ is the inverse of the Dirichlet Laplacian. 
The coefficients $c_{\pm}(\lambda; k)$ are chosen to enforce the boundary condition
\begin{align}
	\int_{-1}^{1} e^{\pm ky}  w^I_k\,dy=0
\end{align}  
in \eqref{eq:w_I}$_2$. Therefore, we take
\begin{align}\label{c_pm_frml}
	c_{+} (\lambda;k)=-\pa_y\phi_{{\rm Na};k}(y=1),
	\quad	c_{-} (\lambda;k)=-\pa_y\phi_{{\rm Na};k}(y= -1).
\end{align}
The boundary conditions of $\phi_{\pm;k}$ are chosen so that the elementary boundary corrector corrects one boundary condition while leaving the other intact. 
The main difficulty is that $\U (y)$ is not the Couette flow, leading to the presence of the nonlocal term $\U'' ik \phi_{\pm;k}$. 
A related problem was addressed in \cite{CWZ23}, where numerous estimates for similar boundary correctors were established.
%\jacob{But there's still a step we do differently}

The goal of the decomposition is to obtain the following estimates on the inhomogeneous component.
\begin{proposition}\label{pro:w_I}
%Suppose the $\mathcal{U}$ is close to the Couette flow in $C^3$. 
%Assume that $\nu |k|^2\leq \|\U'\|_{L^\infty}$. 
Let $ \omega_k^I$ be the solution of \eqref{om_I} in the low frequency regime \eqref{low_reg} with the forcing 
\begin{align}
	\mathbb{F}_k=ik\mathbbm{f}_1+\pa_y \mathbbm{f}_2+\mathbbm{f}_3,\quad \mathbbm{f}_1,\mathbbm{f}_2\in L^2,\,  \mathbbm{f}_3\in H_0^1.
\end{align} 
There exist constants $\delta>0$ and $0<\kappa,\nu_0\ll 1$ such that if
$\norm{\mathcal{U} - y}_{H^4} < \kappa$ and $\nu \in (0,\nu_0)$
then the following estimate holds:
\begin{align}
	\label{w_ I_est}
	\begin{split}
	&	
	\nu^{1/2}|k|\|e^{\delta(\nu k^2)^{1/3}t}
	\omega_k^I\|_{L^\infty L^2}^2
	+
	\Vert e^{\delta (\nu k^2)^{1/3} t} \omega^I_k\Vert_{L^2 \mathcal{Z}}^2
%	+\nu^{1/3}|k|^{2/3}
%	\|e^{\delta(\nu k^2)^{1/3} t} \rho_k^{1/2} \omega_k^I\|_{L^2L^2}^2
%	+
%	|k|^2\|e^{\delta(\nu k^2)^{1/3}t} u_k^I\|_{L^2L^2}^2
%		\\&\qquad
%		+
%	\nu^{1/2}|k|\| e^{\delta(\nu k^2)^{1/3}t}
%	 \omega_k^I\|_{L^2L^2}^2
	\\&
	\lesssim
	\log^2 L
	\big( \nu^{-1/3}|k|^{4/3}
	\|e^{\delta(\nu k^2)^{1/3}t}\mathbbm{f}_1\|_{L^2L^2}^2
	+
	\nu^{-1}\|e^{\delta(\nu k^2)^{1/3}t}\mathbbm{f}_2\|_{L^2L^2}^2
	%$%\\&\qquad
	+
	\|e^{\delta(\nu k^2)^{1/3}t}
	\na_k\mathbbm{f}_3\|_{L^2L^2}^2
	\big).
\end{split}
\end{align}
%\myb{HS: I have changed $|k|^{-1} 
%\log^2 L\|\na_k\mathbbm{f}_3\|_{L^2L^2}^2$ to $
%\log^2 L\|\na_k\mathbbm{f}_3\|_{L^2L^2}^2$.} 
\end{proposition}
%textcolor{red}{[We should remark on how to prove this, especially on why $\mathfrak{J}$ is useful?]} \myb{HS: Addressing this comment below.}
The proof of the proposition requires several preparatory steps and will be completed in Section \ref{sec:space-time}. By Plancherel's theorem and the relation \eqref{phys_to_Fourier}, the problem reduces to establishing estimates for each component of $w^I_k$ in \eqref{w_I_decmp}, namely, the Navier component $w_{{\rm Na};k}$, the unit boundary correctors $w_{\pm;k}$, and the coefficients $c_{\pm;k}$. 
	Compared with \cite[Proposition~6.5]{CLWZ20}, our analysis includes the additional treatment of the
	component $\mathbbm{f}_3\in H_0^1$, which is used to control the nonlocal term $ik \U'' \psi^I_k$ in \eqref{om_I}.

Firstly, we consider the Navier component $w_{\rm Na}$ and present the key estimates. 
For notational simplicity, we omit the $k$-dependence and define 
\begin{align}
	\label{E}
	\mathbb E[w,\phi]:=\nu^{\frac{1}{6}} |k|^{\frac{4}{3}} \| \na_ k \phi \|_{L^2} + \nu^{\frac{2}{3}} |k|^{\frac{1}{3}} \| \na_kw\|_{L^2} + \nu^{\frac{1}{3}} |k|^{\frac{2}{3}} \| w \|_{L^2}.
\end{align}
We summarize the key estimates as follows.

\begin{proposition}
\label{pro_Nav} Let  $w_{\mathrm{Na}; k}$ be the solution of \eqref{Nav_eqn} in the low frequency regime \eqref{low_reg}. Assume the parameter constraints $\|\mathcal{U}-y\|_{H^4} \leq \kappa$. There exists a threshold $0<\kappa_0(\U)=\mathcal{O}(1)$ and $0<\delta_0 <1$ such that if $0<\kappa\leq \kappa_0$ and $k \operatorname{Im} \lambda \geq-\delta_0 \nu^{1 / 3}|k|^{2 / 3}$, then the following estimates hold:
\begin{align}
	\label{F_in_L2_est}
	& \mathbb{E}\left[w_{\mathrm{Na}; k}, \phi_{\mathrm{Na}; k}\right]+\sqrt{|k|(\mathrm{Im}\lambda)_+}\|w_{\nv; k} \|_{L^2}+\nu\left\|\Delta_k w_{\mathrm{Na}; k}\right\|_{L_y^2}%\cancel{+\nu^{\frac{1}{4}}|k|\left\|w_{\mathrm{Na}}\right\|_{L_y^1}} 
	\lesssim
	\|F_k\|_{L_y^2}; 
	\\ 
	\label{F_in_H1_est}
	& \mathbb{E}\left[w_{\mathrm{Na}; k}, \phi_{\mathrm{Na}; k}\right]+\sqrt{|k|(\mathrm{Im}\lambda)_+}\|w_{\nv; k} \|_{L^2}+\nu\left\|\Delta_k w_{\mathrm{Na}; k}\right\|_{L_y^2}%\cancel{+\nu^{\frac{1}{4}}|k|\left\|w_{\mathrm{Na}}\right\|_{L_y^1}} 
	\lesssim \nu^{\frac{1}{6}}|k|^{-\frac{2}{3}}\left\|\nabla_k F_k\right\|_{L_y^2}; \\  \label{F_in_H-1_est}
	& \mathbb{E}\left[w_{\mathrm{Na}; k}, \phi_{\mathrm{Na}; k}\right]+\sqrt{|k|(\mathrm{Im}\lambda)_+}\|w_{\nv; k} \|_{L^2}%\cancel{+\nu^{\frac{1}{4}}|k|\left\|\nabla_k \phi_{\mathrm{Na}}\right\|_{L^{\infty}}}
	\lesssim \nu^{-\frac{1}{3}}|k|^{\frac{1}{3}}\|F_k \|_{H^{-1}_y}.
\end{align}	Here, the $H^{-1}$-norm is defined as \[ \| F \|_{ {H}^{-1}} := \sup \left\{ \langle F, g \rangle : g \in H_0^1([-1,1]), ~\| \na_k g\|_{L^2} :=(\|\pa_y g\|_{L^2}^2+|k|^2\|g\|_{L^2}^2)^{1/2}= 1 \right\}. \]
%	\linfeng{In CLWZ, the $H^{-1}$ norm is defined as $\Vert \partial_y g\Vert_{L^2}=1$ which is bigger. Check if it agrees with CWZ CMP paper}
\end{proposition}
Proposition \ref{pro_Nav} is analogous to \cite[Proposition 3.13]{CWZ23}. 
In \cite{CWZ23}, the authors employ a delicate compactness argument to prove these estimates under the assumption that the ambient shear profile $\mathcal{U}(y)$ is spectrally stable. 
When the shear flow is sufficiently close to the Couette flow, however, the singular integral operator $\mathfrak{J}_k$ introduced in \cite{BHIW23I} yields a simpler and more quantitative proof. 
For completeness, we present the proof of Proposition~\ref{pro_Nav} in Appendix \ref{Applemma}.

Next, we discuss the coefficients $c_\pm$ in \eqref{w_I_decmp}. Thanks to the formula \eqref{c_pm_frml}, one might expect that the estimates developed in Proposition \ref{pro_Nav} are sufficient to provide estimates for the coefficients. However, it turns out that these two coefficients are closely related to the $L^1_y$-bound of $w_{\rm Na}$, and the $\nu^{1/4}\|w_{\rm Na}\|_{L^1_y}$ control in \cite[Proposition 3.13]{CWZ23} is insufficient. To compensate for the loss of information, we keep track of an augmented quantity
\begin{align}
    \|(\U-\lambda) w_{\rm Na}-\U''\phi_{\rm Na}\|_{L^2}
\end{align}
in \eqref{Key_bound_2}.
A suitable control over this quantity implies the desirable bounds on $c_\pm$, see, e.g., Lemma \ref{lem:cpm_est}.  
This step is basically quantifying some inviscid damping in the linearized Euler equations which will be used to estimate the coefficients. 
Finally, we follow the argument in \cite{CWZ23} to develop bounds on the unit boundary correctors $w_\pm$. 

\ifx
\myr{$\nu^{1/4}\|w_{\rm Na}\|_{L^1}$ is not enough to guarantee the $\nu^{1/2}$-threshold. We need the $\nu^{1/6}\|w_{\rm Na}\|_{L^1}$ or $\nu^{1/6}|\pa_y \phi_{\rm Na}(\pm 1)|\lesssim \|F\|_{L^2}$. Moreover, it looks like the weaker estimate that we get on the $|k|\|(U-\lambda)w_{\rm Na}\|_{L^2}$ is enough to bound the $c_1,c_2$ in CLWZ.} 
\fi

Once the estimates mentioned above are carried out, we combine them to conclude the estimate for $w^I. $ 
\subsubsection{Decomposition of the Homogeneous Solution $ \omega_k^H$}
\label{sec:tmnd_H}
Next, we solve the homogeneous part of the solution $ \omega^H_{k}$. 
We begin by introducing the decomposition
\begin{align}
	& \omega^{H}_{k} (t,y)
	= \omega^{(1)}_{k} (t,y)
	+ \omega^{(2)}_{k} (t,y)
	+ \omega^{(3)}_{k} (t,y),
	\\&
	\Delta_k
	\phi_{k}^{(j)}
	=  \omega_k^{(j)},\quad 
	\phi^{(j)}\big|_{y=\pm 1}=0,
	\quad j=1,2,3.
\end{align}
We set $\omega_k^{(1)}$, which can be considered an approximation for the passive scalar evolution, as follows 
\begin{align}\label{eq:om_1}
	 \omega^{(1)}_k(t,y)=  \omega_{{\rm in};k}(y)\exp\lf\{- it\mathcal{U}(y)k-\frac{1}{3}(\mc U')^2\nu k^2 t^3-\nu k^2t\rg\}.
\end{align} 
The component $\omega^{(2)}_{k}$ corrects the residual which solves
\begin{align}
	\label{eq:om_2}
	\begin{split}
		&
		\pa_t  \omega_k^{(2)}
		+
		\U ik  \omega_k^{(2)}-\U''
		ik\psi_{k}^{(2)}
		-	\nu\de_k  \omega_k^{(2)}
	=
	\nu\pa_y^2 \omega_k^{(1)}+(\mc U')^2 \nu k^2t^2 \omega_k^{(1)}
	+
	ik\U''\phi^{(1)}_k,
	\\&
	\phi^{(2)}_k = \de_{k}^{-1}  \omega^{(2)}_k,
	\quad
	\int_{-1}^1  \omega_k^{(2)}e^{\pm k y}dy=0,\qquad
	\\&
	 \omega_{k}^{(2)}(t=0, y)=0.
	\end{split}
\end{align}
Note that the first two terms on the right-hand side of \eqref{eq:om_2}$_1$ equal 
\begin{align}
	(\nu\pa_y^2+(\mc U')^2 \nu k^2t^2) \omega_k^{(1)}
	&=
	\nu(\pa_y-i \U'kt )(\pa_y+i\U'kt) \omega_k^{(1)}-\nu i\U''kt  \omega_k^{(1)}	\\
	&=	\nu\Gamma_k^2 \omega_k^{(1)}-2\nu i\U'kt \Gamma_k \omega_k^{(1)}-\nu i \U'' kt \omega_k^{(1)}, 
\end{align}
where we denote the vector field $\Gamma = \partial_y - i\U' kt$, which commutes with the inviscid passive scalar equation. 
Since $ \omega_k^{(1)}$ experiences transport and enhanced dissipation, the term $\Gamma_k^2 \omega_k^{(1)}$ also experiences enhanced dissipation and this quantity remains small. We further remark that if the initial data is smooth enough, then the term $i\U'' \phi_k^{(1)}$ decays like $\brak{kt}^{-2}$ (independent of $\nu$) due to inviscid damping.   
The component $ \omega^{(3)}$ corrects the boundary condition which solves 
\begin{align}\label{eq:om_3}
	\begin{split}
		&
		\pa_t  \omega_k^{(3)}
		+
		\U ik  \omega_k^{(3)}-\U''ik\psi_{k}^{(3)}
	-\nu\de_k  \omega_k^{(3)}=0,
	\\&
	\int_{-1}^1  \omega_k^{(3)}e^{\pm k y}dy=
	-\int_{-1}^1  \omega_k^{(1)}e^{\pm k y}dy,
		\qquad
	 \omega_{k}^{(3)}(t=0, y)=0.
	\end{split}
\end{align} 

%Now, we can do a Laplace-Fourier transform in $t$ and use ideas ($c_\pm w_\pm$) in the inhomogeneous solution to provide estimates for $ w_k^{(3)}$.
%As a result, we have the following proposition.

Most of the estimates on $ \omega_k^H$ follow quickly from the corresponding estimates used in the proof of the inhomogeneous problem. This will eventually yield the following proposition. 
\begin{proposition}
	\label{pro:w_H}
Let $\omega_k^H$ be the solution of \eqref{om_H} in the low frequency regime \eqref{low_reg}. 
There exist constants $\delta>0$ and $0<\kappa,\nu_0\ll 1$ such that if $\Vert \U -y\Vert_{H^4}<\kappa$ and $\nu \in (0,\nu_0)$, then the following estimate holds:
\begin{align}\label{Est_omg_H_intro}
\begin{split}
	&\nu^{1/2}|k|
	\|e^{\delta(\nu k^2)^{1/3}t}
	\omega_k^H\|_{L^\infty L^2}^2
	+
	\Vert e^{\delta (\nu k^2)^{1/3} t} \omega^H_k\Vert_{L^2 \mathcal{Z}}^2
%	\nu^{1/3}|k|^{2/3}
%	\|e^{\delta(\nu k^2)^{1/3}t}
%	\rho_k^{1/2} \omega_k^H\|_{L^2L^2}^2
%	+
%	|k|^2\| e^{\delta(\nu k^2)^{1/3}t}
%	u^H_k\|_{L^2L^2}^2
%	+
%	\nu^{1/2}|k|
%	\|e^{\delta(\nu k^2)^{1/3}t}
%	 \omega_k^H\|_{L^2L^2}^2
%	\\&\quad
%	\\&
	\lesssim 
	(\log^2 L)
	|k|^{-2} \| \omega_{{\rm in};k}\|_{H^{2}_k}^2.
\end{split}
\end{align}
\end{proposition} 
The proof of the proposition is presented in Section~\ref{sec:space-time}. Combining the decomposition $\omega_k=\omega_k^I+\omega_k^H$ in \eqref{deco:1} and the estimates presented in Propositions \ref{pro:w_I} and~\ref{pro:w_H}, one obtains the key linear estimate \eqref{Lnr_st_est_low}. Hence, the proof for Proposition \ref{pro:lin_st:calU} is concluded.

\subsection{Linear Estimate: Time-dependent Shear Flows}
\label{sec:Frozen}
In this section, we present the frozen-time scheme to upgrade the estimate \eqref{Lnr_st_est_low} for time-independent shear flows to the bound \eqref{Lnr_low_est} presented in the key Proposition~\ref{pro:lin_low}, following the general framework outlined in \cite{ChenWeiZhang20}. It turns out that to prove the proposition, the key is to control the $L_t^2\mathcal{Z}$-norm of the vorticity. Once this time integral bound is developed, the uniform-in-time estimate of the vorticity in \eqref{Lnr_st_est_low} is straightforward to achieve. Hence, we will focus on $L_t^2\mathcal{Z}$-norm estimate throughout this section. 

We first rewrite the linearized system \eqref{lnr_Ns_k} equivalently as
%We treat the nonlinear term $u\cdot\na w$ in \eqref{vort-orig}$_1$ 
%as an external forcing and apply the Fourier transform in the $x$-variable.
%In this way, we are led to
%study the linearized Navier-Stokes equation
\begin{align}
	\label{Frz_1}
	\begin{aligned}
		&\pa_t \omega_k 
		+
		U(t,y)ik \omega_k 
		-
		U''(t,y)ik\psi_k 
		-
		\nu \de_k \omega_k =\na_k\cdot(\mathbbm{ f}_1,\mathbbm{f}_2),
		\\&
		\de_k\psi_k= \omega_k,\quad\psi_k |_{y=\pm 1}=0, \quad \int_{-1}^1 e^{\pm k y} \omega_k(t,y)dy = 0,\\&  \omega_k(t=0,y)= \omega_{{\rm in};k}(y).
    \end{aligned}
\end{align} 
Let $t_j = j\nu^{-1/3}$ and $\mathcal{I}_{\jj}= [t_j, t_{j+1})\cap [0,T]$.
There exists some $N\in\mathbb{N}$ such that
\begin{align}
	[0,T]
	=
	\bigcup_{j=0}^{N}
	\mathcal{I}_{[j]}.
\end{align}
%where $\mathcal{I}_{[j]}=[j \nu^{-1/3}, (j+1) \nu^{-1/3})$ for $0\leq j \leq N-1$ and $\mathcal{I}_{[N]}=[N \nu^{-1/3}, T]$.
%
%\begin{align}
%	\label{Frz_Ij}
%	[0,T]
%	=
%	\bigcup_{j=0}^\infty [j\nu^{-1/3},(j+1)\nu^{-1/3})=:\bigcup_{j=0}^\infty\mathcal{I}_{[j]}.
%\end{align}
%For each time horizon $[0,T]$, there exists $N\in\mathbb{N}$ such that
%\begin{align}
%	\bigcup_{j=0}^{N-1}\mathcal{I}_{[j]}
%	\subset [0,T]
%	\subset
%	\bigcup_{j=0}^{N}\mathcal{I}_{[j]} .
%\end{align}
For $0\leq j\leq N$, we define the ``frozen shear'' within each time interval $\mathcal{I}_{\jj}$ as
\begin{align}
	\label{Frz_Uj}
	U_\jj(y):=U(t_j,y).
\end{align}
We introduce the following layer-cake decomposition.
Let $ \omega_{[0]}$ be the solution of
\begin{align}
	\begin{aligned}
		&\pa_t  \omega_{[0]}+U_{[0]}ik  \omega_{[0]}-U''_{[0]}ik \psi_{[0]}
		-\nu \de_{k}  \omega_{[0]}
		=\na_k\cdot(\mathbbm{f}_1+\mathbbm{f}_{{\rm Disc} [0]}, \mathbbm{f}_2)\mathbbm{1}_{t\in\mathcal{I}_{[0]}},\label{Frz_om_0}
		\\&
		\de_k\psi_{[0]}= \omega_{[0]},
		\quad
		\psi_{[0]} |_{y=\pm 1}=0, \quad \int_{-1}^1 e^{\pm k y} \omega_{[0]} (t,y)dy= 0,
		\\& \omega_{[0]}(t=0,y)= \omega_{{\rm in};k}(y),
	\end{aligned}
\end{align}
where 
\begin{align}
	\mathbbm{f}_{{\rm Disc}[0]}
	:=
	(U_{[0]}(y)-U(t,y)) \omega_{[0]}
	-
	(U''_{[0]} (y)-U''(t,y))\psi_{[0]}
	.
	\label{discj2}
\end{align}
For $1\leq j \leq N$, let  $ \omega_{[j]}$ be the solution of
\begin{align}
\begin{aligned}
	&\pa_t  \omega_{\jj}+U_{\jj}ik  \omega_{\jj}-U''_{\jj}ik \psi_{\jj}
	-\nu \de_{k}  \omega_{[j]}
	=\na_k\cdot(\mathbbm{f}_1+
	\mathbbm{f}_{{\rm Disc} [j]} + \mathbbm{f}_{{\rm Frozen} [j]}, \mathbbm{f}_2)\mathbbm{1}_{t\in \mathcal{I}_{\jj}},
	\\&
	\de_k\psi_{[j]}= \omega_{[j]},
	\quad
	\psi_{[j]} |_{y=\pm 1}=0, \quad \int_{-1}^1 e^{\pm k y} \omega_{[j]} (t,y)dy= 0,
	\quad
	 \omega_{[j]}(t=0,y)=0,
\end{aligned}
\label{Frz_om_j}	
\end{align}
where
\begin{align}
	\mathbbm{f}_{{\rm Disc}[j]}
	&:=
	(U_\jj (y)-U(t,y))
	\sum_{j'=0}^{j}
	 \omega_{[j']}
	-
	(U''_\jj (y)-U''(t,y))
	\sum_{j'=0}^{j}
	\psi_{[j']}
	,
	\label{discj3}\\
	\mathbbm{f}_{\mathrm{Frozen} [j]}
	&:=
	\sum_{j'=0}^{j-1}
	\lf[
	(U''_\jj-U_\jjj'')\psi_\jjj
	-(U_\jj-U_\jjj) \omega_\jjj
	\rg].
	\label{frozenj}
\end{align}One can check that since the frozen shear $U_\jj (y)$ remains close to $U(t,y)$ on the time interval $\mathcal{I}_{\jj}$, the discretization error $\mathbbm{f}_{{\rm Disc} [j]}$ and frozen error $\mathbbm{f}_{{\rm Frozen}\jj}$ are small in the space  $L^2_t(\mathcal{I}_{\jj};L^2_y)$ (see Lemma~\ref{lem:Frzfrc}). 
%We emphasize that the domains of the systems \eqref{Frz_om_0} and \eqref{Frz_om_j} are $[-1,1]\times [0,T]$. eth
On each time interval $\mathcal{I}_\jj$ ($0\leq j\leq N$), it is readily checked that the solution of the linearized Navier-Stokes equations \eqref{Frz_1} can be written as
\begin{align}
	 \omega_k (t,y)
	=
	\sum_{j'=0}^{j}  \omega_{\jjj}.
	\label{Frz_om_dmp}
\end{align}The key observation here is that each component $\omega_{\jj}$ experiences a time-independent shear $U_\jj$. Hence, the linear estimate \eqref{Lnr_st_est_low} developed for \emph{time-independent} shear flows are valid. 
%For $0\leq j\leq N$,
%we consider the linearized Navier-Stokes equation
%\begin{align}
%	\label{Frz_2}
%	\begin{aligned}
%		&\pa_t w_{k} 
%		+
%		U_\jj(y)ik w_{k} 
%		-
%		U_\jj''(y)ik\psi_{k} 
%		-
%		\nu \de_k w_k
%		= \nabla_k \cdot (\mathbbm{f}_1 + \mathbbm{f}_{\rm Disc [j]}, \mathbbm{f}_2),
%		\\&
%		\de_k\psi_k= w_k,\quad\psi_k |_{y=\pm 1}=0, \quad \int_{-1}^1 e^{\pm k y} w_k(t,y)dy= 0,\\&  w_k(t=0,y)= w_{\mathrm{in};k}(y),
%	\end{aligned}
%\end{align} 
%where
%\begin{align}
%	\mathbbm{f}_{\rm Disc[j]}
%	:=
%	\lf[(U_\jj (y)-U(t,y)) w_k
%	-
%	(U''_\jj (y)-U''(t,y))\psi_k
%	\rg].
%	\label{discj}
%\end{align}
%Note that \eqref{Frz_2} is equivalent to \eqref{Frz_1}, with the ambient shear flow $U_\jj (y)$ frozen in time.

%Note that for each $j\in\mathbb{N}_0$, $ w_{[j]}$ is supported on the time interval $[j\nu^{-1/3},\infty)$ and hence $\mathbbm{f}_{\rm Frozen [j]}$ involves contributions from all $ w_{[j']}$, where $j' \leq j$. 
%Motivated by Proposition \ref{pro:lin_st}, 
Based on the decomposition \eqref{Frz_om_dmp}, we would like to design two types of functional. The first type of functional (``horizontal functional" $X_\jj$) captures the sharp enhanced dissipation and inviscid damping for each $\omega_\jj$-component. Bounds on the $X_\jj$-functional can be derived from the linear estimate for the \emph{time-independent shear flows} (Proposition \ref{pro:lin_st:calU}).  On the other hand, we need another type of functional (``vertical functional" $Y_\jj$) to collect all the contributions from the  decomposition \eqref{Frz_om_dmp} on each time interval $\mathcal I_\jj$. 
%The horizontal functional $X_\jj$ and the vertical functional $Y_{\jj}$. The horizontal norm $X_\jj$ captures the ``mass'' of each layer of the cake, while the vertical norm $Y_{\jj}$ quantifies the weighted mass of the solution on $\mathcal{I}_{[j]}$.

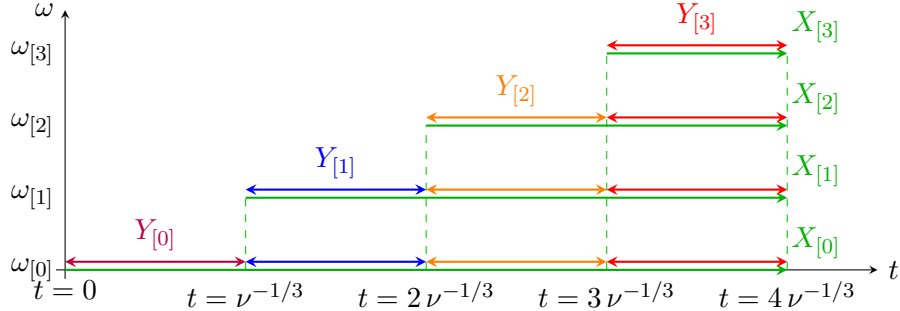
\begin{figure}[h]
\centering
\begin{tikzpicture}[%
    scale=1.2,     % Overall scale
    >=stealth,     % For nicer arrow tips
    every node/.style={font=\small}
]

%-- Variables for easy adjustments
\def\ttick{2}    % Horizontal spacing to represent \nu^{1/3}
\def\vtick{0.8}    % Vertical spacing between w_{[i]}
\def\gap{0.15}     % Small vertical gap for pink horizontal arrows

%-- Axes
\draw[->] (-0.1,0) -- (4.5*\ttick,0) node[right] {$t$};
\draw[->] (0,-0.1) -- (0,3.6*\vtick) node[left] {$ \omega$};

%-- Vertical dashed lines at t = nu^{-1/3}, 2 nu^{-1/3}, ...
\draw[dashed,green!80!black] (\ttick,0) -- (\ttick,\vtick);
\draw[dashed,green!80!black] (2*\ttick,0) -- (2*\ttick,2*\vtick);
\draw[dashed,green!80!black] (3*\ttick,0) -- (3*\ttick,3*\vtick);
\draw[dashed,green!80!black] (4*\ttick,0) -- (4*\ttick,3*\vtick);

%-- Labels along t-axis
\node[below] at (0,0)           {$t=0$};
\node[below] at (\ttick,0)      {$t=\nu^{-1/3}$};
\node[below] at (2*\ttick,0)    {$t=2\,\nu^{-1/3}$};
\node[below] at (3*\ttick,0)    {$t=3\,\nu^{- {1}/{3}}$};
\node[below] at (4*\ttick,0)    {$t=4\,\nu^{-{1}/{3}}$};

%-- w-level labels on the left
\node[left]  at (0,0)         {$ \omega_{[0]}$};
\node[left]  at (0,1*\vtick)  {$ \omega_{[1]}$};
\node[left]  at (0,2*\vtick)  {$ \omega_{[2]}$};
\node[left]  at (0,3*\vtick)  {$ \omega_{[3]}$};

%-- Green horizontal step lines (just drawn full-length here)
\draw[thick,green!70!black,->] (0,0) -- (4*\ttick,0)node[right,above] {\qquad $X_{[0]}$};
\draw[thick,green!70!black,->] (1*\ttick,1*\vtick)-- (4*\ttick,1*\vtick) node[right,above] {\qquad $X_{[1]}$};
\draw[thick,green!70!black,->] (2*\ttick,2*\vtick) -- (4*\ttick,2*\vtick)node[right,above] {\qquad $X_{[2]}$};
\draw[thick,green!70!black,->] (3*\ttick,3*\vtick) -- (4*\ttick,3*\vtick)node[right,above] {\qquad $X_{[3]}$};

%-- Pink vertical lines
%\draw[thick,pink] (2*\ttick,0) -- (2*\ttick,\vtick);
%\draw[thick,pink] (3*\ttick,0) -- (3*\ttick,2*\vtick);
%\draw[thick,pink] (4*\ttick,0) -- (4*\ttick,2*\vtick);

%-- Pink horizontal arrows (labeled Y[...] for illustration)

\draw[<->,thick,purple]
  (0*\ttick,0*\vtick+0.6*\gap) -- (1*\ttick,0*\vtick+0.6*\gap)
  node[midway,above] {$Y_{[0]}$};

\draw[<->,thick,blue]
  (1*\ttick,1*\vtick+0.6*\gap) -- (2*\ttick,1*\vtick+0.6*\gap)
  node[midway,above] {$Y_{[1]}$};
\draw[<->,thick,blue]
  (1*\ttick,0*\vtick+0.6*\gap) -- (2*\ttick,0*\vtick+0.6*\gap);

\draw[<->,thick,orange]
  (2*\ttick,2*\vtick+0.6*\gap) -- (3*\ttick,2*\vtick+0.6*\gap) node[midway,above] {$Y_{[2]}$}
;
\draw[<->,thick,orange]
  (2*\ttick,1*\vtick+0.6*\gap) -- (3*\ttick,1*\vtick+0.6*\gap) ;
\draw[<->,thick,orange]
  (2*\ttick,0*\vtick+0.6*\gap) -- (3*\ttick,0*\vtick+0.6*\gap);

\draw[<->,thick,red]
  (3*\ttick,3*\vtick+0.6*\gap) -- (4*\ttick,3*\vtick+0.6*\gap)
  node[midway,above] {$Y_{[3]}$};
\draw[<->,thick,red]
  (3*\ttick,2*\vtick+0.6*\gap) -- (4*\ttick,2*\vtick+0.6*\gap);
\draw[<->,thick,red]
  (3*\ttick,0*\vtick+0.6*\gap) -- (4*\ttick,0*\vtick+0.6*\gap);
  \draw[<->,thick,red]
  (3*\ttick,1*\vtick+0.6*\gap) -- (4*\ttick,1*\vtick+0.6*\gap);
  
\end{tikzpicture}
\caption{The Frozen Norms} 
\end{figure}
Let us start with the horizontal functional $X_\jj$, which measures the size of the solution $ \omega_\jj$. 
Since $\omega_\jj$ initiates from time $t_j$, the suitable enhanced dissipation time weight from Proposition \ref{pro:lin_st:calU} is $e^{\delta\nu^{1/3}(t-t_j)}$. 
This leads to the definition:
\begin{align}
	\label{Frz_X_norm}
	X_{\jj}
	%\norm{ \omega_{\jj}}_{X_{\jj}}^2
	:=
	\|e^{\delta\nu^{1/3}(t-t_j)} \omega_{\jj}\|_{L_t^2([t_j,T];\mathcal{Z})},
	\quad j\geq 0.
\end{align} 
We observe that the $X_\jj$-estimate can be naturally derived from Proposition \ref{pro:lin_st:calU}. 

Once we have suitable controls over the horizontal functional $X_\jj$, we can bound the following vertical functional (note that the norm is over $\mathcal{I}_{[j]}$)
\begin{align}
	\label{Frz_Y_norm}
	Y_{\jj}
	:=
	\sum_{j'=0}^{j}(j-j'+1)
	\| \omega_{\jjj}\|_{L^2(\mathcal{I}_{[j]};\mathcal{Z})},
	\quad
	j\geq 0.
\end{align} 
Here, the linear-in-$j$ weight is introduced to control the frozen error $\mathbbm{f}_{\rm Frozen}$.  We postpone the rigorous estimate of $Y_\jj$ to Section \ref{sec:frozen}. 
%Note that the the time-integral is over $\mathcal{I}_j$ regardless of $j'$.  note that there is no exponential-in-time weight in the functional and

Finally, to complete the $L_t^2 \mathcal{Z}$ estimate for the vorticity $\omega_k$ in \eqref{Frz_om_dmp}, we collect all the $Y_\jj$-functional estimates together to obtain control over the following total energy
\begin{align}
	\label{Frz_E}
	\EE
	:=\sum_{j=0}^Ne ^{2\delta_\ast j}Y_{\jj}^2.
\end{align} Here, the $e ^{2\delta_\ast j}$ is a suitable enhanced dissipation weight  
and $\delta_\ast>0$ is a constant related to the $\delta$ in  Proposition~\ref{pro:lin_st:calU}.
One observes that the total energy $\mathcal E$ controls the $L_t^2 \mathcal{Z}$-norm of $\omega_k$ (Lemma~\ref{LemmaE}), i.e.,
\begin{align}
	\|e^{\delta_\ast\nu^{1/3}t}
	 \omega_k\|_{L^2([0,T];\mathcal{Z})}^2
	\lesssim 
	\EE.
	\label{EQE}
\end{align}
%where $ \omega_k$ denotes the solution of the linearized Navier-Stokes equations \eqref{Frz_1} which is given in \eqref{Frz_om_dmp} on the time interval $\mathcal{I}_{\jj}$.
% \textcolor{red}{[so its accurate to say we want the double sum]}
% \begin{align*}
% E & = \sum_{j=0}^N \sum_{j'= 0}^j e ^{6\delta_\ast j} (|j'-j|+1)\| w_{\jjj}\|_{L^2(\mathcal{I}_{j};\mathcal{Z})} \\ 
% & = \sum_{j'=0}^N \| w_{\jjj}\|_{L^2(\mathcal{I}_{j};\mathcal{Z})}^2 \left(\sum_{j = j'}^N e^{6\delta_\ast j}(|j'-j|+1)\right).  
% \end{align*}
% \textcolor{red}{[That seems slightly wrong to me, since]}
% \begin{align*}
% 	\left(\sum_{j = j'}^N e^{6\delta_\ast j}(|j'-j|+1)\right). \approx e^{6\delta_\ast N} N???
% \end{align*}
This completes the estimate of the $L^2\mathcal Z$ norm of $ \omega_k$. We postpone the estimate of the quantity $\|(1-|y|)^{1/2}\omega_k\|_{L^\infty_t L_y^2}$ to Section \ref{sec:frozen} because it follows similar line of argument. This concludes the sketch for the proof of Proposition \ref{pro:lin_st}.

\subsection{Nonlinear Stability}
\label{sec:nlnr_outln}
Once we establish suitable estimates for the linearized Navier-Stokes equations \eqref{lnr_Ns_k}, we exploit the quadratic nature of the nonlinearity in the full equation to absorb its effects via the coercive terms in the linear estimates, provided that the perturbation remains sufficiently small. %Finally, to establish the nonlinear stability of the system \eqref{vort-orig}, we invoke Proposition \ref{pro:lin_dy}, which is a consequence of Proposition \ref{pro:lin_st}, and 
An adaptation of the nonlinear stability arguments in \cite[Section 7]{CLWZ20} yields Theorem \ref{T01}; see Section \ref{sec:NL} for details. 

%\newpage

\section{The Linearized Problem}  
\label{sec:Lin}
\subsection{Linear Estimates of the Navier Part}\label{sec:lnr_Nav}
In this section, we develop various linear estimates for the Navier part of the inhomogeneous system \eqref{Nav_eqn}.

To estimate the boundary corrector, the following lemma is crucial.   

\begin{lemma}
\label{L02}
Let $w_{\rm Na}$ be the solution of \eqref{Nav_eqn}.
Let $\lambda=\lambda_r +i \lambda_i$.
Assume the parameter constraints $\|\U-y\|_{H^4}\ll 1$, $\nu|k|^2\leq \|\U'\|_{L^\infty}$ and spectral constraint $0\geq k{\rm Im} \lambda\geq -\delta_0 \nu^{1/3} |k|^{2/3}$, we have
\begin{align}
	\label{Key_bound}
	&\|(\mathcal{U}-\lambda_r) w_{\rm Na}-\mathcal{U}''\phi_{\rm Na}\|_{L^2}\lesssim \|F_k\|_{L^2},
	\\
	\label{Key_bound_2} 
	&
	\|(\mathcal{U}-\lambda) w_{\rm Na}-\mathcal{U}''\phi_{\rm Na}\|_{L^2}
	\lesssim \|F_k\|_{L^2}.
\end{align}
\end{lemma}
%\linfeng{remove the comments below?}
%With the help of this quantity, we might be able to control the $\nu^{1/6}\|w_{\rm Na}\|_{L^1}$ for the parameter range $y\in [y_\lambda-L^{-1}, y_\lambda+L^{-1}]\cup [y_\lambda-\kappa|k|^{-1}, y_\lambda+\kappa|k|^{-1}]^c$, where $\mathcal{U}(y_\lambda)=\lambda.$ \myr{(Check??)} We run the same argument as before.

\begin{proof}
Without loss of generality, we assume that $k>0$.
Let $\lambda_i=-\delta\nu^{1/3}|k|^{-1/3}$, where $\delta\in [0,\delta_0]$.
We rewrite \eqref{Nav_eqn} as
\begin{align}  
\begin{split}
		&-\nu\de_k w_{\rm Na} 
	+
	ik\lf(\mathcal{U}-\lambda_r\rg)w_{\rm Na} 
	-
	i\mathcal{U}''k \phi_{\rm Na}-\delta\nu^{1/3}|k|^{2/3}w_{\rm Na} =F_k,
	\\&
	\phi_{{\rm Na}}
	=\de_{D}^{-1} w_{\rm Na},
	\quad
	w_{\rm Na} \big|_{y=\pm 1}=0.
\end{split}
	\label{lemma_2_4}
\end{align}
We test \eqref{lemma_2_4}$_1$ by $\overline{(\mathcal{U}-\lambda_r)w_{\rm Na}-\mathcal{U}''\phi_{\rm Na}}$ and take the imaginary part to obtain
\begin{align}
\begin{split}
	&
	|k|
	\int |(\mathcal{U}-\lambda_r)w_{\rm Na}-\mathcal{U}''\phi_{\rm Na}|^2 dy
	\\&
	\les
	\left| \int F_k \overline{((\mathcal{U}-\lambda_r)w_{\rm Na}
 	-
 	\mathcal{U}''\phi_{\rm Na})}dy	 	\right|
	+
	\left|
	\nu{\rm Im} \int\de_kw_{\rm Na } \overline{((\mathcal{U}-\lambda_r)w_{\rm Na}-\mathcal{U}''\phi_{\rm Na})}dy  \right|
	\\&\quad
	+
	\delta 
	\nu^{1/3}|k|^{2/3}
	\left|
	\int w_{\rm Na} \overline{((\mathcal{U}-\lambda_r)w_{\rm Na}
	-
	\mathcal{U}'' \phi_{\rm Na})}dy
	\right|
	\\&
	\les
	\big(
	|k|^{-1/2}\|F_k\|_{L^2}
	+
	\delta\nu^{1/3}
	|k|^{1/6}
	\| w_{\rm Na}\|_{L^2}
	\big)
	|k|^{1/2}
	\|(\mathcal{U}-\lambda_r)w_{\rm Na}
	-
	\mathcal{U}''\phi_{\rm Na}\|_{L^2}
	\\&\quad 
	+
	\nu
	\left|
	\int \pa_yw_{\rm Na} \mathcal{U}'
	\overline{w_{\rm Na}}
	\right|
	+
	\nu
	\left|
	\int\na_kw_{\rm Na}\cdot\overline{\na_k(\mathcal{U}''\phi_{\rm Na})}
	dy\right|.	 		
\end{split}
\end{align}
From Young's inequality it follows that
\begin{align}
\begin{split}
	&	
	|k|
	\Vert (\mathcal{U}-\lambda_r)w_{\rm Na}-\mathcal{U}''\phi_{\rm Na}
	\Vert_{L^2}^2
	\\& 
	\les
	C_{\eps} |k|^{-1}\|F_k\|_{L^2}^2+
	\eps
	|k|
	\|(\mathcal{U}-\lambda_r)w_{\rm Na}-\mathcal{U}''\phi_{\rm Na}\|_{L^2}^2
	+
	\nu^{ 4/3}|k|^{-1/3}\|\na_kw_{\rm Na}\|_{L^2}^2
	\\&\qquad
	+
	C_{\eps}
	\nu^{2/3}|k|^{1/3}\|w_{\rm Na}\|_{L^2}^2
	%	\\ &
	%	\les |k|^{-1}\|\mathbb{F}\|_{L^2}^2
	%	 +
	%	|k|
	%	 \|((U-\lambda_r)w_{\rm Na}-U''\phi_{\rm Na})\|_{L^2}^2
	%	 +
	%	 |k|^{-1}\|\mathbb{F}\|_{L^2}\nu^{1/3}|k|^{2/3}\|w_{\rm Na}\|_{L^2}
	\\&
	\les
	C_{\eps} |k|^{-1}\|F_k\|_{L^2}^2+
	\eps
	|k|
	\|(\mathcal{U}-\lambda_r)w_{\rm Na}-\mathcal{U}''\phi_{\rm Na}\|_{L^2}^2,
\label{lemma_2_10}
\end{split}
\end{align}
for any $\eps \in (0,1)$, where we used Proposition~\ref{pro_Nav} in the last step. We take $\eps>0$ in \eqref{lemma_2_10} sufficiently small to get
\begin{align}
	{|k|}^{2}
	\| 
	(\mathcal{U}-\lambda_r)w_{\rm Na}
	-
	\mathcal{U}''\phi_{\rm Na}
	\|_{L^2}^2
	\les 
	\|F_k\|_{L^2}^2,
\end{align}
which completes the proof of \eqref{Key_bound}.
Finally, we use Proposition~\ref{pro_Nav} to obtain
\begin{align}
\begin{split}
	&
	{|k|}^{2}
	\|(\mathcal{U}-\lambda)w_{\rm Na}
	-
	\mathcal{U}''\phi_{\rm Na}\|_{L^2}^2
	\\&
	\les
	{|k|}^{2}\|(\mathcal{U}-\lambda_r)w_{\rm Na}-\mathcal{U}''\phi_{\rm Na}\|_{L^2}^2
	+
	\nu^{2/3}|k|^{4/3}\|w_{\rm Na}\|_{L^2}^2
%	\\&
	\les
	\|F_k\|_{L^2}^2.\label{Key_bnd}
\end{split}
\end{align}
This completes the proof of the lemma.
\end{proof}

\begin{lemma}
\label{pro_Nav_hi} 
Let $w_{\rm Na}$ be the solution of \eqref{Nav_eqn}.
Assume the parameter constraints $\|\U-y\|_{H^4}\ll 1$, $\nu|k|^2\geq  \|\mathcal{U}'\|_{L^\infty}$ and the spectral constraint $0\geq k{\rm Im}\lambda\geq -\delta_0 \nu^{1/3}|k|^{2/3}$, we have
\begin{align}
	\label{F_in_L2_hi}
	&\nu k^2\norm{w_{\rm Na}}_{L^2}
	\lesssim \enorm{F_k},
	%		\mathbb{E}[w_{\rm Na;k}, \phi_{\rm Na;k}]+\nu\|\de_k w_{\rm Na;k}\|_{L_y^2}+\nu^{\frac{1}{4}}|k|\|w_{\rm Na;k}\|_{L_y^1}%+\myr{\underbrace{\frac{|k|\|(U-\lambda)w_{\rm Na;k}\|_{L^2}}{\kappa\nu^{-\frac{1}{6}}|k|^{-\frac{4}{3}}+1}}_{Check!??}}
	%		\lesssim \|\mathbb{F}_k\|_{L^2};
	\\
	\label{F_in_H1_hi} &
	\nu k^2 \enorm{u_{\rm Na}} \lesssim \norm{F_k}_{H^{-1}},
	\\&
	\nu |k| \enorm{w_{\rm Na}} \lesssim \norm{F_k}_{H^{-1}}.
	\label{F_in_H-1_w_hi}
\end{align}
\end{lemma}
\begin{proof}
	Without loss of generality, we assume that $k>0$.
	We test \eqref{Nav_eqn}$_1$ with $\overline{\phi_{{\rm Na}}}$ and integrate by parts to obtain	
	\begin{align*}
		-\nu\norm{w_{\rm Na}}_{L^2}^2 
		+ 
		ik\int (\mathcal{U} - \lambda) w_{{\rm Na}} \overline{\phi_{{\rm Na}} }\, dy  
		- 
		ik \int  \mathcal{U}'' \phi_{{\rm Na}} \overline{\phi_{{\rm Na}}}\,dy 
		= 
		\int F_k \overline{\phi_{{\rm Na}}}\,dy.
	\end{align*}
	Taking the real part of the above equation, we get
	\begin{align*}
		\nu&\norm{w_{\rm Na}}_{L^2}^2 
		+ 
		{\rm Im} 
		 k\int (\mathcal{U} 
		- \lambda) w_{{\rm Na}} \overline{\phi_{{\rm Na}} }\, dy 
		\\&
		=\nu\norm{w_{\rm Na}}_{L^2}^2 - {\rm Im} k\int\mathcal{U}' \phi'_{{\rm Na}} \overline{\phi_{{\rm Na}} }\, dy 
		+
		{\rm Im}\lambda k\int |\phi'_{{\rm Na}}|^2 \, dy + {\rm Im}\lambda k^3 \int |{\phi_{{\rm Na}} }|^2\, dy 
		\\&
		= -{\rm Re} \int F_k \overline{\phi_{{\rm Na}}}\,dy,
	\end{align*}
	from where
	\begin{align*}
		\nu\norm{w_{\rm Na}}_{L^2}^2 
		\leq &
		\abs{\int{F_k}\overline{\phi_{{\rm Na}}}\,dy} 
		+ k\norm{\mathcal{U}'}_{L^\infty_y}\enorm{\phi_{{\rm Na}}'}\enorm{\phi_{{\rm Na}}} 
		+ \delta_0\nu^{1/3}|k|^{2/3} \enorm{\phi_{{\rm Na}}'}^2
		\\& 
		+ \delta_0\nu^{1/3}|k|^{8/3} \enorm{\phi_{{\rm Na}}}^2.
	\end{align*}
	Using $ \|\mathcal{U}'\|_{L^\infty}\leq \nu|k|^2$ and the identity
	\begin{equation}
		\label{identity}
		\norm{w_{\rm Na}}_{L^2}^2 = \enorm{\phi_{{\rm Na}}''}^2 
		+ 
		k^4\enorm{\phi_{{\rm Na}}}^2 
		+ 
		2 k^2\enorm{\phi_{{\rm Na}}'}^2,
	\end{equation}
	we obtain 
	\begin{align}
		\label{mid:step}
		\nu
		\norm{w_{\rm Na}}_{L^2}^2 
		\lesssim 
		\abs{\int{F_k}
		\overline{\phi_{{\rm Na}}}\,dy} 
		\lesssim
		\enorm{F_k}
		\enorm{\phi_{{\rm Na}}}.
	\end{align}
%\linfeng{we used $\Vert U'\Vert_{L^\infty}$ close to 1}
	Note that \eqref{identity} implies
	\begin{equation*}
		\enorm{\phi_{{\rm Na}}} 
		\le 
		\frac{1}{k^2} \enorm{w_{\rm Na}},
	\end{equation*}
and thus
\eqref{F_in_L2_hi} follows.
To prove~\eqref{F_in_H1_hi}, we use~\eqref{identity} and \eqref{mid:step} to obtain
\begin{align*}
	&\nu k^2\norm{u_{\rm Na}}_{L^2}^2
	\lesssim 
	\nu\norm{w_{\rm Na}}_{L^2}^2 
	\lesssim \norm{F_k}_{H^{-1}}
	\Vert \phi_{\rm Na}\Vert_{H^1}
	\les  
	\norm{F_k}_{H^{-1}}
	\enorm{u_{{\rm Na}}}.
\end{align*}
Finally, we note that
\begin{align*}
	\nu\norm{w_{\rm Na}}_{L^2}^2 
	\lesssim 
	\norm{F_k}_{H^{-1}} \enorm{u_{{\rm Na}}} 
	\lesssim 
	\frac{1}{|k|}
	\norm{F_k}_{H^{-1}} \enorm{w_{{\rm Na}}},
\end{align*}
and~\eqref{F_in_H-1_w_hi} follows. 
\end{proof}

\subsection{Estimates of the Coefficients $c_{+}$ and $c_{-}$}
In this section, we establish several estimates on the coefficients $c_{\pm}$ that were introduced in the decomposition \eqref{w_I_decmp}.
We aim to derive some analogous estimates for the coefficients $c_1$ and $c_2$ in \cite[Lemmas 3.6, 4.1, and~4.2]{CLWZ20}. 
In addition, we present the necessary modifications to treat the extra nonlocal term.

Since $\U\in C^3$ and $\Vert \U -y\Vert_{H^4}\ll 1$, there exists a strictly increasing function $\wt\U\in C^3(\rr)$ such that $\wt \U (y)= \U (y)$ on $[-1,1]$, $\lim_{y\to -\infty} \wt \U (y)=-\infty$, and $\lim_{y\to \infty} \wt \U (y)=\infty$.
Hence,
for any $\lambda_r\in\mathbb{R}$,
there exists a unique $y_\lambda\in\rr$ such that $\wt \U (y_\lambda)=\lambda_r$.
For $\lambda_r\in\rr$, we define a cutoff function $\chi$ on $\rr$ such that
\begin{align}
	\label{defn_chi}
	\chi(y;\lambda_r)
	=
		\begin{aligned}
			&\frac{1}{\wt{\mathcal{U}}(y)-\lambda_r},\quad \text{if~} |y-y_\lambda|\geq L^{-1}.
		\end{aligned}
\end{align} 

%\begin{align}
%	\label{defn_chi}
%	\chi(y;\lambda_r)
%	=
%	\begin{cases}
%		\begin{aligned}
%	 &\frac{1}{\wt{\mathcal{U}}(y)-\lambda_r},\quad \text{if~} |y-y_\lambda|\geq L^{-1};
%	\\[2ex]&
%	2 L^2 (\wt{\U} (y) -\lambda_r)- L^4 (\wt{\U} (y) -\lambda_r)^3
%	,
%	\quad\text{if~}
%	|y-y_\lambda| 
%	\leq L^{-1}.
%	\end{aligned}
%\end{cases}
%\end{align} 

%\begin{align}
%	\label{defn_chi}
%	\chi(y;\lambda_r)
%	=
%	\begin{cases} \frac{1}{\wt{\mathcal{U}}(y)-\lambda_r},\quad \text{if~} |y-y_\lambda|\geq L^{-1},\\ \text{smooth},\quad \text{otherwise.}
%	\end{cases}
%\end{align} 

\begin{lemma}
\label{L01}	
Let $w_{\rm Na}$ be the solution of \eqref{Nav_eqn} and $g_\pm\in H^{1}(-1,1)$ be test functions with the boundary condition $g_{+} (-1)=0$ and $g_-(+1)=0$. 
%Assume that the test functions $g_\pm\in H^{1}(-1,1)$ with the boundary condition $g_\pm(\mp 1)=0$. 
 For $\nu|k|^2\leq \Vert \mathcal{U}'\Vert_{L^\infty}$ and $0\geq k {\rm Im} \lambda \geq -\delta_0 \nu^{1/3} |k|^{2/3}$, we have
\begin{align}
\begin{split}
	&
	|\lan w_{\rm Na}, g_\pm\ran|
	%	\\&
	\lesssim 
	|k|^{-1}
	\| F_k\|_{H^{-1}}
	\big(
	L^{3/2} 
	\|g_\pm \|_{L^\infty ((-1,1)\cap (y_\lambda-L^{-1}, y_\lambda+L^{-1}))}
	\\&\hspace{2.3cm}
	+
	|g_\pm (\pm 1)|(|\mathcal{U}(\pm 1)-\lambda_r|
	+L^{-1})^{-3/4}L^{3/4}
	+
	\|g_\pm \chi\|_{H^1}
	\\&\hspace{2.3cm}
	+L\|g_\pm \chi\|_{L^2}+\nu^{-1/2}\|g_{\pm}\chi\|_{L^1}
	\big).
	\label{pairing}
\end{split}
\end{align}
\end{lemma}

The above lemma is analogous to \cite[Lemma~3.6]{CLWZ20}. 
The main modification relative to \cite[Lemma~3.6]{CLWZ20} lies in the appearance of the $L^1$-norm of $g_{\pm}\chi$ on the right-hand side.

\begin{proof}
Let $\eta_1\in H_0^1 (-1,1)$ be a test function.
We test \eqref{Nav_eqn}$_1$ with $\overline{\eta_1}$ to obtain
\begin{align}
\begin{split}
	&
	|k|
	\lf|\lf\lan (\mathcal{U}-\lambda_r) w_{\rm Na}-
	\mathcal{U}'' \phi_{\rm Na},\eta_1\rg\ran\rg|
	\\&\quad
	\lesssim
	\|F_k \|_{H^{-1}}
	\| \nabla_k \eta_1\|_{L^2}
	+
	\nu\|\pa_y w_{\rm Na}\|_{L^2}\|\pa_y \eta_1\|_{L^2}
	+
	\nu^{1/3}|k|^{2/3}\|w_{\rm Na}\|_{L^2}\|\eta_1\|_{L^2}
	\\&\quad
	\les
	\|F_k \|_{H^{-1}}
	(\|\pa_y \eta_1\|_{L^2}
	+
	L
	\Vert \eta_1\Vert_{L^2}
	).
	\label{ast_L1}
\end{split}
\end{align}
where we used \eqref{F_in_H-1_est} in the last step.
Let $G_k$ be the Green's function for $\de_k=-|k|^2+\partial_{yy}$ with homogeneous Dirichlet boundary condition, which has the explicit form
\begin{align}
\begin{split}
	G_k (y,y')
	=
	-\frac{1}{k \sinh (2k)}
	\begin{cases}
	\sinh (k(1-y')) \sinh (k(1+y))
	,\qquad y\leq y';
	\\
	\sinh (k(1-y)) \sinh (k(1+y'))
	,\qquad y\geq y'.
	\end{cases}	
	\label{green}
\end{split}
\end{align}
From the symmetric property $G_k(y,z)=G_k(z,y)$ it follows that
\begin{align}
\begin{split}
	\lan \de_k^{-1}g,h\ran
	&=
	\int\int G_k(y,z)g(z)dz\overline{h(y)} dy
	%	\\&
	=\int g(z)\overline{\int G_k(z,y)h (y)dy} dz
	\\&
	=\lan g, \de_k^{-1}h\ran.
	\label{symmetry}
\end{split}
\end{align} 
Using \eqref{ast_L1} and \eqref{symmetry}, we conclude
\begin{align}
\begin{split}
	\lf|\lf\lan w_{\rm Na} , (\mathcal{U}-\lambda_r)\eta_1 -\de_k^{-1}(\mathcal{U}''\eta_1)
	\rg\ran\rg|
	&=
	\lf|\lf\lan (\mathcal{U}-\lambda_r) w_{\rm Na}
	-
	\mathcal{U}'' \phi_{\rm Na},\eta_1\rg\ran\rg|
	\\&
	\les
	|k|^{-1}\| F_k \|_{H^{-1}}(\|\pa_y \eta_1\|_{L^2}
	+ 
	L\|\eta_1\|_{L^2}),\quad \label{pairing_1}
\end{split}
\end{align}
for all $\eta_1\in H^1_0(-1,1)$.

Let $\eta_2\in H^1 (-1,1)$ be a test function with $\eta_2 (-1)=0$.
Let $L_{\ast} \in [L,\infty)$ be a constant to be determined later.
Let
$
	\chi^\ast_{+} (y)
	:=
	\max\{0,L_\ast(y-1+L^{-1}_\ast)\}
$ be defined on $[-1,1]$.
% Moreover, we also choose $\chi^\ast:=\chi^\ast_+$. 
It follows that
\begin{align}
\begin{split}
	\|(\mathcal{U}-\lambda_r)\chi^\ast_{+}\|_{L^\infty}
	&\leq
	\Vert \mathcal{U}-\lambda_r\Vert_{L^\infty ([1-L_{\ast}^{-1}, 1])}
	\Vert \chi_{+}^{\ast} \Vert_{L^\infty}
	%	\\&
	\leq |\mathcal{U}(1)-\lambda_r|+L_\ast^{-1} \|\mathcal{U}'\|_{L^\infty}.
	\label{pairing_34} 
\end{split}
\end{align} 
Since $w_{\rm Na}|_{y=\pm 1}=0$, the fundamental theorem of calculus implies
\begin{align*}
	|w_{\rm Na}(y)|
	\leq 
	L_\ast^{-1/2}
	\|\pa_y w_{\rm Na}\|_{L^2},
	\quad
	\forall y\in [1-L_\ast^{-1},1],
\end{align*}
from where
\begin{align}
\|w_{\rm Na}\|_{L^1([1-L_\ast^{-1},1])}\lesssim L_\ast^{-3/2}\|\pa_y w_{\rm Na}\|_{L^2}.\label{pairing_3/2}
\end{align}
%\myr{For a test function $f\in H^1(-1,1)$, we define $f^\ast_{+}:=f-f(1)\chi^\ast_{+}$ and $f^\ast_{-}:=f-f(-1)\chi^\ast_{-}$.}
Denote $\eta_{2}^{\ast}
=\eta_{2}-\eta_{2} (1) \chi^{\ast}_{+}$.
Direct computation gives
\begin{align}
\begin{split}
	&\lf|
	\lf
	\lan w_{\rm Na}, (\mathcal{U}-\lambda_r)\eta_2-\de_k^{-1}(\mathcal{U}''\eta_2)\rg\ran\rg|
	\\&
	\leq 
	\lf|\lf\lan w_{\rm Na}, (\mathcal{U}-\lambda_r)\eta_{2}^{\ast}
	-
	\de_k^{-1}(\mathcal{U}''\eta_{2}^{\ast})\rg\ran\rg|
	+
	|\eta_2(1)|\lf\lan w_{\rm Na}, (\mathcal{U}-\lambda_r)\chi_{+}^{\ast}\rg\ran|
	\\&\qquad
	+
	|\eta_2(1)|\lf|\lf\lan w_{\rm Na},  \de_k^{-1}(\mathcal{U}''\chi_{+}^{\ast})\rg\ran\rg| 
	\\&=:
	T_1+T_2+T_3.
	\label{pairing_2}
\end{split}
\end{align}
Note that $\eta_2^\ast\in H^1_0$.
From
\eqref{pairing_1} it follows that
\begin{align}
\begin{split}
	|T_1|&
	\lesssim 
	|k|^{-1}\| F_k\|_{H^{-1}}(\|\pa_y \eta_2^\ast\|_{L^2}+L\|\eta_2^\ast\|_{L^2})\\
	&\lesssim  
	|k|^{-1}\|
	F_k \|_{H^{-1}}(\|\pa_y \eta_2 \|_{L^2}+L\|\eta_2 \|_{L^2})
	\\&\quad
	+ 
	|\eta_2(1)||k|^{-1}\|F_k\|_{H^{-1}}(\|\pa_y \chi_{+}^\ast\|_{L^2}
	+
	L\|\chi_{+}^\ast\|_{L^2})
	\\&
	\lesssim  |k|^{-1}\|F_k\|_{H^{-1}}(\|\pa_y \eta_2 \|_{L^2}+L\|\eta_2 \|_{L^2})
	\\&\quad
	+ 
	|\eta_2(1)||k|^{-1}\|F_k\|_{H^{-1}}(L_\ast^{1/2}+L L_\ast^{-1/2}).\label{pairing_31}
\end{split}
\end{align}
For the term $T_2$, we appeal to \eqref{pairing_34}, \eqref{pairing_3/2} and \eqref{F_in_H-1_est} to get
\begin{align}
\begin{split}
	|T_2|
	&\lesssim
	|\eta_{2} (1)|
	\Vert w_{\rm Na}\Vert_{L^1 ([1-L^{-1}_{\ast}, 1])}
	\Vert (\mathcal{U}-\lambda_r) \chi^{\ast}_{+}\Vert_{L^\infty}	
	\\&
	\les
	|\eta_2(1)|L_\ast^{-3/2}\|\pa_y w_{\rm Na}\|_{L^2}(|\mathcal{U}(1)-\lambda_r|+L_{*}^{-1})
	\\&
	\lesssim
	|\eta_2(1)|L_\ast^{-3/2}\nu^{-1}\|F_k\|_{H^{-1}}(|\mathcal{U}(1)-\lambda_r|+L_{*}^{-1}).
	\label{pairing_32}
\end{split}
\end{align}
For the term $T_3$, we use the explicit expression of the Green's function \eqref{green}, yielding
\begin{align}
\begin{split}
	&\|\de_k^{-1}(\mathcal{U}''\chi_{+}^\ast)\|_{L^2}
	\\&
	\les\max_{y\in[-1,1]}\lf|\int_{-1}^1\frac{G_k(y,y')}{1-\max\{y,y'\}} (1-\max\{y,y'\})
	\mathcal{U}''(y')\chi_{+}^\ast(y')dy'\rg|
	\\&
	\les
	\lf\|\frac{G_k(y,y')}{1-\max\{y,y'\}} \rg\|_{L_{y,y'}^\infty}\lf\|(1-\max\{y,y'\})
	\mathcal{U}''(y')\chi_{+}^\ast(y')
	\rg\|_{L_{y}^\infty L_{y'}^1}
	\\&
	\lesssim
	\|\mathcal{U}''\|_{L^\infty}L_\ast^{-1}\|\chi_{+}^\ast\|_{L^1}
	\lesssim  L_\ast^{-2}.
	\label{F_in_H-1_est2}
\end{split}
\end{align}
Combining \eqref{F_in_H-1_est}, \eqref{F_in_H-1_est2}, and the constraint $L_\ast\geq L$, we obtain
\begin{align}
\begin{split}
		|T_3|
		&\lesssim |\eta_2(1)|\|w_{\rm Na}\|_{L^2}L_\ast^{-2}
	\lesssim 
		|\eta_2(1)| \nu^{-2/3}|k|^{-1/3}
	\|F_k\|_{H^{-1}}L_\ast^{-2}
	\\&
	 \lesssim
	 |k|^{-1}|\eta_2(1)|\|F_k\|_{H^{-1}}. \label{pairing_33}
\end{split}
\end{align}
Let $L_\ast:=(|\mathcal{U}(1)-\lambda_r|+L^{-1})^{1/2}L^{3/2}$.
From \eqref{pairing_2}, \eqref{pairing_31},  \eqref{pairing_32}, and \eqref{pairing_33}, we conclude that
\begin{align}
\begin{split}
	\label{pairing_4}
	&\lf|\lf\lan w_{\rm Na}, (\mathcal{U}-\lambda_r)\eta_2-\de_k^{-1}(\mathcal{U}''\eta_2)\rg\ran\rg|
	\\&
	\lesssim 
	|k|^{-1}\|F_k\|_{H^{-1}}(\|\pa_y \eta_2\|_{L^2}
	+L\|\eta_2\|_{L^2})
	\\&\quad
	+|k|^{-1}|\eta_2(1)|\| F_k \|_{H^{-1}}(|\mathcal{U}(1)-\lambda_r|+L^{-1})^{1/4}L^{3/4}.
%	\label{pairing_51}
\end{split}
\end{align}
%Note that \eqref{pairing_4} holds
%for all $\eta_2\in H^1(-1,1)$ with $\eta_2(-1)=0. $

%\step{3} 
Let $g_+\in H^{1}(-1,1)$ with $g_+(-1)=0$. 
Denote
$\eta_3:=g_+\chi$, where $\chi$ is defined in \eqref{defn_chi}.
It is clear that $\eta_3 (-1)=0$ and $\eta_3 \in H^1 (-1,1)$.
%where $\chi$ is defined in \eqref{defn_chi}. 
%\linfeng{shouldn't it be outside of $[\lambda-L^{-1}, \lambda+ L^{-1}]$}
%Note that $(\mathcal{U}-\lambda_r)\eta_3=g_+$ on the interval $[-1,1]\cap [\lambda_r-L^{-1}, \lambda_r+L^{-1}]^c$. 
%Hence we can view the $\eta_3$ as an approximate preimage of the function $ g_+$ with respect to the operator $(\mathcal{U}-\lambda_r)\mathbbm{1}-\de_k^{-1}(\mathcal{U}''\mathbbm{1})$. 
We invoke the estimate \eqref{pairing_4} for $\eta_3$ to get
\begin{align}
\begin{split}
	&|\lan w_{\rm Na}, g_+\ran|
	\\&
	\les
	|\lan w_{\rm Na}, g_+-[(\mathcal{U}-\lambda_r) \eta_3-\de_k^{-1}(\mathcal{U}''\eta_3)]\ran| +
	|\lan w_{\rm Na},(\mathcal{U}-\lambda_r) \eta_3-\de_k^{-1}(\mathcal{U}''\eta_3)\ran |
	\\&
	\les
	\|w_{\rm Na}\|_{L^2}\|g_+-(\mathcal{U}-\lambda_r) \eta_3
	\|_{L^2}
	+
	|\lan w_{\rm Na},\de_k^{-1}(\mathcal{U}''\eta_3)\ran|
	\\&\quad
	+
	|k|^{-1}\|F_k \|_{H^{-1}}(\|\pa_y \eta_3 \|_{L^2}
	+
	L\|\eta_3 \|_{L^2})
	\\&\quad
	+
	|k|^{-1}|\eta_3 (1)|\| F_k \|_{H^{-1}}(|\mathcal{U}(1)-\lambda_r|+L^{-1})^{1/4}L^{3/4}
	\\&
	=
	:T_4+T_5+T_6+T_7 .
	\label{pairing_5}
\end{split}
\end{align}
%\linfeng{here}
We bound the term $T_4$ using \eqref{F_in_H-1_est} as
\begin{align}
\begin{split}
	T_4
	%	&\lesssim 
	%	\|w_{\rm Na}\|_{L^2}\|g_+-(U-\lambda_r)\chi g_+\|_{L^2}
	%	\\&
	&\lesssim \nu^{-2/3}|k|^{-1/3}
	\|F_k \|_{H^{-1}}
	\|g_+\|_{L^\infty([-1,1]\cap [y_\lambda-L^{-1}, y_\lambda+L^{-1}])} L^{-1/2}
	\\&
	\lesssim |k|^{-1}L^{3/2}\|F_k\|_{H^{-1}}\|g_+\|_{L^\infty([-1,1]\cap [y_\lambda-L^{-1}, y_\lambda+L^{-1}])}.
	\label{pairing_61}
\end{split}
\end{align} 
For the term $T_5$,  we appeal to \eqref{F_in_H-1_est} and \eqref{symmetry}, yielding
\begin{align}
\begin{split}
	T_5
	&=
	|\lan \de_k^{-1}w_{\rm Na}, \mathcal{U}'' \eta_3\ran|
	\lesssim 
	\|\phi_{\rm Na}\|_{L^\infty}\|\mathcal{U}''\|_{L^\infty}
	\|g_+\chi \|_{L^1}
	\\&
	\lesssim 
	\|\na_k\phi_{\rm Na}\|_{L^2}\|g_+\chi\|_{L^1}
	\\&
	\lesssim
	\nu^{-1/2}|k|^{-1}
	\|F_k\|_{H^{-1}}\|g_+\chi\|_{L^1}. \label{pairing_62}
\end{split}
\end{align}
Note that the term $T_6$ appears on the right-hand side of \eqref{pairing}.
Finally, for the term $T_7$, we note that
\begin{align*}
\begin{split}
	|\chi(1)|
	&\les \min\{|\mathcal{U}(1)-\lambda_r|^{-1}, L\}
	\les \max\{|\mathcal{U}(1)-\lambda_r|,L^{-1}\}^{-1}
	\\&
	\les (|\mathcal{U}(1)-\lambda_r|+L^{-1})^{-1}.
\end{split}
\end{align*}
Thus, we get
\begin{align}
	T_7
	\lesssim |g_+(1)|\| F_k\|_{H^{-1}}|k|^{-1}(|\mathcal{U}(1)-\lambda_r|+L^{-1})^{-3/4}L^{3/4}.
	\label{pairing_63}
\end{align}
From \eqref{pairing_5},  \eqref{pairing_61}, \eqref{pairing_62}, and \eqref{pairing_63}, we conclude the estimate \eqref{pairing} for $g_{+}$. 
The proof for $g_{-}$ is similar and thus we omit the details.
\end{proof}

\begin{lemma}
\label{lem:cpm_est}
Let $w_{\rm Na}$ be the solution of \eqref{Nav_eqn} with $0\geq k {\rm Im} \lambda \geq -\delta_0 \nu^{1/3} |k|^{2/3}$.
For $\nu k^2\leq \Vert \U'\Vert_{L^\infty}$, we have
\begin{align}
	\label{cest_FL2}
	&
	(1+|k||\lambda_r-\mathcal{U}(1)|)|c_{+}|
	+
	(1+|k||\lambda_r-\mathcal{U}(-1)|)|c_{-}|
	\lesssim
	\nu^{-1/6}|k|^{-5/6}\|F_k\|_{L^2} \log L,
	\\&
	\label{cest_FH1}
	(1+|k||\lambda_r-\mathcal{U}(1)|)|c_{+}|+
	(1+|k||\lambda_r-\mathcal{U}(-1)|)|c_{-}|
	\lesssim
	k^{-1/2}\|\na_k F_k \|_{L^2}\log L,
\end{align}
where $c_{\pm}$ are the coefficients introduced in \eqref{w_I_decmp}.
\end{lemma}

\begin{proof}
	Without loss of generality, we assume that $k>0$.
	Let $\lambda=\lambda_r -i \delta L^{-1}$, where $\delta \in [0,\delta_0]$.
We rewrite \eqref{Nav_eqn}$_1$ as
\begin{align}
	\label{cest_FL2_1}
	(\mathcal{U}-\lambda_r+i L^{-1} )(\pa_{y}^2-k^2)\phi_{\rm Na}
	-
	\mathcal{U}'' \phi_{\rm Na}
	=
	\frac{g_k}{ik}
\end{align} 
where
\begin{align}
	g_k
	=
	\nu\de_k w_{\rm Na}
	+
	F_k
	- 
	k(1-\delta) L^{-1} w_{\rm Na}.
\end{align}
From \eqref{F_in_L2_est} it follows that $\|g_k\|_{L^2}\lesssim \|F_k \|_{L^2}$.
%From \eqref{c_pm_frml} it follows that
%\begin{align}
%\begin{split}
%	-c_{+}&=\int_{-1}^{1}\frac{\sinh(k(y+1))}{\sinh(2k)} (\pa_y^2-k^2)\phi_kdy
%	\\&
%	= \frac{\sinh(k(y+1))}{\sinh(2k)}\pa_y\phi_k\Big|_{y=-1}^{1}-\int_{-1}^1\frac{k\cosh(k(y+1))}{\sinh(2k)}\pa_y\phi dy   
%	\\&\qquad
%	-k^2  \int_{-1}^{1}\frac{\sinh(k(y+1))}{\sinh(2k)}\phi_kdy\\
%	&=\pa_y\phi_k(y=1).
%\end{split}
%\end{align}
From \eqref{c_pm_frml} and \eqref{cest_FL2_1} it follows that
\begin{align}
\begin{split}
	\label{Key_rel}
	|c_{+}|
	&=|\pa_y \phi_{\rm Na}(y=1)|
	\\&
	=\lf|\int_{-1}^1 \frac{1}{\mathcal{U}-\lambda_r
	+
	iL^{-1}}\lf(\mathcal{U}''\phi_k +\frac{ g_k}{ik}\rg)\frac{\sinh(k(y+1))}{\sinh(2k)}  dy\rg|.
\end{split}
\end{align} 
Recall that $\mathcal{U}(y_{\lambda})=\lambda_r$. 
Direct computation gives
\begin{align}
	\begin{split}
		&\lf\|\frac{1}{\mathcal{U}
			-
			\lambda_r+i L^{-1}}\rg\|_{L_y^1}
		= \int_{-1}^1\frac{1}{|\mathcal{U}-\lambda_r+i L^{-1}|}dy
		\\&\quad
		\lesssim 
		\int_{[-1,1]\cap [y_{\lambda}-L^{-1},y_{\lambda}+L^{-1}]}\frac{1}{ L^{-1}}dy 
		+
		\int_{[-1,1]\backslash[y_{\lambda}-L^{-1},y_{\lambda}+L^{-1}]}\frac{1}{|\mathcal{U}-\lambda_r|}dy\\
		&\quad
		\lesssim 
		\log L
		\label{EQ120a}
	\end{split}
\end{align}
and
\begin{align}
	\begin{split}
		&\lf\|\frac{1}{\mathcal{U}-\lambda_r+i L^{-1}}\rg\|_{L_y^2}^{ 2}
		=
		\int_{-1}^1\frac{1}{|\mathcal{U}-\lambda_r|^2+L^{-2}}dy 
		\\&
		\quad
		\lesssim \int_{[-1,1]\cap [y_{\lambda}-L^{-1},y_{\lambda}+L^{-1}]}\frac{1}{ L^{-2}}\,dy +\int_{[-1,1]\backslash[y_{\lambda}-L^{-1},y_{\lambda}+L^{-1}]}\frac{1}{|\mathcal{U}-\lambda_r|^2}dy
		\\&
		\quad 
		\lesssim 
		L.
	\end{split}
	\label{cest_FL2_2}
\end{align}

$\bullet$ \textbf{The case $|\lambda_r-\mathcal{U}(1)|\leq |k|^{-1}$.}
From \eqref{Key_rel} it follows that
\begin{align}
\begin{split}
	|c_+|
	&
	\lesssim 
	\|\phi_k\|_{L^\infty}\lf\|\frac{1}{\mathcal{U}-\lambda_r+i L^{-1}}\rg\|_{L_y^1}
	+
	|k|^{-1}\|g_k\|_{L^2}\lf\|\frac{1}{\mathcal{U}-\lambda_r+i L^{-1}}\rg\|_{L_y^2}.
	\label{EQ151a} 
\end{split}
\end{align}
Combining \eqref{EQ120a}, \eqref{cest_FL2_2}, \eqref{EQ151a}, and \eqref{F_in_L2_est}, we deduce that
\begin{align}
\begin{split}
	|c_{+}|
	&\lesssim \|\na_k\phi_k\|_{L^2}\log L
	+
	|k|^{-1}\|F_k\|_{L^2}\nu^{-1/6}|k|^{1/6}
%	\\&
%	\lesssim 
%	|\log L|
%	\nu^{-1/6}|k|^{-4/3}
%	\|F_k\|_{L^2}+|k|^{-5/6}\nu^{-1/6}
%	\|\mathbb{F}\|_{L^2}
	\\&
	\lesssim 
	\nu^{-1/6} |k|^{-5/6} \log L
	\|F_k \|_{L^2}.
\end{split}
\end{align}

Suppose $F_k \in H^1_0$.
%
%
%
%\textit{The case $|\lambda_r-\mathcal{U}(1)|\leq |k|^{-1}$ and $\mathbb{F}\in H_0^1$.}
We integrate by parts of \eqref{Key_rel} to obtain
\begin{align}
\begin{split}
	|c_{+}|
	&=
	\Big|\int_{-1}^1 \frac{1}{\mathcal{U}-\lambda_r+iL^{-1}}\lf(\mathcal{U}''\phi_k 
	+
	\frac{ \nu\de_kw_{\rm Na}+F_k -k(1-\delta) L^{-1}w_{\rm Na}}{ik}\rg)
	\\&\quad\times
	\frac{\sinh(k(y+1))}{\sinh(2k)}  dy\Big|\\
	&\les 
	\|\phi_k\|_{L^\infty}\lf\|\frac{1}{\mathcal{U}-\lambda_r+iL^{-1}}\rg\|_{L_y^1}
%	\\&\quad
	+
	|k|^{-1}\nu\|\de_k w_{\rm Na}\|_{L^2}\lf\|\frac{1}{\mathcal{U}-\lambda_r+iL^{-1}}\rg\|_{L_y^2}
	\\&\quad
	+
	\lf|\int_{-1}^1 \log\lf(\mathcal{U}-\lambda_r+iL^{-1}\rg) \pa_y\lf(\frac{F_k}{\mathcal{U}' ik}\frac{\sinh(k(y+1))}{\sinh(2k)} \rg) dy\rg|
	\\&\quad
	+
	L^{-1}
	\Vert w_{\rm Na}\Vert_{L^2}\lf\|\frac{1}{\mathcal{U}-\lambda_r+iL^{-1}}\rg\|_{L_y^2}.
\end{split}
\label{EQ130a}
\end{align}
By the complex natural logarithmic rule we have
\begin{align}
\label{log_bnd}
	&\lf|\log(\mathcal{U}
	-
	\lambda_r
	+
	iL^{-1})\rg|
	\lesssim 
	\log L,
\end{align}
for all $y\in[-1,1]$ and $\lambda_r\in[\mathcal{U}(1)-|k|^{-1}, \mathcal{U}(1)+|k|^{-1}]$.
Using \eqref{F_in_H1_est}, \eqref{EQ120a}, \eqref{cest_FL2_2}, \eqref{EQ130a}, \eqref{log_bnd}, we obtain
\begin{align} 
	|c_{+}|
	&
	\les 	
	|k|^{-1/2}\|\na_k F_k \|_{L^2}
	\log L.
\end{align}
%completing the proof of \eqref{cest_FH1} for the case $|\lambda_r-\mathcal{U}(1)|\leq |k|^{-1}$ and $\mathbb{F}\in H_0^1$.
%\linfeng{here}

$\bullet$ \textbf{The case  $\lambda_r\geq \U(1)+|k|^{-1}$}. 
Using \eqref{F_in_L2_est} and
\begin{align}
	\frac{1+|k||\mathcal{U}(1)-\lambda_r|}{|\mathcal{U}(y)-\lambda_r|+
	L^{-1}}
	\lesssim 
	|k|,
	\quad
	y\in[-1,1],
\end{align}
we obtain
\begin{align}
\begin{split}
	&(1+|k||\mathcal{U}(1)-\lambda_r|)|c_{+}|
	\\&\quad
	\lesssim 
	\lf|\int_{-1}^1
	\frac{ 1+|k||\mathcal{U}(1)-\lambda_r|}{\mathcal{U}-\lambda_r+iL^{-1}}
	\lf(\mathcal{U}''\phi_k +\frac{ g_k}{ik}\rg)\frac{\sinh(k(y+1))}{\sinh(2k)}  dy \rg|
	\\&\quad
	\lesssim
	\|\na_k\phi_k\|_{L^2}
	\lf\|\frac{\sinh(k(y+1))}{\sinh(2k)} \rg\|_{L^1}
	\\&\qquad \qquad
	+
	\|g_k\|_{L^2} (|k|^{-1}+|\mathcal{U}(1)-\lambda_r|)\lf|\int_{-1}^1\lf(\frac{\sinh(k(y+1))}{(\mathcal{U}-\lambda_r)\sinh(2k)}\rg)^2dy \rg|^{1/2}
	\\&\quad
%	\lesssim \nu^{-1/6}|k|^{-4/3}\|F_k \|_{L^2}+
%	\nu^{-1/6}|k|^{-5/6}\|F_k \|_{L^2}
%	\\&\quad
	\lesssim
	\nu^{-1/6}|k|^{-5/6}\|F_k\|_{L^2}.
	\label{EQ152a}
\end{split}
\end{align}

Suppose $F_k \in H^1_0$. We proceed as in \eqref{EQ152a} and use \eqref{F_in_H1_est} instead of \eqref{F_in_L2_est}, yielding
\begin{align}
	(1+|k||\mathcal{U}(1)-\lambda_r|)
	|c_{+}|
	&
	\lesssim 
	|k|^{-3/2}\|\mathbb \nabla_k F_k
	\|_{L^2}.
\end{align}
%For Lemma 4.1 of \cite{CLWZ20}, we can use the above formulation to treat 

$\bullet$ \textbf{The case $\lambda_r<\U(1)-|k|^{-1}$.}
Let $\wt{y_\lambda}\in\rr$ be such that $\wt{\U} (\wt{y_\lambda}) = (1+\lambda)/2$.
Let $E=[-1,1]\cap [-1, \wt{y_\lambda})$ and $E^c=[-1,1]\cap [\wt{y_\lambda},1]$.
% For $y\in E^c$, it is clear that
%\begin{align}
%	\frac{1}{|\U(y)-\lambda_r|}
%	\lesssim \frac{1}{|\U(1)-\lambda_r|+|k|^{-1}}.
%\end{align}
Using \eqref{Key_rel}, we have
\begin{align}
	|c_{+}|
	=:\int_{-1}^1 
	\mathcal{I} \,dy
	=
	\int_{E^c} 	\mathcal{I} \,dy
	+
		\int_{E} 	\mathcal{I} \,dy.
		\label{EQ120b}
\end{align}
For the first term on the right-hand side of \eqref{EQ120b}, we use \eqref{F_in_L2_est} and the bound
\begin{align}
	\lf(\int_{E^c}
	\left|
	\frac{\sinh(k(y+1))}{(\U-\lambda_r)\sinh(2k)}
	\right|^2dy \rg)^{1/2} \lesssim |k|^{-1/2}|\lambda_r-\U(1)|^{-1}
\end{align}
to get
\begin{align}
	\begin{split}
	\int_{E^c} 	\mathcal{I} \,dy
	&\lesssim \nu^{-1/6}|k|^{-11/6}\|F_k\|_{L^2}(|\U(1)-\lambda_r|+|k|^{-1})^{-1}
	\\&\quad
	+ |k|^{-3/2}\|F_k\|_{L^2}(|\U(1)-\lambda_r|+|k|^{-1})^{-1}
	\\&
	\lesssim 
	\nu^{-1/6}|k|^{-5/6}
	\| F_k \|_{L^2}(|k| |\U(1)-\lambda_r|+{1})^{-1}.
	\end{split}
\label{EQ130c}
\end{align}
For the second term on the right-hand side of \eqref{EQ120b}, we split into two cases.
If $\lambda_r \leq -1-|k|^{-1}$, we have
\begin{align}
	\begin{split}
		&\int_{E} 	\mathcal{I} \,dy
	\les	
	\| \U'' \phi_k+(ik)^{-1}g_k\|_{L^2(E)}
	\lf\|\frac{\sinh(k(y+1))}{\sinh(2k)}\rg\|_{L^\infty(E)}
	(|\U(-1)-\lambda_r|+L^{-1})^{-1}
	\\&
%	\les
%	\nu^{-1/6} |k|^{-4/3}
%	e^{-|k|/4}(|\U(-1)-\lambda_r|+\delta_0L^{-1})^{-1}
%	\Vert F_k \Vert_{L^2}
%	\\&
	\lesssim
	\nu^{-1/6} |k|^{-4/3} e^{-|k|/4}(|\U(-1)-\lambda_r|+2|k|^{-1})^{-1}
	\Vert F_k \Vert_{L^2}
	\\&
	\lesssim
		\nu^{-1/6} |k|^{-4/3} e^{-|k|/4}(|\U(-1)-\lambda_r|+|k|^{-1}|\U(1)-\U(-1)|+|k|^{-1})^{-1}
			\Vert F_k \Vert_{L^2}
	\\&
	\lesssim 
		\nu^{-1/6} |k|^{-1/3}
	e^{-|k|/4}(|\U(-1)-\lambda_r|
	+ 
	|\U(1)-\U(-1)|+ 1)^{-1}
	\Vert
	F_k \Vert_{L^2}
	\\&
	\lesssim (|k||\U(1)-\lambda_r|+ 1)^{-1}
	\nu^{-1/6}|k|^{-5/6}\|F_k \|_{L^2}.
	\end{split}
\label{EQ130d}
\end{align}
If $-1-|k|^{-1} < \lambda_r < \U(1)-|k|^{-1}$, we appeal to \eqref{EQ120a} and \eqref{cest_FL2_2} to obtain
\begin{align}
\begin{split}
		\int_{E} 	\mathcal{I} \,dy
	&\lesssim 
	\|\phi_k\|_{L^\infty}
	\lf\|\frac{1}{|\U-\lambda_r|
	+
	L^{-1}}\rg\|_{L^1(E)}
	\lf\|\frac{\sinh(k(y+1))}{\sinh(2k)}\rg\|_{L^\infty(E)}
	\\&\quad
	+
%	\mathbbm{1}_{\lambda_r\in[ -1-|k|^{-1}, 1-|k|^{-1}]}
	k^{-1}
	\|g_k\|_{L^2(E)}
	\lf\|\frac{\sinh(k(y+1))}{\sinh(2k)}\rg\|_{L^\infty(E)}
	\lf\|
		\frac{1}{| \U-\lambda_r|+ L^{-1}}\rg\|_{L^2 (E)}
	\\&
	\les
	\nu^{-1/6}|k|^{-5/6}\|F_k\|_{L^2}(|k||\U(1)-\lambda_r|+{1})^{-1}
	\log L.
\end{split}
\label{EQ131a}
\end{align}
%Now we use the estimate \eqref{F_in_L2_est} to obtain that 
%\begin{align}
%\|\na_k\phi_k\|_{L^2}\leq \nu^{-1/6}|k|^{-4/3}\|\mathbb{F}\|_{L^2},\quad \|(ik)^{-1}g_k\|_{L^2} \lesssim |k|^{-1}\|\mathbb{F}\|_{L^2}.
%\end{align}  
%First we estimate the term $T_1$.
%It is readily checked that
%\begin{align}
%	\lf(\int_{E^c}
%	\left|
%	\frac{\sinh(k(y+1))}{(\U-\lambda_r)\sinh(2k)}
%	\right|^2dy \rg)^{1/2} \lesssim |k|^{-1/2}|\lambda_r-\U(1)|^{-1}.
%\end{align}
%Since $\nu|k|^2\lesssim1$, we use \eqref{F_in_L2_est} to get
%\begin{align}
%\begin{split}
%	T_1
%	\lesssim& \nu^{-1/6}|k|^{-11/6}\|\mathbb{F}\|_{L^2}(|\U(1)-\lambda_r|+|k|^{-1})^{-1}
%	\\&\qquad
%	+ |k|^{-3/2}\|\mathbb{F}\|_{L^2}(|\U(1)-\lambda_r|+|k|^{-1})^{-1}
%	\\
%	\lesssim &\nu^{-1/6}|k|^{-5/6}\|\mathbb{F}\|_{L^2}(|k| |\U(1)-\lambda_r|+{1})^{-1}.
%\end{split}
%\end{align}
%For the term $T_2$, we use \eqref{F_in_L2_est} to get
%\begin{align}
%	\begin{split}
%	T_2
%	&\lesssim
%	\nu^{-1/6} k^{-4/3} \|\mathbb{F}\|_{L^2}
%	|\log L|	
%	e^{-|k|/4}
%	\\&
%	\lesssim
%	\nu^{-1/6} k^{-5/6}\|\mathbb{F}\|_{L^2}(|k||\U(1)-\lambda_r|+{1})^{-1}
%	|\log L|.
%	\end{split}
%\end{align}
%Using \eqref{EQ120a} and \eqref{cest_FL2_2}, we obtain
%\begin{align}
%	T_2
%	+T_3
%	\lesssim
%	\nu^{-1/6}|k|^{-5/6}\|\mathbb{F}\|_{L^2}(|k||\U(1)-\lambda_r|+{1})^{-1}
%	|\log L|.
%\end{align}
Combining \eqref{EQ120b}, \eqref{EQ130c}, \eqref{EQ130d}, and \eqref{EQ131a}, we conclude that
\begin{align}
	|c_{+}|
	\lesssim \nu^{-1/6}|k|^{-5/6}
	\|F_k \|_{L^2}(|k||\lambda_r-\U(1)|+L^{-1})^{-1}
	\log L.
	\label{EQ121a}
\end{align}

For the case when $F_k \in H^1_0$, we have
\begin{align}
\begin{split}
	|c_{+}|
	&
	\leq 
	\Big|
	\int_{-1}^1 \frac{1}{\U-\lambda_r+iL^{-1}
	}
	\lf(\U''\phi_k 
	+
	\frac{\nu\de_kw_{\rm Na}- k(1-\delta) L^{-1} w_{\rm Na}}{ik} \rg)
	\\&\quad\quad
	\times
	\frac{\sinh(k(y+1))}{\sinh(2k)}  dy\Big| 
	\\&\quad
%	+
%	\lf|\int_{-1}^1 \frac{1}{\U-\lambda_r+iL^{-1}
%	}
%	\lf(\U''\phi_k +\frac{ \nu\de_kw_{\rm Na}}{ik} \rg)\frac{\sinh(k(y+1))}{\sinh(2k)}  dy\rg|
%	\\&\quad
	+ 
	\lf| \int_{-1}^1 \frac{1}{\U-\lambda_r+iL^{-1}}  \frac{F_k }{ik}\frac{\sinh(k(y+1))}{\sinh(2k)}  dy\rg|\\
	&=:
	\mathcal{I}_1+\mathcal{I}_2.
\end{split}
\label{EQ132d}
\end{align} 
For the term $\mathcal{I}_1$, we proceed in a similar fashion as in \eqref{EQ120b}--\eqref{EQ121a} and use \eqref{F_in_H1_est} to obtain
\begin{align}
	\mathcal{I}_1
	\lesssim k^{-1/2}\|\na_k F_k \|_{L^2}
	(k|\lambda_r-\U(1)|+L^{-1})^{-1}
	\log L.%23
	\label{EQ132a}
\end{align} 
For the term $\mathcal{I}_2$, we split into two cases.
If $\lambda_r \leq \U(-1)-k^{-1}$, we use the bound
\begin{align}
	\max_{y\in[-1,1]}\lf|
	\frac{1+k| \U(1)-\lambda_r|}{ \U (y)-\lambda_r+iL^{-1}}\rg|
	\lesssim 
	k^2
\end{align} 
to get
\begin{align}
\begin{split}
	\mathcal{I}_2
	&\lesssim 
	\frac{k^2}{	(1+k|\mathcal{U}(1)-\lambda_r|)}
	\lf| \int_{-1}^1  \frac{F_k}{ik }\frac{\sinh(k(y+1))}{\sinh(2k)}  dy\rg|
	\\&
	\lesssim 
	k^{-1/2}
	(1+k|\mathcal{U}(1)-\lambda_r|)^{-1}
	\|\na_k F_k\|_{L^2}. 
\end{split}
	\label{EQ132b}
\end{align}
If $\U(-1)-k^{-1}\leq \lambda_r \leq \U(1)-k^{-1}$, we use the bound $(1+|k||\mathcal{U}(1)-\lambda_r|)\lesssim |k|$ and \eqref{log_bnd} to obtain
\begin{align}
\begin{split}
	\mathcal{I}_2
	&\lesssim
	 \frac{k}{(1+k|\mathcal{U}(1)-\lambda_r|)}
	\lf| \int_{-1}^1 \log\lf(\mathcal{U}-\lambda_r+i\delta_0 L^{-1}\rg) \pa_y\lf( \frac{F_k }{ik\mathcal{U}' }\frac{\sinh(k(y+1))}{\sinh(2k)} \rg) dy\rg|
	\\&
	\lesssim 
	k^{-1/2}
	(1+k|\mathcal{U}(1)-\lambda_r|)^{-1}
	\|\na_k F_k \|_{L^2}\log L.
\end{split}
\label{EQ132c}
\end{align}
Combining \eqref{EQ132d}, \eqref{EQ132a}, \eqref{EQ132b}, and \eqref{EQ132c}, we conclude
\begin{align}
	\begin{split}
	|c_{+}|	
	\lesssim 
	|k|^{-1/2}
	\|\nabla_k F_k \|_{L^2}(|k||\lambda_r-\U(1)|+L^{-1})^{-1}
	\log L.
	\end{split}
\end{align}

The estimates for $|c_{-}|$ in \eqref{cest_FL2} and \eqref{cest_FH1} are similar and thus we omit the further details.
\end{proof}

%The following lemma is an analog of \cite[Lemma~4.2]{CLWZ20}.
%Lemma 4.2 in \cite{CLWZ20}, more work is needed.
\begin{lemma}\label{lem:cpm_est_2}
	Let $w_{\rm Na}$ be the solution of \eqref{Nav_eqn} with $0\geq k {\rm Im} \lambda \geq -\delta_0 \nu^{1/3} |k|^{2/3}$.
For $F_k \in H^{-1}$, we have
\begin{align}
\begin{split}
		&
	\label{cest_F_H-1}
	(1+|k| |\lambda_r-\U(1)|)^{3/4}|c_{+}|
	+
	(1+|k||\lambda_r-\U(-1)|)^{3/4}|c_{-}|
	\\&\qquad
	\lesssim
	\nu^{-1/2} |k|^{-1/2}
	\|F_k\|_{H^{-1}}
	\log L.
\end{split}
\end{align}
\end{lemma} 

\begin{proof}
	Let
\begin{align}
	g_+ (y)=\frac{\sinh(k(y+1))}{\sinh(2k)},
	\quad
	g_- (y)=\frac{\sinh(k(1-y))}{\sinh(2k)}.
\end{align}
It is clear that $g_{+} \in H^1 (-1,1)$ and $g_{+} (-1)=0$.
Using \eqref{c_pm_frml} and Lemma~\ref{L01}, we get
\begin{align}
\begin{split}
	|c_{+}|&=|\lan w_{\rm Na}, g_+\ran |
	\lesssim |k|^{-1}\|F_k \|_{H^{-1}}\Big(L^{3/2}\|g_+ \|_{L^\infty ((-1,1)\cap (y_r-L^{-1}, y_r+L^{-1}))}
		\\&\quad
		+
		(|\U( 1)-\lambda_r|+L^{-1})^{-3/4}L^{3/4}
	+\|g_+ \chi\|_{H^1}
%	\\&\quad
	+L\|g_+ \chi\|_{L^2}
	+
	\nu^{-1/2}\|g_{+}\chi\|_{L^1}\Big),
	\label{EQ121b}
\end{split} 
\end{align} 
where $\chi$ is a smooth cutoff function defined in \eqref{defn_chi}.
To estimate the last term on the right side of \eqref{EQ121b}, 
we aim to prove
\begin{align}
	\label{cest_F_H-1_1}
	(1+|k| |\lambda_r-\U(1)|)^{3/4}\|g_{+}\chi\|_{L^1}
	\lesssim  
	|k|^{ 1/2}\log L,
\end{align}
which concludes that the $L^1$-contribution on the right-side of \eqref{EQ121b}.

$\bullet$ \textbf{The case $\lambda_r-\U(1)\geq 2k^{-1}$.} 
Using the mean value theorem, we obtain $y_\lambda -1 \geq k^{-1}\geq L^{-1}$. 
Hence, we deduce that
\begin{align}
	\|g_+\chi\|_{L^1}
	\lesssim 
	\Vert \chi\Vert_{L^\infty}
	\|g_+\|_{L^1}
	\lesssim \frac{1}{|k|^{-1}+|\U(1)-\lambda_r|}
	\frac{1}{|k|}
	=
	\frac{1}{1+ |k| |\U (1)-\lambda_r|}.
\end{align}

$\bullet$
\textbf{The case $|\lambda_r-\U(1)|\leq 2k^{-1}$.} 
We use \eqref{EQ120a} to get
\begin{align}
	\|g_+\chi\|_{L^1}
	\lesssim \|g_+\|_{L^\infty}
	\Vert \chi\Vert_{L^1}
	\lesssim 
	\int_{-1}^1\frac{1}{|\U-\lambda_r|+L^{-1}}dy
	\les
	\log L. 
\end{align}

$\bullet$
\textbf{The case $\lambda_r\in[-1-2|k|^{-1}, \U(1)-2|k|^{-1}]$.} 
Let $\wt{y_\lambda}= \wt{\U}^{-1} ((1+\lambda_r)/2)$.
Let $E=[-1,1]\cap [-1, \wt{y_\lambda})$ and $E^c=[-1,1]\cap[ \wt{y_\lambda}, 1]$.
It is readily checked that
\begin{align}
	\Vert \chi\Vert_{L^\infty (E^c)}
	\lesssim \lf(|k|^{-1}+|\U(1)-\lambda_r|\rg)^{-1}. 
\end{align}
Using \eqref{EQ120a}, we obtain
\begin{align}
\begin{split}
	\|g_+\chi\|_{L^1}
	&
	\leq \Vert g_+ \Vert_{L^\infty (E)} 
	\Vert
	\chi\|_{L^1(E)}
	+
	\Vert g_+\Vert_{L^1 (E^c)}
	\Vert \chi\Vert_{L^\infty(E^c)}
%	\\&
%	\lesssim 
%	\lf\|g_{+}\rg\|_{L^\infty (E)}
%	\int_{-1}^1\frac{1}{|\U-\lambda_r|+L^{-1}}dy 
%	+
%	\|g_+\|_{L^1(E^c)}
%	\lf(k^{-1}+|\U(1)-\lambda_r|\rg)^{-1}
%	\\&
%	\lesssim 
%	e^{-k(1-\wt{y_\lambda})/2}
%	\log L +\frac{1}{1+|k||\U(1)-\lambda_r|} 
	\\&
	\lesssim e^{-|k|(1-\wt{y_\lambda})}\log L 
	+
	|k|^{-1}
	\lf(|k|^{-1}+|\U(1)-\lambda_r|\rg)^{-1}
	\\&
	\lesssim \lf(1+|k||\U(1)-\lambda_r|\rg)^{-1}
	\log L.
\end{split}
\end{align}

$\bullet$
\textbf{The case $\lambda_r\leq -1-2 |k|^{-1}$.}
It is clear that 
\begin{align}
	|k||\U(1)-\lambda|
	=
	|k||\U(1)-\U(-1)|
	+
	|k||\U(-1)-\lambda_r|
	\lesssim 
	k^2|\U(-1)-\lambda_r|.
\end{align}
Therefore,
\begin{align}
\begin{split}
	\|g_+\chi\|_{L^1}
	&\leq \Vert g_+ \Vert_{L^1 (E)}
	\Vert	\chi\Vert_{L^\infty (E)}
	+\Vert g_+ \Vert_{L^1 (E^c)}
	\Vert \chi\Vert_{L^\infty (E^c)}
	\\&
	\lesssim \frac{e^{-|k|/2}}{|\U(-1)-\lambda_r|+|k|^{-1}}
	+
	\frac{1}{|k|}
	\frac{1}{|\U(-1)-\lambda_r|}
	\\&
	\lesssim \frac{e^{-|k|/4}}{1+|k||\U(1)-\lambda_r|}
	+
	\frac{1}{1+|k||\U(1)-\lambda_r|}
	\\&
	\lesssim \frac{1}{1+|k||\U(1)-\lambda_r|}.	
\end{split}
\end{align}
This concludes the proof of \eqref{cest_F_H-1_1}. 
The rest of the terms on the right side of \eqref{EQ121b} are estimated in a similar fashion as in \cite[Lemma~4.2]{CLWZ20}.
Since the estimates for $c_{-}$ is similar, we conclude the proof of the lemma.
\end{proof}

\subsection{Boundary Corrector}
In this section, we provide the estimates for the boundary corrector. Most of the estimates are carried out in the works \cite{CLWZ20,CWZ23}. For the sake of completeness, we decide to review their main results and key estimates.  

Recall that the elementary boundary correctors introduced in \eqref{w_I_decmp} satisfy
\begin{align} 
	\label{bdy_crt2}
	\left\{
	\begin{aligned}
		&    -\nu(\pa_y^2-k^2)w_{\pm;k}+ik\lf(\mathcal{U}-\lambda\rg)w_{\pm;k}-\mc U''ik\phi_{\pm;k}	=0,\\
		&\phi_{\pm ;k} = \de_{k}^{-1} w_{\pm;k},\quad \pa_y\phi_{\pm ;k}(y=\pm 1)=1, \quad \pa_y\phi_{\pm ;k}(y=\mp 1)=0.  
	\end{aligned}
	\right.
\end{align} 

Let $Ai (y)$ be the Airy function, which is a nontrivial solution of $f''-yf=0$.
%We recall the homogeneous solutions 
%\begin{align}
%H_-(z)=Ai(e^{i\pi/6}z), \quad H_+(z)=Ai(e^{i5\pi/6}z),
%\end{align}
%are two linearly independent solutions of the Airy equation  $u''(z)-iz u(z)=0$.
For convenience, we denote
\begin{align}
	L_\pm:=
	\lf(\frac{\nu}{\mc U'(\pm 1)k}\rg)^{-1/3}, 
	\quad d_\pm := \frac{\mc U(\pm 1)\mp \mc U'(\pm 1)-\lambda-i\nu k}{\mc U'(\pm 1)}.
\end{align}
and define the approximate boundary correctors
\begin{align}
	W_{\pm;\rma}(y)
	=
	Ai
	\left(
	L_{\pm} (y+d_{\pm})
	\exp 
	\left(
	i \lf(\frac{\pi}{2} \pm \frac{\pi}{3}\rg)
	\right)
	\right).
\end{align}
%
% $(W_{+;\rma}, W_{-;\rma})$ as
%\begin{align}
%	W_{-;\rma} (y)
%	=
%	Ai(e^{i\pi/6}L_-(y+d_-)),
%	\qquad 
%	W_{+;\rma}(y)
%	=
%	Ai(e^{i5\pi/6}L_+(y+d_+)).
%\end{align}
%\linfengnote{Here, `$\rma$' stands for the main Airy approximation component and `$\rme$' stands for the remainder. }
It is easy to check that
\begin{align}
	L_{\pm}^{-3} \pa_y^2 W_{\pm;\rma}
	-
	i
	\lf(y+d_{\pm}
	\rg)
	W_{\pm;\rma}=0,
\end{align}
which is equivalent to
\begin{align} 
	\label{bdy_crt_A} -\nu(\pa_y^2-k^2)
	W_{\pm;\rma}
	+
	ik
	\lf(
	\mc U(\pm 1)
	+
	\mc U'(\pm 1)(y\mp 1)-\lambda 
	\rg) 
	W_{\pm;\rma}=0.
\end{align}
Let $W_{\pm;\rme}$ be the solution of
\begin{align}
\label{bdy_crt_E}
\left\{
\begin{aligned}
	&  
	  -\nu(\pa_y^2-k^2)
	  W_{\pm;\rme}+
	  ik\lf(\mc U-\lambda\rg)
	  W_{\pm;\rme}
	  -
	  \mc U''ik
	  \Phi_{\pm;\rme}	
	  \\&\qquad
	  =
	  \mc U''ik\Phi_{\pm; \rma}-ik[\mc U(y)-(\mc U(\pm 1)
	+
	\mc U'(\pm 1)(y \mp 1 )]W_{\pm;\rma},
	\\&
	\Phi_{\pm ;\rma}=\de_{k}^{-1}W_{\pm;\rma},
	\quad
	\Phi_{\pm ;\rme}=\de_{k}^{-1}W_{\pm;\rme},
	\quad \quad 
	W_{\pm ;\rme}
	\big|_{y=\pm 1}=0.
	\end{aligned}  
\right.
\end{align}
Here, the expression $\mc U(\pm 1)+\mc U'(\pm 1)(y\mp 1)$ is the linear approximation of the strictly monotone shear flow on the upper and lower boundaries. 
%The $(\de_{D;k})^{-1}$ is the inverse of the Dirichlet Laplacian $\pa_y^2-k^2$ on the channel $y\in[-1,1]$.   
Define
\begin{align}
	&
	W_{\pm} =W_{\pm;\rma}
	+
	W_{\pm; \rme},
\quad
	\Phi_{\pm}= \Delta_{k}^{-1} W_{\pm}.
	\quad
\end{align}
It is clear that
\begin{align}
	-\nu (\partial_y^2 - k^2)
	W_{\pm} 
	+ik
	(\mc U(y) - \lambda)W_{\pm} 
	- \mc U''(y) ik \Phi_{\pm}
	=0.
\end{align}
%Now we can follow their argument \cite{CWZ23} and derive similar estimates. 
%\myr{(Check!?)}
%\linfeng{This section is quite similar to CWZ section 4.
%}
%We denote
%\begin{align}
%	W_1=W_{a,1}+W_{e,1}
%\end{align}
We assume that the boundary correctors solution $(w_{+}, w_{-})$ to \eqref{bdy_crt2} can be written as 
\begin{align}
	w_{+}
	=
	\mathfrak{C}^{-}_{+}
	W_{-} (y)
	+ 
	\mathfrak{C}^{+}_{+}
	W_{+} (y)
	,\quad
	w_{-}=
	\mathfrak{C}_{-}^{-}
	W_{-} (y) +
	\mathfrak{C}^{+}_{-} W_{+} (y).
	\label{EQ101a}
\end{align}
It is easy to see that
\begin{align}
	&\int_{-1}^{1}\frac{\sinh(k(y\pm 1))}{\sinh(2k)} w_{+} dy
	=\int_{-1}^{1}\frac{\sinh(k(y\pm1))}{\sinh(2k)} (\pa_y^2-k^2)\phi_{+} dy
	=\pa_y \phi_{+} (y=\pm 1)
\end{align}
and
\begin{align}
	&\int_{-1}^{1}\frac{\sinh(k(y\pm 1))}{\sinh(2k)} w_{-} dy
	=\int_{-1}^{1}\frac{\sinh(k(y\pm1))}{\sinh(2k)} (\pa_y^2-k^2)\phi_{-} dy
	=\pa_y \phi_{-} (y=\pm 1).
\end{align}
Therefore, the boundary conditions on $\partial_{y} \phi_{\pm}$ in \eqref{bdy_crt2} are satisfied by setting 
\begin{align}
	\int_{-1}^{1}
	w_{+} (y)
	\frac{\sinh (k(y+1))}{\sinh (2k)}\,dy
	=1
	\quad
	\text{and}
	\quad
	\int_{-1}^{1}
	w_{+} (y)
	\frac{\sinh (k(y-1))}{\sinh (2k)}\,dy
	=0
\end{align}
and
\begin{align}
	\text{and}
	\quad
	\int_{-1}^{1}
	w_{-} (y)
	\frac{\sinh (k(y+1))}{\sinh (2k)}\,dy
	=0,
	\quad
		\int_{-1}^{1}
	w_{-} (y)
	\frac{\sinh (k(y-1))}{\sinh (2k)}\,dy
	=1.
\end{align}
Using \eqref{EQ101a}, we obtain
\begin{align}
\begin{split}
		&\mathfrak C_{+}^{-}
	\underbrace{	\int_{-1}^{1}
	W_{-} (y)
	\frac{\sinh (k(y-1))}{\sinh (2k)}
	\,dy}_{=:A_{--}}
	+\mathfrak C_{+}^{+}
	\underbrace{\int_{-1}^{1}
	W_{+}(y)
	\frac{\sinh (k(y-1))}{\sinh (2k)}\,dy}_{=:A_{+-}}=0,
	\\
	&\mathfrak C_{+}^{-}
	\underbrace{\int_{-1}^{1}
	W_{-} (y)
	\frac{\sinh (k(y+1))}{\sinh (2k)}
	\,dy}_{=:A_{-+}}
	+\mathfrak C_{+}^{+}
	\underbrace{\int_{-1}^{1}
	W_{+}(y)
	\frac{\sinh (k(y+1))}{\sinh (2k)}\,dy}_{=:A_{++}}=1.
\end{split}
\label{EQ101}
\end{align}
Similarly, we have
\begin{align}
	\mathfrak{C}^-_-
	A_{--}
	+
	\mathfrak{C}^+_-
	A_{+-}
	=1,
	\quad
	\mathfrak{C}_-^-A_{-+}
	+
	\mathfrak{C}^+_- A_{++}
	=0
\end{align}
If $\mathfrak{D}:=A_{--}A_{++}-A_{+-} A_{-+}\neq 0$, then $w_{+}$ and $w_{-}$ are linearly independent and
\begin{align}
	(\mathfrak C_{+}^-,\mathfrak C_{+}^+)=
	\frac{(-A_{+-}, A_{--})}{\mathfrak D},
	\quad(\mathfrak C_{-}^-, \mathfrak C_{-}^+)
	=
	\frac{(A_{++}, -A_{-+})}{\mathfrak D}.\label{C_coeff_exp}
\end{align}
We shall prove that $\mathfrak{D} \neq 0$ after Lemma~\ref{lem:4.4CWZ23} in \eqref{Evan's_func_bnd}.
To ease the notation in forthcoming arguments, we further denote
\begin{align}
	\label{A_ap_exp}
	A_{s\pm;\rma}:=\int_{-1}^{1}
	W_{s;\rma}(y) \frac{\sinh (k (y\pm 1))}{\sinh (2k)}
	dy,\quad A_{s\pm;\rme}:=\int_{-1}^{1}
	W_{s;\rme}(y) \frac{\sinh (k (y\pm 1))}{\sinh (2k)}
	dy,
\end{align}
where $s\in\{+,-\}$ and 
\begin{align}
	A_0 (z) 
	=
	\int_{e^{i\pi/6}z}^{\infty}
	Ai (t)dt
	= 
	e^{i\pi/6}
	\int_{z}^\infty
	Ai (e^{i\pi/6} t)dt.
\end{align}
It is easy to see that
\begin{align}
	A_0'(L_{-} (y+d_{-}))
	=-e^{i\pi/6} Ai (L_{-} (y+d_{-}) e^{i\pi/6})
	=-e^{i\pi/6}
	W_{-;\rma}
	(y).
	\label{EQ133a}
\end{align}
%and
%\begin{align}
%	A_0'(L_{-} (y+d_{-}))
%	=-e^{i\pi/6} A_i (L_{-} (y+d_{-}))
%	=-W_{-;\rma}
%	(y)
%\end{align}
%\linfengnote{The above notations for $A_0$ is used in \cite{CWZ} and we may erase them if they're not explicitly used here.}

\begin{lemma}\label{lem:5.2}
%	\footnote{This is similar to \cite[Lemma 5.2]{CLWZ20}.}
	Let $\rho_k=\min\{1,L(1-|y|)\}$.
	For $\nu k^2\les 1$, we have the following bounds
	\begin{align}
		\|W_-\|_{L_y^\infty}\lesssim& \nu^{-1/2}L^{-1}(1+|k(\lambda-\mathcal{U}(-1))|)^{1/2}|A_0(L_-(d_--1))|,
		\label{est:5.2infty}
		\\
		\|W_+\|_{L_y^\infty}\lesssim& \nu^{-1/2}L^{-1}(1+|k(\lambda-\mathcal{U}(1))|)^{1/2}|A_0(-L_+(\overline{d_+}+1))|,
		\label{est:5.3infty}
		\\ 
		\label{est:5.21}
		&\frac{L\|W_-\|_{L^1_y}}{|A_0(L_-(d_--1))|}+\frac{L\|W_+\|_{L^1_y}}{|A_0(-L_+(\overline{d_+}+1)|}\lesssim 1,\\
		&\frac{L^{1/2}\|\rho^{1/2}_kW_-\|_{L^2_y}}{|A_0(L_-(d_--1))|}+\frac{L^{1/2}\|\rho_k^{1/2}W_+\|_{L^2_y}}{|A_0(-L_+(\overline{d_+}+1))|}
		\lesssim 1.
			\label{AiryCan_L2}
	\end{align}
\end{lemma}

\begin{proof}
	The estimates
	\eqref{est:5.2infty} and \eqref{est:5.3infty} are similar to \cite[Lemma~4.3]{CWZ23} and thus we omit the details.
To prove \eqref{est:5.21},
	we use \eqref{EQ133a} and \cite[Lemma~B.3]{CWZ23} to obtain
	\begin{align}
	\begin{split}
		&
		\left|
			\frac{W_{-;\rma}(y)}{A_0(L_-(d_--1))}
			\right|
		=
		\left|
		\frac{A_0' (L_{-} (y+d_{-}))}{A_0 (L_{-} (y+d_{-}))}
		\right|
		\left|
		\frac{A_0 (L_{-} (y+d_{-}))}{A_0(L_-(d_--1))}
	\right|
		\\&\quad
		\lesssim 
		(1+|L_-( d_--1)|)^{1/2}\exp\left\{-c L_{-} (y+1)(1+|L_-(d_--1)|)^{1/2}\right\},
		\label{ptwise_bnd}
	\end{split}
	\end{align}
	which implies 
	\begin{align}
	\begin{split}
			\|W_{-;\rma}
			\|_{L^1_y}
		&\lesssim
		|A_0(L_-(d_--1))| (1+|L_-(d_--1)|)^{1/2}
		\\&\quad
		\times
		\int_{-1}^1\exp\left\{-c L_{-} (y+1)(1+|L_-(d_--1)|)^{1/2}\right\}dy
		\\&
		\lesssim
		\frac{1}{L_{-}(1+|L_-(d_--1)|)^{1/2}}|A_0(L_-(d_--1))| (1+|L_-(d_--1)|)^{1/2}
		\\&
		\lesssim \frac{1}{L}|A_0(L_-(d_--1))|.
	\end{split}
	\end{align}
From \cite[Lemma~4.2]{CWZ23} it follows that
	\begin{align}
	\begin{split}
		\|W_{-;\rme}
		\|_{L^1_y}
		\lesssim \frac{1}{L}|A_0(L_-(d_--1))|
	\end{split}
\end{align}
Hence, we conclude the proof of \eqref{est:5.21} for the term $\Vert W_{-}\Vert_{L^1}$.
The estimate for $\Vert W_{+}\Vert_{L^1_y}$ is similar and thus we omit the details.

Next we prove the weighted $L^2$ estimates in \eqref{AiryCan_L2}.
Using \eqref{ptwise_bnd}, we obtain
\begin{align}
	\begin{split}
	&\frac{L}{|A_0 (L_{-} (d_{-}-1))|}
	\Vert \rho_k^{1/2} W_{-;\rma}
	\Vert_{L^2}
	\\&
	=
	\frac{L}{|A_0 (L_{-} (d_{-}-1))|}
	\left(
	\int_{-1}^{1}
	\rho_k (y)
	|
	A'_0 (
	L_{-} (y+d_{-}) )
	|^2\,dy
	\right)^{1/2}
	\\&
	=
	L
	\left(
	\int_{-1}^{1}
	\rho_k (y)
	\frac{|A_0' (
		L_{-} (y+d_{-}) |^2}{|A_0 (
		L_{-} (y+d_{-}) |^2}
	\frac{|A_0 (
		L_{-} (y+d_{-}) |^2}{|A_0 (L_{-} (d_{-}-1))|^2}
	\,dy
	\right)^{1/2}
	\\&
	\les
	L 
	\lf(
	\int_{-1}^1 \rho_k (y)
	(1+|L_-( d_--1)|)
	\exp\left\{-2c L_{-} (y+1)(1+|L_-(d_--1)|)^{1/2}\right\}
	dy
	\rg)^{1/2}
	\end{split}
\label{EQ132e}
\end{align}
$\bullet$ \textbf{The case $|L_{-} (d_{-}-1)|\leq 1$.} From \eqref{EQ132e} it follows that
\begin{align}
	\begin{split}
	&\frac{L}{|A_0 (L_{-} (d_{-}-1))|}
	\Vert \rho_k^{1/2} W_{-;\rma}
	\Vert_{L^2}
	\\&\quad
	\lesssim
	L
	\left(
	\int_{-1}^{1}
	(1+|L_{-} (y + d_{-})|)
	e^{-C(L_{-} (1+y))}
	\,dy
	\right)^{1/2}
	\lesssim 
	L^{1/2}.
	\end{split}
\label{EQ134a}
\end{align}
$\bullet$
\textbf{The case $|L_{-} (d_{-}-1)|>1$.} 
From \eqref{EQ132e} it follows that
\begin{align}
	&
	\frac{L}{|A_0 (L_{-} (d_{-}-1))|}
	\Vert \rho_k^{1/2} W_{-;\rma}
	\Vert_{L^2}
	\\&\quad
	\lesssim
	L
	\left(
	\int_{-1}^{1}
	\rho_k (y)
	(1+|L_{-} (y + d_{-})|)
	e^{-CL_{-} (1+y)
	|L_{-}(d_{-}-1)|^{1/2}}
	\,dy
	\right)^{1/2}
	\\&\quad
	\lesssim L
	\left(
	\int_{-1}^{1}
	(1+ L_{-} (1+y) + L(y+1) |L_{-} (d_{-}-1)|)
	e^{-C L_{-} (1+y) |L_{-} (d_{-}-1)|^{1/2}}
	\,dy
	\right)^{1/2}
	\\&\quad
	\lesssim
	L^{1/2}.
	\label{EQ100}
\end{align}
Combining \eqref{EQ134a}
and \eqref{EQ100}, we conclude the estimate for $\Vert \rho_k^{1/2} W_{-;\rma} \Vert_{L^2}$ in \eqref{AiryCan_L2}.
For the reminder term $\Vert \rho_k^{1/2} W_{-;\rme}\Vert_{L^2}$,
we have
\begin{align}
	\begin{split}
		&
		\|\rho_k^{1/2}W_{-;\rme}\|_{L^2}
	\lesssim
	\|W_{-;\rme}\|_{L^2}
	\lesssim \|W_{-;\rme}\|_{L^1}^{1/2}\|W_{-;\rme}\|^{1/2}_{L^\infty}.
%	\\&\quad
%	\lesssim
%	\left(
%	L_{-}^{-7/4} |k|^{-1/4}
%	|A_0(L_-(d_--1))|
%	\right)^{1/2}
%	|A_0(L_-(d_--1))|^{1/2}
%	\\&\quad
%	\lesssim L^{-7/8}|k|^{-1/8}|A_0(L_-(d_--1))|
%	\lesssim
%	L^{-1/2}
%	|A_0(L_-(d_--1))|.
%	\label{EQ100a}
	\end{split}
\end{align}
Using \cite[Lemma~4.2]{CWZ23}, we obtain
\begin{align}
	\begin{split}
		&
		\|\rho_k^{1/2}W_{-;\rme}\|_{L^2}
%		\lesssim
%		\|W_{-;\rme}\|_{L^2}
%		\lesssim \|W_{-;\rme}\|_{L^1}^{1/2}\|W_{-;\rme}\|^{1/2}_{L^\infty}
%		\\&\quad
		\lesssim
		\left(
		L_{-}^{-7/4} |k|^{-1/4}
		|A_0(L_-(d_--1))|
		\right)^{1/2}
		|A_0(L_-(d_--1))|^{1/2}
		\\&\quad
		\lesssim L^{-7/8}|k|^{-1/8}|A_0(L_-(d_--1))|
		\lesssim
		L^{-1/2}
		|A_0(L_-(d_--1))|.
		\label{EQ100a}
	\end{split}
\end{align}
Similar arguments prove the $\Vert \rho_k^{1/2} W_{+;\rma}\Vert_{L^2_y}$ part in \eqref{AiryCan_L2} and thus
we conclude the proof of \eqref{AiryCan_L2}.
\end{proof}

We recall the estimates for the coefficients $A_{\pm\pm}$ and $A_{\pm\mp}$ from \cite{CWZ23}. 
\begin{lemma}\label{lem:4.4CWZ23}
	\cite[Lemma~4.4]{CWZ23}
	There exists a constant $K>0$ such that if $\min \{L_{+}, L{-}\} \geq K$, then
	\begin{align}
		&		|A_{--}|\geq \frac{1}{CL_{-}} |A_0 (L_{-} {(d_{-}-1)})|
		,\quad
		|A_{++}|\geq \frac{1}{CL_{+}} |A_0 (-L_{+} (\overline{d_{+}}+1))|;\\
		&	\left|
		\frac{A_{-+}}{A_{--}}
		\right|	
		\leq
		\frac{\sqrt{2}}{2}
		,\quad
		\left|
		\frac{A_{+-}}{A_{++}}
		\right|
		\leq \frac{\sqrt{2}}{2}.
	\end{align}
	Moreover, the following two estimates hold
	\begin{align}
		&|A_{-+}|\lesssim (L_-^{-7/4}|k|^{-1/4}+L_-^{-2})|A_0(L_{-} (d_{-}-1))|,\\
		&|A_{+-}|\lesssim (L_+^{-7/4}|k|^{-1/4}+L_+^{-2})|A_0(-L_+(\overline{d_+}+1))|.
	\end{align}
\end{lemma}
A direct consequence of the lemma is that the following lower bound for the Evan's function $\mathfrak{D}$
\begin{align}
\begin{split}
		|\mathfrak{D}|
		&=
		|A_{--}A_{++}-A_{-+}A_{+-}|\geq \frac{|A_{--}A_{++}|}{2}
		\\&
	\geq
	\frac{|A_0(L_-(d_--1))| |A_0(-L_+(\overline{d_+}+1))|}{2CL^{2}}.
\end{split}
	\label{Evan's_func_bnd}
\end{align}
This lower bound will turn out to be crucial to estimate the coefficients $\mathfrak C^\pm_\pm,$ and $\mathfrak C^{\pm}_\mp$ introduced in  \eqref{EQ101a}. 

%With these information, we can estimate the unit boundary corrector $w_\pm$ in \eqref{EQ101a}.
\begin{lemma}
	\label{lem:uni_bcrrctr} 
%\footnote{This is similar to \cite[Proposition 4.3, 4.4]{CLWZ20}.}
There exists a constant $K>0$ such that if $\min \{L_{+}, L{-}\} \geq K$ and $\nu k^2 \les 1$, then
\begin{align}
	\label{est:P4.3ARMA1}
	&\frac{\nu^{1/2}\|w_-\|_{L^\infty}}{(1+|k||\lambda-\mathcal{U}(-1)|)^{1/2}}+\|w_-\|_{L^1}\lesssim 1,
	\\
	\label{est:P4.3ARMA2}
	&
	\frac{\nu^{1/2}\|w_+\|_{L^\infty}}{(1+|k||\lambda-\mathcal{U}(1)|)^{1/2}}+\|w_+\|_{L^1}
	\lesssim 
	1,
	\\
	\label{est:P4.4AR}%54
	& \|\rho_k^{-1/4}w_+\|_{L^2}\lesssim
	\nu^{-7/24}
	k^{-1/12}
	(1+|k||\lambda-\mathcal{U}(1)|)^{3/8},
	\\&
	\label{est:P4.4AR3}
	\|\rho_k^{-1/4}w_-\|_{L^2}\lesssim
	\nu^{-7/24}
	k^{-1/12}
	(1+|k||\lambda-\mathcal{U}(-1)|)^{3/8},
	\\ 
	\label{est:P4.4AR2}
	&\|\rho_k^{1/2}w_-\|_{L^2}+\|\rho_k^{1/2}w_+\|_{L^2}\lesssim L^{1/2} .
\end{align}
\end{lemma}
\begin{proof}
%	Without loss of generality, we focus on the $w_-$. 
	First we prove the estimate \eqref{est:P4.3ARMA1}. Thanks to Lemma \ref{lem:4.4CWZ23}, \eqref{C_coeff_exp}, and the lower bound on the Evan's function \eqref{Evan's_func_bnd}, we obtain
	\begin{align}\label{CMP2411}
		|\mathfrak C_{-}^-|= \left|\frac{A_{++}}{\mathfrak D}\right|\lesssim \left|\frac{A_{++}}{A_{++}A_{--}}\right|\lesssim \frac{L}{|A_0(L_-(d_--1))|}
	\end{align}
%and
%\begin{align}
%	|\mathfrak{C}^{+}_{-}|
%	\lesssim
%	\left|
%	\frac{A_{-+}}{A_{++} A_{--}}
%	\right|
%	\lesssim
%	\frac{L}{A_0 (-L_{+} \overline{d_{+}} +1)}.
%\end{align}
%
%\coly
	and
	\begin{align}\label{CMP2412}
		|\mathfrak C_-^+|\lesssim \left|\frac{A_{-+}}{A_{++}A_{--}}\right|\lesssim \frac{L^2(L^{-7/4}|k|^{-1/4}+L^{-2})}{|A_0(-L_+( \overline{d_+}+1))|}\lesssim \frac{\nu^{-1/12}|k|^{-1/6}}{|A_0(-L_+( \overline{d_+}+1))|},
	\end{align}	
where we used $\nu k^2 \les 1$ in the last step.
	Using Lemma~\ref{lem:5.2}, \eqref{EQ101a}, \eqref{est:5.2infty}, \eqref{est:5.3infty}, \eqref{CMP2411}, and \eqref{CMP2412}, we arrive at
	\begin{align}
		&
		\|w_-\|_{L^\infty}
		\leq
		|\mathfrak C_-^-|\|W_-\|_{L^\infty}
		+
		|\mathfrak C_-^+|\|W_+\|_{L^\infty}
		\\&\quad
%		\lesssim  \frac{L}{|A_0(L_-(d_--1))|}\|W_-\|_{L^\infty} 
%		+
%		\frac{\nu^{-1/12} |k|^{-1/6}}{|A_0(-L_+( \overline{d_+}+1))|}\|W_+\|_{L^\infty}
%		\\&\quad
		\lesssim
		\nu^{-1/2}(1+|k||\lambda-\mathcal{U}(-1)|)^{1/2}
		+
		\nu^{-1/4} k^{-1/2}(1+|k||\lambda-\mathcal{U}(1)|)^{1/2}
		\\&\quad
		\lesssim
		\nu^{-1/2}(1+|k||\lambda-\mathcal{U}(-1)|)^{1/2}.\label{est:quote}
	\end{align}
	For the $L^1$ norm, we use \eqref{est:5.21}, \eqref{CMP2411}, and \eqref{CMP2412} to get
	\begin{align}
		\|w_-\|_{L^1}
		&\leq
		|\mathfrak C_-^-|\|W_-\|_{L^1}+|\mathfrak C_-^
		+
		|\|W_+\|_{L^1}
%		\lesssim \frac{L}{|A_0(L_-(d_--1))|}\|W_-\|_{L^1} 
%		+
%		\frac{\nu^{-1/12} |k|^{-1/6}}{|A_0(-L_+(\overline{d_+}+1))|}\|W_+\|_{L^1}
%		\\&
		\lesssim 1.
	\end{align} 
	This completes the proof of  \eqref{est:P4.3ARMA1}. The proof of \eqref{est:P4.3ARMA2} is similar and thus we omit the details.

	The proofs of the estimates \eqref{est:P4.4AR}, \eqref{est:P4.4AR3}, and \eqref{est:P4.4AR2}
	are similar to 
	\cite[Proposition~4.4]{CLWZ20}  and thus we omit the details.
\end{proof}

\subsection{Space-time Estimates}\label{sec:space-time}
In this section, we prove the space-time estimates for the
Navier-Stokes equations linearized around a time-independent shear flow.
%First, we recall that the inhomogeneous solution $w_I$ satisfies 
%\begin{align}
%	\label{eq:w_I2}
%	\left\{
%	\begin{aligned}
%		&-\nu(\pa_y^2-k^2)w_{I;k}
%		+
%		ik\lf(\mathcal{U}-\lambda\rg)w_{I;k}-ik\;\mathcal{U}''\psi_{I;k}	
%		=	F_k(\lambda,y),
%		\\&
%		\psi_{I;k}
%		=(-\Delta_{k})^{-1} w_{I;k},\quad
%		\int_{-1} ^1 e^{\pm ky}w_{I;k} dy
%		=0. 
%	\end{aligned}
%	\right.
%\end{align}
%\begin{proposition}
%	\label{pro:w_I}
%%	\footnote{This is similar to  \cite[Proposition 6.5]{CLWZ20}}
%Let $w^I$ be the solution of \eqref{om_I} 
%with the forcing 
%\begin{align}
%	\mathbb{F}_k
%	=
%	ik\mathbbm{f}_1+\pa_y\mathbbm{f}_2+\mathbbm{f}_3,\quad \mathbbm{f}_1,\mathbbm{f}_2\in L^2,\,  \mathbbm{f}_3\in H_0^1.
%\end{align} 
%Assume that $\|\U-y\|_{H^4}\ll 1$ and $\nu|k|^2\leq \|\U'\|_{L^\infty}$.
%Then we have
%\begin{align}
%\begin{split}
%	&\nu^{1/3}|k|^{2/3}\|\rho_k^{1/2}w^I\|_{L^2L^2}^2+
%	|k|^2\|u^I\|_{L^2L^2}^2
%	+
%	\nu^{1/2}|k|\|w^I\|_{L^2L^2}^2
%	+
%	\nu^{1/2}|k|\|w^I\|_{L^\infty L^2}^2
%	\\&
%	\lesssim
%	\log^2 L
%	\left( \nu^{-1/3}|k|^{4/3}\|\mathbbm{f}_1\|_{L^2L^2}^2+\nu^{-1}\|\mathbbm{f}_2\|_{L^2L^2}^2
%	+
%	\|\na_k\mathbbm{f}_3\|_{L^2 L^2}^2
%	\right).
%\end{split}
%\end{align}
%\end{proposition}

\begin{proof}[Proof of Proposition~\ref{pro:w_I}]
	Without loss of generality, we assume that $k>0$.
	Let $\lambda= \lambda_r + i \lambda_i$ where $\lambda_i=  -\delta_0 \nu^{1/3} k^{-1/3}$.
We take the Fourier transform in $t$ to get
\begin{align}
	\label{Fourier}	 
	\begin{split}
		&w_{k} (\lambda,y) :=\int_0^\infty \omega_{k}^I(t,y)e^{it  k \lambda}dt, 
		\quad
		\phi_{k}(\lambda,y):=\int_0^\infty \psi_{k}^I(t,y)e^{it  k\lambda}dt, 
		\\&
		F_j (\lambda,y):=\int_0^\infty \mathbbm{f}_j (t,y) e^{it  k\lambda}dt,\quad
		j=1,2,3.
	\end{split}
\end{align}
Thus, from \eqref{om_I} it follows that
\begin{align}
	\begin{split}
	&	-\nu \de_k
	w_k
	+ik (\U -\lambda)w_k
	-ik \U'' \phi_k
	=
	ik F_1 + \pa_y F_2 + F_3,
	\\&
	\phi_k = \de_k^{-1} w_k,\quad
	\int_{-1}^{1} 
	w_k e^{\pm ky}dy=0.
	\end{split}
\end{align}
We further decompose $w_k$ as
\begin{align}
\begin{split}
	w_k
	=
	w_{\rm Na}^{(1)}
	+
	w_{\rm Na}^{(2)}
	+
	w_{\rm Na}^{(3)}
	+
	(c^{(1)}_{+}+c^{(2)}_{+}+c^{(3)}_{+}) w_{+;k}
		+
	(c^{(1)}_{-}+c^{(2)}_{-}+c^{(3)}_{-}) w_{-;k},
\end{split}
\end{align}
where $w_{\rm Na}^{(1)}$, $w_{\rm Na}^{(2)}$, and $w_{\rm Na}^{(3)}$ solves \eqref{Nav_eqn} with the right-hand side forcing $ik F_1$, $\pa_y F_2$, and $F_3$, respectively, and
\begin{align}
	c_{+}^{(j)}=
	-\pa_y \phi_{\rm Na}^{(j)}(y=1),
	\quad
	c_{-}^{(j)} = -\pa_y \phi_{{\rm Na}}^{(j)} (y=-1),
	\quad j=1,2,3.
\end{align}
Using the Plancherel's theorem, we obtain
\begin{align}
\begin{split}
	&\int_{0}^{\infty}
	\Vert
	e^{-tk \lambda_i}
	 \omega_k^I (t)\Vert_{L^2}^2
	dt \sim
	\int_{\rr} \Vert w_k 
	(\lambda_r)\Vert_{L^2}^2
	d \lambda_r,
	\\&
	\int_{0}^{\infty}
	\Vert
	e^{-tk \lambda_i}
	\psi_k^I (t)\Vert_{L^2}^2 dt
	\sim
	\int_{\rr} \Vert \phi_k (\lambda_r)\Vert_{L^2}^2
	d \lambda_r,
	\\&
	\int_{0}^{\infty}
	\Vert
	e^{-tk \lambda_i}
	u_k^I (t)\Vert_{L^2}^2 dt
	\sim
	\int_{\rr} \Vert 
	(ik, \pa_y)
	\de_{k}^{-1} w_k
	\Vert_{L^2}^2
	d \lambda_r,
	\\&
	\int_{0}^{\infty}
	\Vert
	e^{-tk \lambda_i}
	\mathbbm{f}_i (t)\Vert_{L^2}^2 dt
	\sim
	\int_{\rr} \Vert F_j\Vert_{L^2}^2
	d \lambda_r,
	\quad j=1,2,3.
\end{split}
\end{align}
The proofs of the components $(w_{\rm Na}^{(1)}, c_{+}^{(1)} w_{+;k}, c_{-}^{(1)} w_{-;k})$ and $(w_{\rm Na}^{(2)}, c_{+}^{(2)} w_{+;k}, c_{-}^{(2)} w_{-;k})$ are analogous to \cite[Proposition~6.5]{CLWZ20} which rely on Proposition~\ref{pro_Nav}, Lemmas~\ref{lem:cpm_est}, \ref{lem:cpm_est_2} and~\ref{lem:uni_bcrrctr}, thus we omit the details.

Next we estimate the contribution of the component $(w_{\rm Na}^{(3)}, c_{+}^{(3)} w_{+;k}, c_{-}^{(3)} w_{-;k})$.
From \eqref{cest_FH1} and Lemma~\ref{lem:uni_bcrrctr}, we obtain
\begin{align}
	\begin{split}
		&
	|c^{(3)}_{+}| 
	\Vert
	w_{+}	\Vert_{L^2}
	+	|c^{(3)}_{+}| 
	\Vert
	w_{+}	\Vert_{L^2}
	\les
	|k|^{-1/2} \nu^{-1/4}
	\log L 
	\Vert \nabla_k F_3\Vert_{L^2},
	\\&
		|c^{(3)}_{+}| 
	\Vert
	\rho_k^{1/2}
	w_{+}	\Vert_{L^2}
	+	|c^{(3)}_{+}| 
	\Vert
	w_{+}	\Vert_{L^2}
	\les
	|k|^{-1/3} \nu^{-1/6}
	\log L 
	\Vert \nabla_k F_3\Vert_{L^2},
	\\&
		|c^{(3)}_{+}| 
	\Vert
	w_{+}	\Vert_{L^1}
	+	|c^{(3)}_{+}| 
	\Vert
	w_{+}	\Vert_{L^1}
	\les
	|k|^{-1/2} 
	\log L 
	\Vert \nabla_k F_3\Vert_{L^2}.
	\end{split}
\end{align}
Let $\omega_{\rm Na}^3$ be the solution of
\begin{align}
	&
	\partial_t \omega^{(3)}_{\rm Na}
	-\nu \de_k  \omega^{(3)}_{\rm Na}
	+
	ik \U  \omega^{(3)}_{\rm Na}
	-ik \U'' \psi^{(3)}_{\rm Na}
	= \mathbbm{f}_3,
	\\&
	\omega^{(3)}_{\rm Na}|_{t=0}=0,\quad
	\omega^{(3)}_{\rm Na}|_{y=\pm 1}=0.
\end{align}
By Plancherel's theorem and \eqref{F_in_H1_est}, we have
\begin{align}
	\begin{split}
		&
		\nu^{1/3}|k|^{2/3}
		\Vert e^{-tk\lambda_i} \omega_{\rm Na}^{(3)} (t)\Vert_{L^2 L^2}^2
		\les
		\nu^{1/3}|k|^{2/3}
		\Big\Vert
		\Vert w^{(3)}_{\rm Na} (\lambda_r)\Vert_{ L^2_y}
		\Big\Vert^2_{\rr}
		\\&
		\les
			|k|^{-2}
		\Big\Vert
		\Vert \nabla_k F_3\Vert_{L^2}
		\Big\Vert^2_{\rr}
%		\\&
		\les
		\Vert e^{-tk\lambda_i}
		\nabla_k \mathbbm{f}_3\Vert_{L^2 L^2}^2.
	\end{split}
\end{align}
Hence, the bounds on $\nu^{1/3} |k|^{2/3} \Vert e^{-tk \lambda_i} \rho_k^{1/2} \omega^I\Vert_{L^2 L^2}^2$ and $\nu^{1/2} |k| \Vert e^{-tk \lambda_i} \omega^I\Vert_{L^2 L^2}^2$ are consistent with the right-hand side of \eqref{w_ I_est}.
By \eqref{F_in_H1_est},  \eqref{cest_FH1}, \cite[Lemma~9.3]{CLWZ20}, and Lemma~\ref{lem:uni_bcrrctr},
we have
\begin{align}
	\begin{split}
		&
	\Vert (ik, \pa_y) \de_k^{-1} (w^{(3)}_{\rm Na}+c^{(3)}_+ w_{+}+ c^{(3)}_{-} w_{-})\Vert_{L^2}\\&
	\les
	|k|^{-2} \Vert \nabla_k F_3\Vert_{L^2}
	+
	|k|^{-1/2}
	|c^{(3)}_{+}|
	\Vert w_{+}\Vert_{L^1}+
	|k|^{-1/2}
	|c^{(3)}_{-}|
	\Vert w_{-}\Vert_{L^1}
	\\&
	\les
	|k|^{-1}
	\Vert \nabla_k F_3\Vert_{L^2}.
	\end{split}
\end{align}
This completes the proof of the estimate of $|k|^2 \Vert e^{-tk \lambda_i} u_k^I\Vert_{L^2 L^2}^2$.

The proof of the term $\nu^{1/2} |k|\Vert e^{-tk \lambda_i} \omega_k^I\Vert_{L^\infty L^2}^2$ is similar to \cite[Proposition~6.5]{CLWZ20} and thus we omit the details.
\end{proof}

Now we are ready to prove Proposition~\ref{pro:w_H}.
\begin{proof}
	[Proof of Proposition~\ref{pro:w_H}]
We recall the decomposition of the homogeneous component of the solution 
$\omega^H$ in Section \ref{sec:tmnd_H}:
\begin{align}
	\left\{
	\begin{aligned}
		&\omega^{H}_{k} (t,y)
	=\omega^{(1)}_{k} (t,y)
	+\omega^{(2)}_{k} (t,y)
	+\omega^{(3)}_{k} (t,y),
	\\
	&(\partial_{y}^2 - k^2)
	\psi_{k}^{(i)}
	= \omega_k^{(i)},\quad 
	\psi^{(i)}\big|_{y=\pm 1}=0,
	\quad i=1,2,3.
	\end{aligned}
\right.
\end{align}

First we estimate the term $\omega_k^{(1)}(t,y)$. 
Since $\nu k^2 \les 1$, we have
\begin{align}
\begin{split}
	&
		\nu^{1/3}|k|^{2/3}
		\|
		e^{\delta(\nu k^2)^{1/3}t}
		\rho_k^{1/2}\omega_k^{(1)}\|_{L^2L^2}^2
		+
		\nu^{1/2}|k|
		\|
		e^{\delta(\nu k^2)^{1/3}t}
		\omega_k^{(1)}\|_{L^2L^2}^2
		\\&\qquad
		\les
		\nu^{1/3} |k|^{2/3} 
		\Vert
		e^{\delta(\nu k^2)^{1/3}t} \omega^{(1)}_k\Vert_{L^2 L^2}^2.
\end{split}
\end{align}
Direct computation of the explicit expression \eqref{eq:om_1} gives
\begin{align}
		\label{Est_omg_H1}
	\begin{split}
		&
	\nu^{1/3} |k|^{2/3} 
	\Vert
	e^{\delta(\nu k^2)^{1/3}t} \omega^{(1)}_k\Vert_{L^2 L^2}^2
	\les
	\nu^{1/3} |k|^{2/3}
		\Vert \omega_{{\rm{in}};k}\Vert_{L^2}^2
	\int_0^\infty
	e^{2\delta (\nu k^2)^{1/3} t-C\nu k^2 t^3
	-2\nu k^2 t}
	dt
	\\&\quad
	\les
		\nu^{1/3} |k|^{2/3}
	\Vert \omega_{{\rm{in}};k}\Vert_{L^2}^2
	\lf(	\int_0^{(\nu k^2)^{-1/3}}
	e^{2\delta (\nu k^2)^{1/3} t}
	dt
	+
	\int_{(\nu k^2)^{-1/3}}^\infty
	e^{-C\nu k^2 t^3} dt
	\rg)
	\\&\quad
	\les
	\Vert \omega_{{\rm in};k}\Vert_{L^2}^2
	\end{split}
\end{align}
and
\begin{align}
	\label{Est_omg_H5}
	\nu^{1/2}|k|
	\|e^{\delta(\nu k^2)^{1/3}t}
	\omega_k^{(1)}\|_{L^\infty L^2}^2
	\lesssim
	\nu^{1/2} |k|
	\Vert \omega_{{\rm in};k}\Vert_{L^2}^2
	\les \|\omega_{{\rm{in}};k}\|_{L^2}^2.
\end{align}
Next we estimate the term $|k|^2\|e^{\delta(\nu k^2)^{1/3}t} u^{(1)}\|_{L^2L^2}^2$.
Recall that $G_k$ is the Green's function for $\de_k=-|k|^2+\partial_{yy}$ with homogeneous Dirichlet boundary conditions, which has the explicit form \eqref{green}.
%\linfeng{need to revise? we want to add $e^{\delta \nu^{1/3} t}$}
Using integration by parts, we obtain
\begin{align}
	\begin{split}
	u^{(1)}_k 
	&=
	(-\partial_{y}, ik) \Delta_{k}^{-1} \omega^{(1)}
	\\&
	=\int_{-1}^1 (-\pa_y,ik)G_k(y,\wt y)\omega_{{\rm in};k}(\wt y)
	\exb{-\frac{1}{3}(\mc U')^2\nu k^2 t^3-\nu k^2t}
	\\&\quad\quad\times
	\frac{\pa_{\wt y}\exb{- it\mathcal{U}k}}{-it \mathcal{U}'k}d\wt y
	\\&
	=
	\frac{1}{ikt }\int_{-1}^1 \pa_{\wt y}\lf(
	\frac{ \omega_{{\rm in};k}(\wt y)}{  \mathcal{U}' }(-\pa_y,ik)G_k(y,\wt y)
	\exb{-\frac{1}{3}(\mc U')^2\nu k^2 t^3-\nu k^2t}\rg)
		\\&\quad\quad
		\times
	\exb{- it\mathcal{U}k}d\wt y,
	\end{split}
	\label{Est_omg_H2}
\end{align}
where we used $G_k(y,\wt y=\pm 1)=\pa_yG_k(y,\wt y=\pm 1)=0$ in the last step. 
%Note that the expression \eqref{Est_omg_H2} is understood in the sense of distribution since there is a dirac singularity in the expression. 
Direct calculation gives
\begin{align}
	\|e^{\delta(\nu k^2)^{1/3} t}
	u^{(1)}_k(t)\|_{L^2_y}
	\lesssim
	\frac{1}{k \langle t \rangle }\|\omega_{\rmi;k}\|_{H^1_k},
	\label{Est_omg_H3}
\end{align}
which leads to
\begin{align}
	\begin{split}
	|k|^2
	\|e^{\delta (\nu k^2)^{1/3}t}
	u^{(1)}_k\|_{L^2_tL^2_y}^2
	\lesssim \|\omega_{{\rm{in}};k}\|_{H_k^1}^2.
	\label{Est_omg_H4}
	\end{split}
\end{align}

Next we estimate the term $\omega_k^{(2)}(t,y)$. 
Note that $\omega_k^{(2)}(t,y)$ is the solution of \eqref{eq:om_2}, which can be considered as the solution of \eqref{om_I} with the right-hand side forcing 
\begin{align}
	\mathbb{F}_k
	=
	ik \mathbbm{f}_1
	+\pa_y \mathbbm{f}_2+\mathbbm{f}_3,
\end{align}
where 
\begin{align}
	\mathbbm{f}_1= 
	\frac{1}{ik} (\nu\pa_y^2\omega_k^{(1)}+(\mc U')^2 \nu k^2t^2\omega_k^{(1)}),
	\quad
	\mathbbm{f}_2=0,\quad
	\mathbbm{f}_3= ik \U'' \psi^{(1)}_k.
\end{align}
Using the explicit expression of  $\omega_k^{(1)}$ in \eqref{eq:om_1}, we write
\begin{align}
	&\mathbbm{f}_1=
	\frac{1}{ik}
	\exp\lf\{- it\mathcal{U} k-\frac{1}{3}(\mc U')^2\nu k^2 t^3-\nu k^2t\rg\}
	\\&\times
	\Bigg[\nu \lf(\pa_y^2\omega_{{\rm in};k}-\frac{4}{3}\pa_y\omega_{{\rm in};k}\mc U'\mc U'' \nu k^2t^3-\frac{2}{3}(\mc U'\mc U'')'\nu |k|^2t^3\omega_{{\rm in};k}+\frac{4}{9}(\mc U'\mc U'')^2\nu^2|k|^4t^6\omega_{{\rm in};k}\rg)\\
	&\qquad-2\nu i\mc U' kt\lf(\pa_y \omega_{{\rm in};k}-\frac{2}{3}\omega_{{\rm in};k}\mc U'\mc U'' \nu k^2t^3\rg) -\nu i\mc U'' kt\omega_{{\rm in};k}\Bigg].
\end{align} 
Direct computation gives
%\begin{align}
%	\Vert \mathbbm{f}_1 \Vert_{L^2 L^2}^2
%	\les
%	\nu^2 k^2
%\end{align}
%and
%\begin{align}
%	\Vert \mathbbm{f}_2
%	\Vert_{L^2 L^2}^2
%	\les
%	\nu^2 \Vert \pa_y w^{(1)}_k \Vert_{L^2 L^2}^2
%	\les
%	\nu \Vert w_{{\rm in}, k}\Vert_{H^1_k}^2
%\end{align}
\begin{align}
	\| e^{\delta(\nu k^2)^{1/3}t}
	\mathbbm{f}_1 \|_{L_t^2 L_y^2}^2
	\lesssim
 	\nu^{1/3} |k|^{-10/3}
 	\|\omega_{{\rm in};k}\|_{H_k^2}^2.
\end{align}
Using \eqref{Est_omg_H4}, we obtain
\begin{align}
	\Vert e^{\delta(\nu k^2)^{1/3}t}
	\nabla_k \mathbbm{f}_3 \Vert_{L^2 L^2}^2
	=
	\| e^{\delta(\nu k^2)^{1/3}t}
	\nabla_k (i\U'' k\psi_k^{(1)})\|_{L^2 L^2}^2
	\lesssim
	\|\omega_{{\rm in};k}\|_{H_k^1}^2
	\les
	k^{-2} 	\|\omega_{{\rm in};k}\|_{H_k^2}^2.
\end{align}
An application of Proposition~\ref{pro:w_I} yields
\begin{align}
	\begin{split}
		&\nu^{1/3}|k|^{2/3}
		\|e^{\delta(\nu k^2)^{1/3}t} \rho_k^{1/2}
		\omega_k^{(2)}
		\|_{L^2L^2}^2
		+
		|k|^2\|e^{\delta(\nu k^2)^{1/3}t}
		u_k^{(2)}\|_{L^2L^2}^2
		+
		\nu^{1/2}|k|\|e^{\delta(\nu k^2)^{1/3}t}\omega_k^{(2)}\|_{L^2L^2}^2
		\\&\quad
		+
		\nu^{1/2}|k|\|e^{\delta(\nu k^2)^{1/3}t}\omega_k^{(2)}\|_{L^\infty L^2}^2
		\lesssim
		(\log^2 L)
		k^{-2}
		\|
		\omega_{{\rm{in}};k}
		\|_{H^2_k}^2.
	\end{split}
\end{align}

Finally, we estimate the term $\omega_k^{(3)}(t,y)$. We use similar arguments in \cite[Proposition 6.6]{CLWZ20} and replace the estimates  \cite[(6.12) (6.13)]{CLWZ20} by  \cite[(7.21)]{CWZ23}. We omit further details for the sake of brevity.
Therefore, we conclude the proof of the Proposition.
%\myr{HS: Check! It seems that the argument there is independent of Proposition \ref{pro:w_I}}. 
\end{proof}

Next, we derive the space-time estimates of the Dirichlet problem
\begin{align}
	\left\{
	\begin{aligned}
		&\pa_t \omega_{k} -\nu\de_k\omega_k+\mathcal{U}ik \omega_{k}-\mathcal{U}''ik \psi_{k}=	ik \mathbbm{f}_1+\pa_y \mathbbm{f}_2,
		\\&
		\de_k\psi_{k}=\omega_{k},\quad\psi_{k}\big|_{y=\pm1}=0,\quad \int_{-1}^1 e^{\pm k y}\omega_{k} dy=0,\\
		&\omega_{k}(t=0,y)=\omega_{{\rm in};k}(y).
	\end{aligned}
\right.
\label{lnrzd_NS}
\end{align} 
%
%
%\begin{align}
%	\label{lnrzd_NS}
%\begin{split}
%	&\pa_t w_{k}
%	+
%	Uikw_{k}-
%	U''ik\phi_k
%	-
%	\nu\de_kw_{k}
%	=
%	-ik \mathbbm{f}_1+\pa_y \mathbbm{f}_2,\\
%	&\int_{-1}^1 e^{\pm k y}w_k(t,y)dy\equiv 0,\quad
%	w_k(t=0,y)=\omega_{{\rm in};k}(y).
%\end{split}
%\end{align}

\begin{proposition}
\label{pro_lin}
Let $\omega_k$ be the solution of \eqref{lnrzd_NS}.
For $\nu|k|^2\leq \|\U'\|_{L^\infty}$, we have
\begin{align}
\begin{split}
	&\|
	e^{\delta(\nu k^2)^{1/3} t}
	\rho_k^{3/2}\omega_k\|_{L^\infty L^2}^2
	+
	\nu^{1/3}|k|^{2/3}
	\|
	e^{\delta(\nu k^2)^{1/3} t}
	\rho_k^{1/2}\omega_k\|_{L^2L^2}^2
	+
	k^2
	\|
		e^{\delta(\nu k^2)^{1/3} t}
	u_k\|_{L^2L^2}^2
		\\&\quad\quad
	+
	\nu^{1/2}
	|k|
	\|
		e^{\delta(\nu k^2)^{1/3} t}
	\omega_k\|_{L^2L^2}^2
	+
	\nu^{1/2}
	|k|
	\|
		e^{\delta(\nu k^2)^{1/3} t}
	\omega_k\|_{L^\infty L^2}^2
	\\&\quad
	\lesssim 
	|\log \nu|^2
	\big(
	|k|^{-2}
	\|\omega_{{\rm in};k}\|_{H_k^{2}}^2
	+
	\nu^{-1/3}|k|^{4/3}
	\|
		e^{\delta(\nu k^2)^{1/3} t}
	\mathbbm{f}_1\|_{L^2L^2}^2
%		\\&\qquad
	+
	\nu^{-1}
	\| 	e^{\delta(\nu k^2)^{1/3} t}
	\mathbbm{f}_2\|_{L^2L^2}^2
	\big).
	\label{EQ102c}
\end{split}
\end{align}
%\begin{align}
%\begin{split}
%		&\|e^{\delta \nu^{1/3}|k|^{2/3}t}\na_k\phi \|_{L^{\infty}L^2}+\|e^{\delta \nu^{1/3}|k|^{2/3}t}\na_k\phi \|_{L^2L^2}\\
%	&\lesssim |k|^{-1}(\nu^{1/3}|k|^{-1/3}\|\pa_y w_{in}\|_{L^2}+\|w_{in}\|_{L^2})+\nu^{-1/2}|k|^{-1}\|e^{\delta \nu^{1/3}|k|^{2/3}t}(f_1,f_2)\|_{L^2L^2}\\
%	&\quad+\nu^{-1/6}|k|^{-4/3}\|e^{\delta \nu^{1/3}|k|^{2/3}t} f_3\|_{L^2L^2}+|k|^{-2}\|e^{\delta \nu^{1/3}|k|^{2/3}t}\|\na_k f_4\|_{L^2L^2}; \\
%	&f_3\equiv f_4\equiv 0: \quad \|(1-y^2)^{1/2}w_k\|_{L^\infty L^2}+\nu^{1/4}|k|^{1/2}\|w\|_{L^2L^2}+\nu^{1/4}|k|^{1/2}\|w\|_{L^\infty L^2}\\
%	&\hspace{5cm}\lesssim \|w_{in}\|_{L^2}+\nu^{1/6}|k|^{-2/3}\|\pa_yw_{in}\|_{L^2} +\nu^{-1/2}\|(f_1,f_2)\|_{L^2L^2}\label{wght_vor}%:P
%\end{split}
%\end{align}
%\myr{Some estimates come from ChenLiWeiZhang Proposition 6.7. Check page 43 and (6.19) of  ChenLiWeiZhang. There, a different $\rho_k$ is used.}
\end{proposition}
\simhe{\begin{rmk}In fact, it seems that starting from here, we can derive the $\|u_k\|_{L_t^\infty L_y^\infty}$ estimate following the argument in Section 6.4 (equation (6.20)). 
\end{rmk}}

\begin{proof}
We decompose the solution of \eqref{lnrzd_NS} as $\omega_k=\omega_k^H+\omega^I_k$, where $\omega^I_k$ solves \eqref{om_I} with $\mathbb{F}_k=ik\mathbbm{f}_1+\partial_y \mathbbm{f}_2$ and $\omega^H_k$ solves \eqref{om_H}.
Using Propositions~\ref{pro:w_I} and~\ref{pro:w_H}, we obtain
\begin{align}
	\begin{split}
		&\nu^{1/3}|k|^{2/3}
		\|		e^{\delta(\nu k^2)^{1/3} t}
		\rho_k^{1/2}\omega_k\|_{L^2L^2}^2
		+
		|k|^2\|		e^{\delta(\nu k^2)^{1/3} t}
		u_k\|_{L^2L^2}^2
				\\&\quad
				+
		\nu^{1/2}|k|
		\|		e^{\delta(\nu k^2)^{1/3} t} \omega_k\|_{L^2L^2}^2
		+\nu^{1/2}|k|
		\|		e^{\delta(\nu k^2)^{1/3} t} \omega_k\|_{L^\infty L^2}^2
	\\&\quad
	\lesssim 
	|\log\nu|^2
	\big(
	\nu^{1/3}|k|^{-4/3}\|\pa_y\omega_{{\rm in};k}\|_{L^2}^2
	+
	|k|^{-2}
	\|\omega_{{\rm in};k}\|_{H_k^{2}}^2
	\\&\quad\quad
	+
	\nu^{-1/3}|k|^{4/3}
	\|	e^{\delta(\nu k^2)^{1/3} t}
	\mathbbm{f}_1\|_{L^2L^2}^2
	+
	\nu^{-1}
	\|	e^{\delta(\nu k^2)^{1/3} t}
	\mathbbm{f}_2\|_{L^2L^2}^2
	\big).
	\label{EQ102b}
	\end{split}
\end{align}
It remains to estimate the first term on the left-hand side of \eqref{EQ102c}. 
We perform an $L^2$ energy estimate as in \cite[Proposition~6.7]{CLWZ20}, where a mollified version $\wt \rho_k$ of the function $\rho_k$ is introduced in \cite[(6.17)]{CLWZ20}. 
%These functions behaves like $L(1-|y|)$ near the boundary of the domain. 
The only modification lies in the additional term in the energy estimate, which can be estimated as 
\begin{align}
\begin{split}
	\lf|\int \wt \rho_k^{3}\omega_{k}i \mathcal{U}''\overline{\psi_k} dy\rg|
	&= 
	\lf|\int \wt \rho_k^{3}\de_k\psi_{k}  \mathcal{U}''\overline{\psi_k} dy\rg|
	%\\&
	\lesssim \|\mathcal{U}\|_{W^{3,\infty}}|k|\|\na_k\psi_k\|_{L^2}^2.
\end{split}
\end{align}
We integrate in time to get
\begin{align}
\begin{split}
	&\int_0^\infty \|\mathcal{U}\|_{W^{3,\infty}}|k|\|\na_k\psi_k\|_{L^2}^2 dt
	\\&\lesssim
	|\log \nu|^2
	\big(
	\nu^{1/3}|k|^{-4/3}\|\pa_y\omega_{{\rm in};k}\|_{L^2}^2
%	\\&\quad\quad
	+
	\|\omega_{{\rm in};k}\|_{H_k^{2}}^2
	+
	\nu^{-1/3}|k|^{4/3}\|\mathbbm{f}_1\|_{L^2L^2}^2
	+
	\nu^{-1}\|\mathbbm{f}_2\|_{L^2L^2}^2
	\big).
\end{split}
\end{align}
Hence, we obtain
\begin{align}
\begin{split}
	&
	\Vert 	e^{\delta(\nu k^2)^{1/3} t}
	 \rho_k^{3/2} \omega_k\Vert_{L^\infty L^2}^2
	 \\
	&\lesssim
	|\log \nu|^2
	\big(
	k^{-2}
	\|\omega_{{\rm in};k}\|_{H_k^{2}}^2
	+
	\nu^{-1/3}|k|^{4/3}
	\|	e^{\delta(\nu k^2)^{1/3} t}
	\mathbbm{f}_1\|_{L^2L^2}^2
	+
	\nu^{-1}
	\|	e^{\delta(\nu k^2)^{1/3} t}
	\mathbbm{f}_2\|_{L^2L^2}^2
	\big).
	\label{EQ102a}
\end{split}
\end{align}
Combining \eqref{EQ102b} with \eqref{EQ102a}, we conclude the proof of the proposition.
\end{proof}

Finally, we are ready to prove  Proposition~\ref{pro:lin_st:calU}.

\begin{proof}[Proof of Proposition \ref{pro:lin_st:calU}.]
	It remains to estimate the term
	$\Vert e^{\delta \nu^{1/3}t} (1-|y|)^{1/2} \omega_k\Vert_{L^\infty L^2}$.
The proof is similar to \cite[Proposition~6.1]{CLWZ20} by appealing to Proposition~\ref{pro_lin}, thus we omit the details.
%
%\myr{HS: It seems that we will need to add the $e^{\delta\nu^{1/3}t}$ weight to  Proposition \ref{pro_lin} and Proposition \ref{pro:lin_st_hi}. }
%\linfeng{need to check prop~\ref{pro_lin}; prop~\ref{pro:lin_st_hi} is good}
\end{proof}

\subsection{The Frozen-time Scheme}\label{sec:frozen}
In this section, we prove Proposition \ref{pro:lin_low}. We will derive the following more refined version:
\begin{proposition}
	\label{pro:lin_dy} 
		Let $\omega_k$ be the solution of \eqref{lnr_Ns_k} and $U(t,y)$ be the heat extension defined in \eqref{def:barU}. 
		For any $\delta_\ast \in(0, \delta/2)$ and
	 $\nu|k|^2 \leq \inf_{t\geq 0}\|U'(t,\cdot)\|_{L_{y}^\infty}$, we have
		\begin{align}
	\begin{split}
			&\|
			e^{\delta_\ast\nu^{1/3}t}(1-|y|)^{1/2}\omega_k\|_{L^\infty L^2}^2
			+
			\nu^{1/3}|k|^{2/3}
			\|e^{\delta_\ast\nu^{1/3}t}\rho_k^{1/2}\omega_k\|_{L^2L^2}^2
			\\&\qquad
			+
			|k|^2\|e^{\delta_\ast\nu^{1/3}t}u_k\|_{L^2L^2}^2
			+
			\nu^{1/2}|k|\|e^{\delta_\ast\nu^{1/3}t}\omega_k\|_{L^2L^2}^2
			\\&\quad
			\lesssim_{\delta_\ast} 
			|\log \nu|^2
			\lf(
			|k|^{-2}
			\|\omega_{{\rm in};k}\|_{H_k^{2}}^2\rg. 
%			\\&\hspace{2.2cm}
			\lf.
			+\nu^{-1/3}|k|^{4/3}
			\|e^{\delta_\ast\nu^{1/3}t}\mathbbm{f}_1\|_{L^2L^2}^2+
			\nu^{-1}\|e^{\delta_\ast\nu^{1/3}t}\mathbbm{f}_2\|_{L^2L^2}^2\rg).
	\end{split}
			\label{Lnr_dy_est}
	\end{align}
\end{proposition} 

Proposition \ref{pro:lin_low} is a direct consequence of Proposition~\ref{pro:lin_dy} by taking $\zeta=\delta_\ast$ in \eqref{Lnr_dy_est}.

\begin{proof}[Proof of Proposition~\ref{pro:lin_dy}]
\textbf{The $X_{\jj}$ estimates}.
We recall that $X_{[0]}$ is defined in \eqref{Frz_X_norm}.
Since $\omega_{[0]}$ is a solution of the frozen-time system \eqref{Frz_om_0},
we use Proposition \ref{pro:lin_st:calU} to obtain
\begin{align}
\begin{split}
		X_{[0]}^2
	&
	\lesssim  
	|\log \nu|^2 
	\lf( \|\omega_{{\rm in};k}\|_{H_k^{2}}^2 +
	\nu^{1/3}|k|^{-4/3}\|\pa_y\omega_{{\rm in};k}\|_{L^2}^2 
	\rg. 
	\\&\quad
	\lf.
	+
	\nu^{-1/3}|k|^{4/3}
	\|e^{\delta\nu^{1/3}t}(\mathbbm{f}_1+\mathbbm{f}_{{\rm Disc} [0]})\|_{L^2(\mathcal{I}_{[0]};L^2)}^2
	+
	\nu^{-1}
	\|e^{\delta\nu^{1/3}t}\mathbbm{f}_2\|_{L^2(\mathcal{I}_{[0]};L^2)}^2 \rg)
	\\&
	\lesssim  
	|\log \nu|^2
	 \lf( 
	 |k|^{-2}
	 \|\omega_{{\rm in};k}\|_{H_k^{2}}^2 
%	 +
%	 \nu^{1/3}|k|^{-4/3}
%	 \|\pa_y\omega_{{\rm in};k}\|_{L^2}^2 
%	 \rg. 
%	\\&\quad
%	\lf.
	+
	\nu^{-1/3}|k|^{4/3}
	\| \mathbbm{f}_1+\mathbbm{f}_{{\rm Disc} [0]} \|_{L^2(\mathcal{I}_{[0]};L^2)}^2 +
	\nu^{-1}
	\| \mathbbm{f}_2\|_{L^2(\mathcal{I}_{[0]};L^2)}^2 
	\rg).
\end{split}
\label{EQ123a}
\end{align} 
%Here, we observe the one of the key advantage of decomposing the forcing into disjoint intervals $\mathcal{I}_{\jj}$. Since the forces are localized in time, the $e^{\delta\nu^{1/3}t}$-time weight  that we pay to get the enhanced dissipation is nullified on these time intervals $\mathcal{I}_\jj$. This effect becomes apparent when we consider $X_{\jj}$ components, where $j\geq 1$.
Similarly,
for $j\geq 1$,
since $\omega_\jj$ is a solution of \eqref{Frz_om_j}, Proposition~\ref{pro:lin_st:calU} implies
\begin{align}
\begin{split}
		X_{[j]}^2
	&
	\lesssim  
	|\log \nu|^2 
	\lf(\nu^{-1/3}|k|^{4/3}\| \mathbbm{f}_1+\mathbbm{f}_{{\rm Disc} [j]} +\mathbbm{f}_{{\rm Frozen} [j]}\|_{L^2(\mathcal{I}_\jj;L^2)}^2 
	+
	\nu^{-1}\| \mathbbm{f}_2\|_{L^2(\mathcal{I}_\jj;L^2)}^2 
	\rg)
	\\&
	\les
	|\log \nu|^2 
	\Big(\nu^{-1/3}|k|^{4/3}
	\| \mathbbm{f}_1+\mathbbm{f}_{{\rm Disc} [j]}\|_{ L^2(\mathcal{I}_\jj;L^2)}^2
	+
	\nu^{-1}\| \mathbbm{f}_2\|_{L^2(\mathcal{I}_\jj;L^2)}^2 
	+
	\nu^{1/3} Y_{\jj}^2\Big),
\end{split}
	\label{Frz_X_j_est}
\end{align}
where we used Lemma~\ref{lem:Frzfrc} and $\nu k^2 \les 1$ in the last step.
%\begin{align}
%\begin{split}
%		X_{[j]}
%		&\lesssim  
%		|\log \nu| \Big(\nu^{-1/6}|k|^{2/3}\| \mathbbm{f}_1\|_{ L^2(\mathcal{I}_\jj;L^2)}+\nu^{-1/6}|k|^{2/3}\|\mathbbm{f}_{{\rm Disc}[j]}\|_{L^2(\mathcal{I}_\jj;L^2)} 	
%	\\&\quad
%	+
%	\nu^{-1/2}\| \mathbbm{f}_2\|_{L^2(\mathcal{I}_\jj;L^2)}
%	+ {\myr{\nu^{1/2}} |k|^{2/3}} {\sum_{j'=0}^{j-1}(|j-j'|+1)\|\omega_{\jjj}\|_{L^2(\mathcal{I}_{\jj};\mathcal{Z})}} \Big)
%	\\&
%	\lesssim  
%	|\log \nu| 
%	\Big(\nu^{-1/6}|k|^{2/3}\| \mathbbm{f}_1+\mathbbm{f}_{{\rm Disc} [j]}\|_{ L^2(\mathcal{I}_\jj;L^2)}
%	+\nu^{-1/2}\| \mathbbm{f}_2\|_{L^2(\mathcal{I}_\jj;L^2)} 
%	+
%	\nu^{1/6} Y_\jj\Big),
%\end{split}
%	\label{Frz_X_j_est} 
%\end{align}
%for $j\in\mathbb{N}$.
%Here, we have used the constraint $\nu|k|^{2}\lesssim 1$.

\noindent
\textbf{The $Y_\jj$-estimates.} First of all, we observe from the definition \eqref{Frz_X_norm} that 
%Therefore, we have the following estimate
\begin{align}
	\| \omega_{\jjj}\|_{L_t^2(\mathcal{I}_{\jj};\mathcal{Z})}
	\leq 
	e^{-\delta |j-j'|}
	X_{\jjj},\quad \text{for all~} j'\leq j,\label{Frz_X_est1}
\end{align}
(and of course, $\| \omega_{\jjj}\|_{L_t^2(\mathcal{I}_{\jj};\mathcal{Z})} = 0$ if $j' > j$). For $j\geq 0$,
using \eqref{Frz_X_est1} and the Cauchy-Schwarz inequality, we obtain
\begin{align}
\begin{split}
		Y_\jj^2
	&=
	\lf(\sum_{j'=0}^j(j-j'+1)\|\omega_{\jjj}\|_{L^2(\mathcal{I}_{\jj};\mathcal{Z})}\rg)^2
%	\\&
	\lesssim \lf(\sum_{j'=0}^j(j-j'+1)e^{-\delta|j-j'|}X_{\jjj}\rg)^2
	\\&
	\lesssim
	\lf(\sum_{j'=0}^j 
	(j-j'+1)^2
	e^{-2\delta (j-j')+ 7\delta_\ast (j-j')/3}
	\rg)
	\left(\sum_{j'=0}^j e^{-7\delta_\ast (j-j')/3} X_{\jjj}^2, \right)
	\\&
	\les
	\sum_{j'=0}^j e^{-7\delta_\ast (j-j')/3} X_{\jjj}^2.
\end{split}
\label{EQ123f}
\end{align}
%\myr{where we have taken $\delta_\ast \in (0,\delta/12)$ in the last step.}
From the previous estimates \eqref{EQ123a}, \eqref{Frz_X_j_est}, and \eqref{EQ123f} it follows that
\begin{align}
\begin{split}
		Y_\jj^2
		\lesssim 
		&|\log \nu|^2\sum_{j'=0}^j e^{-7\delta_\ast (j-j')/3}  \Big[\nu^{-1/3}|k|^{4/3}\| \mathbbm{f}_1+\mathbbm{f}_{{\rm Disc}[j']}\|_{L^2(\mathcal{I}_\jjj;L^2)}^2 
		+
		\nu^{-1}\| \mathbbm{f}_2\|_{L^2(\mathcal{I}_\jjj;L^2)}^2 
		\\&
		+\nu^{1/3} Y_\jjj^2
		\mathbbm{1}_{j'\geq 1}
		\Big] 
		+
		e^{-7\delta_\ast j/3}|\log \nu|^2
		|k|^{-2}		\|\omega_{{\rm in};k}\|_{H_k^{2}}^2.
%		\lf( \|\omega_{{\rm in};k}\|_{H_k^{2}}^2
%			+
%		\nu^{1/3} |k|^{-4/3}\|\pa_y\omega_{{\rm in};k}\|_{L^2}^2\rg ).
	\end{split}
	\label{Frz_Y_j_est}
\end{align}
For $\EE$ defined \eqref{Frz_E}, we use \eqref{Frz_Y_j_est} to obtain
\begin{align}
\begin{split}
		\EE&
		=
		\sum_{j=0}^N
		e^{2\delta_\ast j} Y_{\jj}^2
%	=
%	\sum_{j=0}^N e^{6\delta_\ast j}Y_\jj^2
\\ &
	\lesssim 
	\sum_{j=0}^N  |\log \nu|^2\sum_{j'=0}^j e^{-\delta_\ast (j-j')/3
		+
		2\delta_\ast j'}  \Big[\nu^{-1/3}|k|^{4/3}\| \mathbbm{f}_1+\mathbbm{f}_{{\rm Disc} [j']}\|_{L^2(\mathcal{I}_\jjj;L^2)}^2 
	\\&\quad
	+
	\nu^{-1}\| \mathbbm{f}_2\|_{L^2(\mathcal{I}_\jjj;L^2)}^2 +\nu^{1/3} Y_\jjj^2 \mathbbm{1}_{j'\geq 1}
	\Big] 
%	\\&\quad
	+
	\sum_{j=0}^N  
	e^{- \delta_\ast j/3}|\log \nu|^2
	|k|^{-2}		\|\omega_{{\rm in};k}\|_{H_k^{2}}^2
%	\lf( \|\omega_{{\rm in};k}\|_{H_k^{2}}^2
%	+ \nu^{1/3}|k|^{-4/3}\|\pa_y\omega_{{\rm in};k}\|_{L^2}^2\rg )
	\\&
	\lesssim 
	\sum_{j'=0}^N  |\log \nu|^2   \Big[\nu^{-1/3}|k|^{4/3}e^{ 2\delta_\ast j'} \| \mathbbm{f}_1
	+
	\mathbbm{f}_{{\rm Disc} [j']}\|_{L^2(\mathcal{I}_\jjj;L^2)}^2 
	\\&\quad
	+
	\nu^{-1} e^{2\delta_\ast j'} \| \mathbbm{f}_2\|_{L^2(\mathcal{I}_\jjj;L^2)}^2 
	+
	\nu^{1/3} e^{2\delta_\ast j'} Y_\jjj^2 
	\mathbbm{1}_{j'\geq 1}
	\Big] 
%	\\&\quad
	+
	| \log \nu|^2
	|k|^{-2}		\|\omega_{{\rm in};k}\|_{H_k^{2}}^2.
%	\lf( \|\omega_{{\rm in};k}\|_{H_k^{2}}^2+ \nu^{1/3}|k|^{-4/3}\|\pa_y\omega_{{\rm in};k}\|_{L^2}^2\rg ).
\end{split}
\label{EQ123b}
\end{align}
Using Lemma~\ref{lem:Frzfrc}, we obtain
\begin{align}
	\EE
	&\lesssim 
	|\log \nu|^2   \Big[\nu^{-1/3}|k|^{4/3} \| e^{ \delta_\ast\nu^{1/3}t}
	\mathbbm{f}_1\|_{L^2([0,T];L^2)}^2 
	+
	\nu^{-1} \|e^{ \delta_\ast \nu^{1/3}t} \mathbbm{f}_2\|_{L^2([0,T];L^2)}^2 \Big]
	\\&\quad
	+
	\nu^{1/3}|\log\nu|^2 \EE
	+
	| \log \nu|^2
	|k|^{-2}		\|\omega_{{\rm in};k}\|_{H_k^{2}}^2.
%	\lf( \|\omega_{{\rm in};k}\|_{H_k^{2}}^2
%	+ \nu^{1/3}|k|^{-4/3}\|\pa_y\omega_{{\rm in};k}\|_{L^2}^2\rg ).
%	\\&\quad
%	+
%	\sum_{j=0}^N
%	|\log \nu|^2
%	\nu^{-1/3} |k|^{4/3}
%	\Vert e^{\delta_\ast \nu^{1/3} t}
%	\mathbbm{f}_{{\rm Disc} [j]}\Vert_{L^2 (\mathcal{I}_{[j]}; L^2)}^2
	\label{E_est1}
\end{align} 
%which leads to
%\begin{align}
%\begin{split}
%		E
%	\lesssim 
%	&
%	|\log \nu|^2   \Big[\nu^{-1/3}|k|^{4/3} 
%	\| e^{\delta_\ast\nu^{1/3}t}
%	\mathbbm{f}_1
%	\|_{L^2([0,T];L^2)}^2 +\nu^{-1} \|e^{ \delta_\ast \nu^{1/3}t} \mathbbm{f}_2\|_{L^2([0,T];L^2)}^2 \Big]
%	\\&
%	+
%	| \log \nu|^2\lf( \|\omega_{{\rm in};k}\|_{H_k^{2}}^2+ \nu^{1/3}|k|^{-4/3}\|\pa_y\omega_{{\rm in};k}\|_{L^2}^2\rg )
%	\\&
%		+
%	\sum_{j=0}^N
%	|\log \nu|^2
%	\nu^{-1/3} |k|^{4/3}
%	\Vert e^{\delta_\ast \nu^{1/3} t} \mathbbm{f}_{{\rm Disc} [j]}\Vert_{L^2 (\mathcal{I}_{[j]}; L^2)}^2,
%\end{split}
%\label{E_est}
%\end{align} 
We take $\nu>0$ sufficiently small to obtain
\begin{align}
	\begin{split}
		&
%	\|e^{\delta_\ast\nu^{1/3}t }\omega\|_{L^2([0,T];\mathcal{Z})}^2
%	\lesssim 
	\EE
	\les
	|\log \nu|^2   \Big[\nu^{-1/3}|k|^{4/3} \| e^{\delta_\ast\nu^{1/3}t}\mathbbm{f}_1\|_{L^2([0,T];L^2)}^2 	
	+
	\nu^{-1} \|e^{\delta_\ast \nu^{1/3}t} \mathbbm{f}_2\|_{L^2([0,T];L^2)}^2 
		\\&\qquad
		+|k|^{-2}		\|\omega_{{\rm in};k}\|_{H_k^{2}}^2
%	+
%	 \|\omega_{{\rm in};k}\|_{H_k^{2}}^2+ \nu^{1/3}|k|^{-4/3}\|\pa_y\omega_{{\rm in};k}\|_{L^2}^2
	\Big].
	\label{Frz_L2t}
	\end{split}
\end{align} 
This concludes the proof of the estimates for the last three terms on the left-hand side of \eqref{Lnr_dy_est}.
%This concludes the $L^2_t\mathcal{Z}$-estimate. 

\noindent
\textbf{The $L^\infty_t$ estimates.} 
It remains to bound the first term on the left-hand side of \eqref{Lnr_dy_est}.
%\begin{align}
%F:=\lf\|\lan \nu^{1/3}t\ran^{-1/2} e^{3\delta_\ast \nu^{1/3}t} (1-|y|)^{1/2}w\rg\|_{L_t^\infty L^2_y}.
%\end{align}
We note that $\omega_{[j']}$ is a solution of \eqref{Frz_om_0} if $j'=0$ and is a solution of \eqref{Frz_om_j} if $j'\geq 1$. The space-time estimates for the time-independent flow problem   
(Proposition \ref{pro:lin_st:calU}) and the $\mathbbm{f}_{{\rm Frozen}[j]}$-estimate in  Lemma~\ref{lem:Frzfrc} imply
\begin{align}
\begin{split}
		&\|
		e^{\delta (j-j')}
		(1-|y|)^{1/2}\omega_{[j']}\|_{L^\infty(\mathcal{I}_{[j]};L^2)}
		\\
&\lesssim  
|\log \nu| 
\Big(
\nu^{-1/6}|k|^{2/3} \| \mathbbm{f}_1
+
\mathbbm{f}_{{\rm Disc} [j']}+
\mathbbm{1}_{j'\geq 1} \mathbbm{f}_{{\rm Frozen} [j']} \|_{L^2(\mathcal{I}_{[j']};L^2)} +\nu^{-1/2} \|\mathbbm{f}_2\|_{L^2(\mathcal{I}_{[j']};L^2)}
\\&\quad
+ |k|^{-1}		\|\omega_{{\rm in};k}\|_{H_k^{2}}
%( \|\omega_{{\rm in};k}\|_{H_k^{2}} + \nu^{1/6}|k|^{-2/3}\|\pa_y\omega_{{\rm in};k}\|_{L^2} 
%)
\mathbbm{1}_{j'=0}
\Big)
	\\
	&\lesssim  
	|\log \nu| 
	\Big(
	\nu^{1/6} Y_{\jjj} 
	\mathbbm{1}_{j'\geq 1}
	+
	\nu^{-1/6}|k|^{2/3} \| \mathbbm{f}_1
	+
	\mathbbm{f}_{{\rm Disc} [j']} \|_{L^2(\mathcal{I}_{[j']};L^2)} +\nu^{-1/2} \|\mathbbm{f}_2\|_{L^2(\mathcal{I}_{[j']};L^2)}\\
	&\quad
	+ 
	|k|^{-1}		\|\omega_{{\rm in};k}\|_{H_k^{2}}
%	( \|\omega_{{\rm in};k}\|_{H_k^{2}} + \nu^{1/6}|k|^{-2/3}\|\pa_y\omega_{{\rm in};k}\|_{L^2} 
%	)
	\mathbbm{1}_{j'=0}
		\Big),
\end{split}
\label{EQ123e}
\end{align}
for $0\leq j'\leq j$.
Hence, using \eqref{EQ123e} and the decomposition of solution $\omega_k$ in \eqref{Frz_om_dmp}, we obtain
\begin{align}
	&\| {e^{\delta_\ast \nu^{1/3}t} }(1-|y|)^{1/2}\omega_k\|_{L^\infty(\mathcal{I}_{\jj};L^2)}
%	\les 
%	\frac{e^{\delta_\ast j}}{1+ j^{1/2}}\|(1-|y|)^{1/2}\omega_k
%	\|_{L^\infty(\mathcal{I}_{[j]};L^2)}
	\\&
	\lesssim 
	\sum_{j'=0}^j 
	e^{\delta_\ast j-\delta (j-j')}
	\|e^{\delta (j-j')}(1-|y|)^{1/2}\omega_\jjj\|_{L^\infty(\mathcal{I}_{[j]};L^2)}
	\\&
	\lesssim
	|\log \nu|
	\sum_{j'=0}^j
	e^{(\delta_\ast - \delta) (j-j')}
	e^{\delta_\ast j'}
	\big[\nu^{1/6}Y_\jjj 
	+
	\nu^{-1/6} |k|^{2/3} \| \mathbbm{f}_1+\mathbbm{f}_{{\rm Disc} [j']}\|_{L^2(\mathcal{I}_{[j']};L^2)} 
	\\&\quad
	+\nu^{-1/2} \|\mathbbm{f}_2\|_{L^2(\mathcal{I}_{[j']};L^2)}
	\big]
	+
	e^{(\delta_\ast -\delta) j}| \log \nu|
	|k|^{-1}		\|\omega_{{\rm in};k}\|_{H_k^{2}}
%	\lf(   \|\omega_{{\rm in};k}\|_{H_k^{2}} 
%	+ \nu^{1/6}|k|^{-2/3}\|\pa_y\omega_{{\rm in};k}\|_{L^2} \rg) 
	.
\end{align}
%Let $\delta_\ast=\delta/12$. Then we have
%\begin{align}
%	&\|
%	e^{\delta_\ast \nu^{1/3}t} (1-|y|)^{1/2}w_k
%	\|_{L^\infty((0,T);L^2)} \\
%	&\lesssim
%	|\log \nu|
%	\sum_{j'=0}^j
%	e^{(\delta_\ast - \delta) (j-j')}
%	e^{\delta_\ast j'}
%	\big[\nu^{1/6}Y_\jjj 
%	+
%	\nu^{-1/6} |k|^{2/3} \| \mathbbm{f}_1+\mathbbm{f}_{{\rm Disc} [j']}\|_{L^2(\mathcal{I}_{[j']};L^2)} 
%	\\&\quad
%	+\nu^{-1/2} \|\mathbbm{f}_2\|_{L^2(\mathcal{I}_{[j']};L^2)}
%	\big]
%	+
%	e^{-\delta_\ast j}| \log \nu|\lf(   \|\omega_{{\rm in};k}\|_{H_k^{2}} 
%	+ \nu^{1/6}|k|^{-2/3}\|\pa_y\omega_{{\rm in};k}\|_{L^2} \rg) 
%	.
%\end{align}
We observe that 
$e^{\delta_\ast j'}\|\mathbbm{f}\|_{L^2(\mathcal{I}_{\jjj};L^2)}\approx \|e^{\delta_\ast \nu^{1/3}t}\mathbbm{f}\|_{L^2(\mathcal{I}_\jjj;L^2)}$ and the sum $\sum_{j'\leq j}e^{2(\delta_\ast-\delta)(j-j')}\leq C_\delta$. Applying the Cauchy-Schwarz inequality and the $\mathbbm{f}_{{\rm Disc}[j]}$-estimate in Lemma~\ref{lem:Frzfrc} yields that
\begin{align} 
	&\|e^{\delta_\ast \nu^{1/3}t} (1-|y|)^{1/2}\omega_k
	\|_{L^\infty((0,T);L^2)}
	\les
	\sup_{0\leq j\leq N}
	\|
	e^{\delta_\ast \nu^{1/3}t} (1-|y|)^{1/2}\omega_k
	\|_{L^\infty(\mathcal{I}_{[j]};L^2)} \\
	&\lesssim 
	|\log\nu|
	\lf(
	\nu^{-1/6} 
	|k|^{2/3} \| e^{\delta_\ast\nu^{1/3} t}\mathbbm{f}_1\|_{L^2([0,T];L^2)} 
	+
	\nu^{-1/2} 
	\|e^{\delta_\ast\nu^{1/3} t}\mathbbm{f}_2\|_{L^2([0,T];L^2)} 
	\rg.
		\\&\quad
	+
	\lf.  |k|^{-1}		\|\omega_{{\rm in};k}\|_{H_k^{2}}
%	 \|\omega_{{\rm in};k}\|_{H_k^{2}} + \nu^{1/6}|k|^{-2/3}\|\pa_y\omega_{{\rm in};k}\|_{L^2} 
	 \rg) 
	+|\log \nu| \nu^{1/6}\sqrt{\EE}
	. 
\end{align}
Using \eqref{Frz_L2t} and taking $\nu>0$ sufficiently small, we complete the proof of the proposition.
\end{proof}

\subsection{Estimates in the High Frequency Regime}\label{sec:high}
%
%\begin{proposition}\label{pro:lin_st_hi}
%Let $w_k$ be the solution of \eqref{lnrzd_NS}. 
%For $\nu|k|^2\ge \|\U'\|_{L^\infty}$, we have
%	\begin{align}
%	\begin{split}
%		&
%		k^2 \Vert e^{\nu k^2 t/2} u_k
%		\Vert_{L^\infty L^2}^2 
%		+
%		\nu k^2 
%		\Vert 
%		e^{\nu k^2 t/2}
%		w_{k} \Vert_{L^2 L^2}^2\ + 
%		k^2 
%		\Vert e^{\nu k^2 t/2}
%		u_k\Vert_{L^2L^2}^2 + 
%		\Vert e^{\nu k^2 t/2}
%		w_k\Vert_{L^\infty L^2}^2
%		  \\&\quad
%		  \lesssim 
%		  \nu^{-1}
%		 \left(
%		  \Vert
%		  e^{\nu k^2 t/2} \mathbbm{f}_1 \Vert_{L^2 L^2}^2 +  
%		  \Vert
%		  e^{\nu k^2 t/2} \mathbbm{f}_2
%		  \Vert_{L^2 L^2}^2\right)
%		 + \enorm{w_{{\rm in}, k}}^2.
%		 \label{sp:t:hi}
%	\end{split}
%	\end{align}
%\end{proposition}
In this section, we prove Proposition \ref{pro:lin_st_hi}, which captures the high mode dynamics.

\begin{proof}[Proof of Proposition~\ref{pro:lin_st_hi}]
	We take the $L^2$ inner product of \eqref{lnr_Ns_k}$_1$ with $-\psi_k$ to obtain
	\begin{align}
		\langle
		\partial_t \omega_k
		-
		\nu (\pa_y^2 -k^2)
		\omega_k
		+
		i U k\omega_k
		-
		U'' ik\psi_k,
		-\psi_k
		\rangle
		=\langle
		ik \mathbbm{f}_{1;k}
		+
		\partial_y \mathbbm{f}_{2;k},
		\psi_k
		\rangle,
	\end{align}
	which leads to
	\begin{align}
		\langle
		\partial_t u_k, u_k
		\rangle
		+
		\nu \Vert \omega_k\Vert_{L^2}^2
		+
		\langle
		iU k \omega_k,
		-\psi_k
		\rangle
		+
		ik \int_{-1}^{1}
		U'' |\psi_k|^2
		=\langle
		ik \mathbbm{f}_{1;k}
		+
		\partial_y  \mathbbm{f}_{2;k},
		\psi_k
		\rangle.
	\end{align}
	We take the real part of the above identity to get
	\begin{align}
		\frac{1}{2}
		\frac{d}{dt}
		\enorm{u_k}^2 
		+
		\nu \enorm{\omega_{k}}^2 
		=
		{\rm Re}
		\int_{-1}^{1} U
		ik\omega_{k} \overline{\psi_{k}} \,dy
		+
		{\rm Re}
		\int_{-1}^{1} 
		(ik \mathbbm{f}_{1;k}
		+
		\pa_y  \mathbbm{f}_{2;k}
		)\overline{\psi_{k}}\,dy.
		\label{EQ129b}
	\end{align}
%	\siming{Direct computation yields that
%	\begin{align*} 
%	&\|w_k\|_{L^2}^2={\rm Re}\int (\pa_y^2-k^2)\psi\overline{(\pa_y^2-k^2)\psi}dy={\rm Re}\int |\pa_y^2\psi|^2-2k^2\pa_y^2\psi\overline{\psi}+|k|^4|\psi|^2dy\\
%	&=\|\pa_y^2\psi_k\|_{L^2}^2+2k^2\|\pa_y\psi_k\|_{L^2}^2+k^4\|\psi_k\|_{L^2}^2.
%\end{align*}	
%As a consequence, we have that 
%\begin{align*}
%\nu\|w_k\|_{L^2}^2= \frac{1}{8}\nu\|w_k\|_{L^2}^2+\frac{7}{8}\nu\|\pa_y^2\psi_k\|_{L^2}^2+\frac{7}{4}\nu\|\pa_y k\psi_k\|_{L^2}^2+\frac{7}{8}\nu\|k^2\psi_k\|_{L^2}^2.
%\end{align*}	}
	For the first term on the right-hand side of \eqref{EQ129b}, we integrate by parts to get
	\begin{align}
		\begin{split}
			\abs{	{\rm Re}
				\int_{-1}^{1} U
				ik\omega_{k} \overline{\psi_{k}} \,dy }
			&= 
			\abs{{\rm Re}~
				ik 
				\int_{-1}^{1} 
				U' \pa_y \psi_{k} \overline{\psi_{k}} \,dy} 
			%\\&
			\leq \norm{U'}_{L^\infty_{t,y}} |k|\enorm{\pa_y \psi_{k}} \enorm{\psi_{k}} 
			\\&
			\le %\frac{1}{2}
		\frac{4}{3} \nu k^2 \enorm{\pa_y \psi_{k}}^2  
			+ 
			%\frac{1}{2}
			\frac{3}{4}\nu k^4 \enorm{\psi_{k}}^2,
		\end{split}
		\label{EQ129c}
	\end{align}
	where we used $\frac{1}{2}\Vert U'\Vert_{L^\infty_{t,y}} \leq \nu k^2$ in the last step.
	For the second term on the right-hand side of \eqref{EQ129b}, we again integrate by parts to obtain
	\begin{align}
		\begin{split}
			\left|
			{\rm Re}
			\int_{-1}^{1} 
			(ik \mathbbm{f}_{1;k}
			+
			\pa_y \mathbbm{f}_{2,k})\overline{\psi_{k}}\,dy 
			\right|
			&
			\le  
			\frac{2}{\nu k^2}\paren{\enorm{\mathbbm{f}_{1;k}}^2 
				+  
				\enorm{\mathbbm{f}_{2;k}}^2} 
			\\&\quad
			+ 
			\frac{1}{8} \nu k^2 \paren{\enorm{\pa_y \psi_{k}}^2  +  k^2 \enorm{\psi_{k}}^2}.
		\end{split}
		\label{EQ129a}
	\end{align}
	From \eqref{EQ129b}, \eqref{EQ129c}, \eqref{EQ129a}, and \eqref{identity}, it follows that
	\begin{align*}
		\frac{d}{dt}
		\enorm{u_k}^2 
		+
		\nu \enorm{\omega_{k}}^2  
		\lesssim
		\frac{1}{\nu k^2} \paren{\enorm{ \mathbbm{f}_{1;k} }^2 
			+
			\enorm{\mathbbm{f}_{2;k}}^2}.
	\end{align*}
	We integrate in time to obtain
	\begin{align}
		\begin{split}
			&
			e^{2\aaa \nu^{1/3} t}
			\enorm{u_k}^2 
			+
			\nu  \int_{0}^{t}
			e^{2\aaa \nu^{1/3} s}
			\enorm{\omega_{k}}^2ds  
			\\&\quad
			\leq
			2 \aaa \nu^{1/3}
			\int_0^t
			e^{2\aaa \nu^{1/3} s}
			\Vert u_k\Vert_{L^2}^2 
			+ 
			\frac{C}{\nu k^2}\int_{0}^{t}
			e^{2\aaa \nu^{1/3} s} \paren{\enorm{ \mathbbm{f}_{1;k} }^2 
				+  \enorm{\mathbbm{f}_{2;k}}^2} 
%			\\&\quad
			+  
			C
			\enorm{u_{{\rm in},k}}^2.
		\end{split}
	\label{EQ141a}
	\end{align}
	Since $k^2 \Vert u_k\Vert_{L^2}^2 \leq \Vert \omega_k\Vert_{L^2}^2$, we infer that
	\begin{align}
		\begin{split}
			2\aaa \nu^{1/3}
			\int_0^t
			e^{2\aaa \nu^{1/3} s}
			\Vert u_k\Vert_{L^2}^2 ds
			\leq	
			\frac{	\nu}{2}  \int_{0}^{t}
			e^{2\aaa \nu^{1/3} s}
			\enorm{\omega_{k}}^2 ds,
		\end{split}
		\label{EQ141b}
	\end{align}
	by using $\nu k^2 \geq 1/4$.
	Combining \eqref{EQ141a} with \eqref{EQ141b}, we obtain
	\begin{align}
		\begin{split}
			&
			k^2	
			\Vert e^{\aaa \nu^{1/3} t}
			u_k\Vert_{L^\infty L^2}^2
			+
			\nu k^2
			\Vert e^{\aaa \nu^{1/3} t} 
			\omega_k\Vert_{L^2 L^2}^2
			\\&\qquad
			\les
			\nu^{-1}
			(
			\Vert	e^{\aaa \nu^{1/3} s} \mathbbm{f}_{1;k} \Vert_{L^2 L^2}^2 
			+  
			\Vert e^{\aaa \nu^{1/3} s} f^2_k\Vert_{L^2 L^2}^2)
			+  
			\enorm{\omega_{{\rm in},k}}^2.
		\end{split}
	\end{align}
	Using $k^2 \Vert u_k\Vert_{L^2}^2 \les \Vert \omega_k\Vert_{L^2}^2 \les \nu k^2 \Vert \omega_k\Vert_{L^2}$, we infer that
	\begin{align}
		\begin{split}
			&
			k^2
			\Vert	e^{\aaa \nu^{1/3} t} u_k \Vert_{L^\infty L^2}^2
			+
			\nu k^2
			\Vert e^{\aaa \nu^{1/3} t} 
			\omega_k\Vert_{L^2 L^2}^2
			+
			k^2
			\Vert e^{\aaa \nu^{1/3} t} 
			u_k\Vert_{L^2 L^2}^2
			\\&\qquad
			\les
			\nu^{-1} (
			\Vert	e^{\aaa \nu^{1/3} s} \mathbbm{f}_{1;k} \Vert_{L^2 L^2}^2 
			+  
			\Vert e^{\aaa \nu^{1/3} s} \mathbbm{f}_{2;k}
			\Vert_{L^2 L^2}^2)
			+  
			\enorm{\omega_{{\rm in},k}}^2.
		\end{split}
		\label{EQ130b}
	\end{align}

	It remains to estimate the term
	$\Vert e^{\aaa \nu^{1/3} t} (1-|y|)^{1/2} \omega_k \Vert_{L^\infty L^2}$.
	Denote $G= \partial_t \psi_k+ik U \psi_k$.
	Direct computation shows that
	\begin{align}
		\left\{
		\begin{aligned}
			&
			\de_{k}
			G
			=
			\partial_t \omega_k
			+
			ik
			U\omega_k
			+2 ik U' \pa_y\psi_k
			+
			ik U'' \psi_k,
			\\&
			G|_{y=\pm 1}
			=
			\partial_y G|_{y=\pm 1}=0.
		\end{aligned}
		\right.
		\label{EQ130e}
	\end{align}
	We take the $L^2$ inner product of \eqref{lnr_Ns_k}$_1$ with $-G$ and use \eqref{EQ130e} to obtain
	\begin{align}
		\begin{split}
			&
			\langle
			ik \mathbbm{f}_{1;k}+\partial_y \mathbbm{f}_{2;k},
			-G
			\rangle
			=
			\langle
			\partial_t \omega_k 
			-\nu \de_{k} \omega_k
			+
			ik U
			\omega_k
			-
			U''ik \psi_k,
			-G
			\rangle
			\\&\quad
			=
			\langle
			\de_{k} G
			-
			2ik U' \pa_y \psi_k
			-
			2ik U'' \psi_k
			-\nu\de_{k} \omega_k,
			-G
			\rangle
			\\&\quad
			=
			\Vert \partial_y G\Vert_{L^2}^2
			+
			k^2 \Vert G\Vert_{L^2}^2
			+
			\nu
			\langle
			\omega_k,
			\de_{k} G
			\rangle
			+
			2ik \langle
			U' \pa_y \psi_k, G
			\rangle
			+2ik 
			\langle U'' \psi_k, G
			\rangle
			\\&\quad
			=
			\Vert \partial_y G\Vert_{L^2}^2
			+
			k^2 \Vert G\Vert_{L^2}^2
			+
			\nu
			\langle
			\omega_k,
			\partial_t \omega_k
			+
			ik
			U\omega_k
			+2 ik U' \pa_y \psi_k
			+
			ik U'' \psi_k
			\rangle
			\\&\quad\quad\quad
			+
			2ik \langle
			U' \pa_y \psi_k, G
			\rangle
			+2ik 
			\langle U'' \psi_k, G
			\rangle.
		\end{split}
	\end{align}
	Taking the real part, we get
	\begin{align*}
		\begin{split}
			&\frac{\nu}{2}
			\frac{d}{dt}
			\Vert \omega_k\Vert_{L^2}^2
			+
			\Vert \partial_y G\Vert_{L^2}^2
			+
			k^2 \Vert G\Vert_{L^2}^2
			\\&
			\lesssim
			\nu |k|
			\Vert \omega_k\Vert_{L^2}
			\Vert \pa_y \psi_k\Vert_{L^2}
			+
			\nu |k| 
			\Vert \pa_y \psi_k\Vert_{L^2}
			\Vert \psi_k\Vert_{L^2}
			+
			|k|
			\Vert \pa_y \psi_k \Vert_{L^2}
			\Vert G\Vert_{L^2}
			\\&\quad\quad
			+
			|k| 
			\Vert \psi_k\Vert_{L^2}
			\Vert G\Vert_{L^2}
			+
			|k| \Vert \mathbbm{f}_{1;k} \Vert_{L^2} \Vert G\Vert_{L^2}
			+
			\Vert \mathbbm{f}_{2;k} \Vert_{L^2}
			\Vert \partial_y G\Vert_{L^2}
			\\&
			\lesssim
			\nu
			\Vert \omega_k\Vert_{L^2}^2
			+
			\nu k^2 
			\Vert \pa_y \psi_k\Vert_{L^2}^2
			+
			|k|
			\Vert \pa_y \psi_k\Vert_{L^2}
			\Vert G\Vert_{L^2}
			+
			|k|
			\Vert \psi_k\Vert_{L^2}
			\Vert G\Vert_{L^2}
			\\&\quad\quad
			+
			|k|
			\Vert \mathbbm{f}_{1;k}\Vert_{L^2}
			\Vert G\Vert_{L^2}
			+
			\Vert \mathbbm{f}_{2;k} \Vert_{L^2}
			\Vert \partial_y G\Vert_{L^2}
			,
		\end{split}
	\end{align*}
	where we used Lemma~\ref{lem:Frz_cf} and \eqref{identity}.
	Therefore, using the Young and Poincar\'{e} inequalities, we obtain
	\begin{align}
		\begin{split}
			\nu
			\frac{d}{dt}	
			\Vert \omega_k \Vert_{L^2}^2
			+
			\Vert \partial_y G\Vert_{L^2}^2
			&
%			\lesssim
%			\Vert \mathbbm{f}_{1;k} \Vert_{L^2}^2
%			+
%			\Vert \mathbbm{f}_{2;k}
%			\Vert_{L^2}^2
%			+
%			\nu \Vert \omega_k\Vert_{L^2}^2
%			+
%			\nu k^2 
%			\Vert \pa_y \psi_k\Vert_{L^2}^2
			%		\\&
			%		\lesssim
			%		\Vert f^1_k \Vert_{L^2}^2
			%		+
			%		\Vert f^2_k\Vert_{L^2}^2
			%		+
			%		\nu \Vert \omega_k\Vert_{L^2}^2
			%		+
			%		\nu k^2
			%		\Vert \psi'_k\Vert_{L^2}^2
%			\\&
			\lesssim
			\Vert \mathbbm{f}_{1;k} \Vert_{L^2}^2
			+
			\Vert \mathbbm{f}_{2;k} \Vert_{L^2}^2
			+
			\nu \Vert \omega_k\Vert_{L^2}^2.
		\end{split}
	\end{align}
	We integrate in time and use \eqref{EQ130b} to obtain
	\begin{align}
		\begin{split}
			\Vert
			e^{\aaa \nu^{1/3} t}
			\omega_k \Vert_{L^\infty L^2}^2
			&
			\lesssim
			\nu^{-1}	
			\Vert
			e^{\aaa \nu^{1/3} t} \mathbbm{f}_{1;k}
			\Vert_{L^2 L^2}^2
			+
			\nu^{-1}	
			\Vert
			e^{\aaa \nu^{1/3} t} \mathbbm{f}_{2;k}
			\Vert_{L^2 L^2}^2
			+
			\Vert \omega_{{\rm in},k}\Vert_{L^2}^2.
			\label{EQ131b}
		\end{split}
	\end{align}
	The proof of the proposition is thus completed.
\end{proof}

\section{Nonlinear Estimates} \label{sec:NL}
In this section, we prove nonlinear stability. We denote
\begin{align}
	&E_k:
	=
	\| e^{\aaa \nu^{1/3} t}
	(1-|y|)^{1/2}
	 \omega_k\|_{L_t^\infty L_y^2}
	+
	\nu^{1/4}|k|^{1/2}\| e^{\aaa \nu^{1/3} t} \omega_k\|_{L_t^2L_y^2}
	+
	|k|\| e^{\aaa \nu^{1/3} t} u_k\|_{L_t^2L_y^2}
\end{align}
for $k\in \ZZ \setminus \{0\}$, and
\begin{align}
	E_0:=\| \omega_0\|_{L_t^\infty L_y^2}.
\end{align}
%The parameter $\aaa\in (0,1)$ shall be determined below in the proof.
%+|k|^{1/2}\|u_k\|_{L^\infty_t L^\infty_y}
%\myr{HS: Here, we have the problem about the derivation of higher regularity norms. For example, how to derive that $\|\na_k w_k\|_{L^2}$ becomes finite after a positive time? We might need to introduce the following component into the energy functional, 
%\begin{align}\nu^{\al }\|\sqrt{\nu t} \na_k w_k\|_{L_t^\infty L^2}. \end{align}When we do the energy estimate of the above quantity, we have an extra quantity 
%\begin{align}
%\nu^{2\al }\nu\|\na_k w_k\|_{L^2}^2,
%\end{align}and maybe we can use the $\nu^{...}\|\na_k w_k\|_{L^2L^2}^2$ control to control it.  }
%We introduce the following sum of the energy $E_k$:
%\begin{align}
%	E:=\sum_{k\in\mathbb{Z}}E_k,\quad E_\nq:=\sum_{k\in \mathbb{Z}\setminus \{0\}}E_k.
%\end{align}
%Note that we are using the $\ell^1_k$ type norm here. 
%We prove the following theorem 
%\linfeng{check the high frequency case. The following the theorem seems to hold only for $\nu k^2 \les 1$.}
%\begin{theorem}
%	Suppose that $\| w_{in}\|_{H^2}\leq \ep \nu^{1/2}|\log\nu|^{-4}$. Then
%	\begin{align}
%	E\leq C\ep \nu^{1/2}|\log\nu|^{-4}.
%	\end{align}
%\end{theorem}

\begin{proof}[Proof of Theorem~\ref{T01}]
We apply the Fourier transform in the $x$-variable and
rewrite the equation \eqref{vort-orig} into the linearized form \eqref{lnr_Ns_k}, where
\begin{align}
	\label{nlnrty}
	\mathbbm{f}_{1;k}=-f^1_k:=\sum_{\ell\in \mathbb Z} u^{(1)}_\ell  \omega_{k-\ell},\quad 
	\mathbbm{f}_{2;k}
	=-f^2_k:
	=
	\sum_{\ell\in \mathbb Z} u^{(2)}_\ell \omega_{k-\ell}.
\end{align}
By the boundary condition $u^{(2)} |_{y=\pm 1}=0$ and incompressibility $\nabla \cdot u=0$,
we have $u_0^{(2)} = 0$. 
Therefore, from \eqref{EQ01}$_1$ it follows that
\begin{align}
\begin{split}
	 (\pa_t-\nu\pa_y^2)
	 u^{(1)}_0
	 =
	 &-\sum_{\ell \in \mathbb{Z}\setminus \{0\}} u^{(2)}_\ell {\pa_y u^{(1)}_{-\ell}}
	=
	-\sum_{\ell \in \mathbb{Z}\setminus \{0\}}
	u^{(2)}_\ell {(\pa_y u^{(1)}_{-\ell}
	-i\ell u^{(2)}_{-\ell})}
	\\
	=&
	\sum_{\ell \in \mathbb{Z}\setminus \{0\}}
	u^{(2)}_\ell { \omega_{-\ell}}
	=
	-
	f^2_0.
	\label{EQ124a}
\end{split}
\end{align}
Multiplying \eqref{EQ124a} by $\partial_y^2 u^{(1)}_0$ and integrating by parts, we obtain
\begin{align}
	\frac{1 }{2} \frac{d}{dt}\|\pa_y u^{(1)}_0 \|_{L^2}^2
	=
	-\nu \|\pa_y^2 u^{(1)}_0 \|_{L^2}^2
	+
	\lan f^2_0,\pa_y^2  u^{(1)}_0\ran
	\leq -\frac{1}{2}\nu
	\|\pa_y^2 u^{(1)}_0 \|
	_{L^2}^2
	+
	(2\nu)^{-1}\| f^2_0 \|_{L^2}^2,
	\label{EQ127b}
\end{align}
where we used the boundary conditions $u^{(1)}_0 |_{y=\pm 1}=0$.
Since $ \omega_0 =- \partial_y u^{(1)}_0$, integrating \eqref{EQ127b} in time yields
\begin{align}\label{E_0_est}
	E_0^2
	=
	\|  \omega_0 \|_{L^\infty_t L^2_y}^2
	\lesssim \nu^{-1}
	\| f^2_0 \|_{L^2L^2}^2
	+
	\|  \omega_{\rm in;0} \|_{L^2}^2. 
\end{align}
For $k\neq 0$, we have
\begin{align}
	\begin{split}
		&
		\lf\|
		\frac{e^{\aaa \nu^{1/3} t} u^{(2)}_k (t,y)}{(1-|y|)^{1/2}}
		\rg\|_{L^2_t L^\infty_y}^2	
		=
		\lf\|
		e^{2\aaa \nu^{1/3} t}
		\sup_{y\in [-1,1]}
		\frac{| u^{(2)}_k (t,y) |^2}{1-|y|}
		\rg\|_{L^1_t}
		\\&	\qquad
		=
		\lf\|
		e^{2\aaa \nu^{1/3} t}
		\max\Big\{
		\sup_{y\in [0,1]}
		\frac{ |\int_1^y \partial_z u^{(2)}_k (t,z)dz|^2}{1-y},
		\sup_{y\in [-1,0]}
		\frac{ |\int_{-1}^y \partial_z u^{(2)}_k (t,z)dz|^2}{1+y}
		\Big\}
		\rg\|_{L^1_t}	
		\\&\qquad
		\leq
		\Vert e^{\aaa \nu^{1/3} t}
		\partial_y u^{(2)}_k\Vert_{L^2_t L^2_y}^2
		=|k|^2 
		\Vert e^{\aaa \nu^{1/3} t} u^{(1)}_k\Vert_{L^2_t L^2_y}^2
		\leq  E_k^2,
	\end{split}
	\label{EQ124b}
\end{align}
where we used the Cauchy-Schwarz inequality.
For $k\in \mathbb Z$, combining \eqref{nlnrty} with \eqref{EQ124b} gives
\begin{align}
	\begin{split}
		&
		\| e^{\aaa \nu^{1/3} t} f^2_k\|_{L^2_t L^2_y}
		\leq \sum_{\ell\in\mathbb{Z} \setminus \{0\}}\| e^{\aaa \nu^{1/3} t}u^{(2)}_\ell  \omega_{k-\ell}\|_{L^2_t L^2_y}
			\\&\quad
		\leq
		C\sum_{\ell\in\mathbb{Z} \setminus \{0\} }
		\lf\|\frac{e^{\aaa \nu^{1/3} t} u^{(2)}_\ell}{(1-|y|)^{1/2}}\rg\|_{L^2_t L^\infty_y} 
		\|
		 (1-|y|)^{1/2} \omega_{k-\ell}\|_{L_t^\infty L_{y}^2}
		\\&\quad
		\leq C\sum_{\ell \in \mathbb{Z}\setminus \{0\}}E_\ell E_{k-\ell},
		\label{f_2_L22}
	\end{split}
\end{align}
where we used $u^{(2)}_0= 0$.
In particular, 
\begin{align}
	\begin{split}
		&
		\| f^2_0\|_{L^2_t L^2_y}
		\leq C\sum_{\ell \in \mathbb{Z}\setminus \{0\}}E_\ell E_{-\ell}.
		\label{f_2_L221}
	\end{split}
\end{align}
Combining \eqref{E_0_est} with \eqref{f_2_L221}, we obtain
\begin{align}
	E_0
	\les
	\nu^{-1/2}
	\sum_{\ell\in \mathbb{Z}\setminus \{0\}}
	E_\ell E_{-\ell}
	+\textbf{}
	\Vert  \omega_{\rm in;0}\Vert_{L^2}.
	\label{E_0_est1}
\end{align}

Next we estimate $E_k$ for $k\neq 0$.
By Proposition \ref{pro:lin_st}, for $k\neq 0$,
\begin{align}
	\begin{split}
		E_k
		&\lesssim 
		|\log\nu|
				\big(
				|k|^{-1}		\|\omega_{{\rm in};k }\|_{H_k^{2}}
%		\| \omega_{{\rm in};k}\|_{H^2_k}
%		+
%		\nu^{1/6} |k|^{-2/3}
%		\Vert \partial_y  \omega_{{\rm in};k}\Vert_{L^2}
%		)
		+
		\min\{\nu^{-1/6} |k|^{2/3}, \nu^{-1/2}\}
		\| e^{\aaa \nu^{1/3} t} f^1_k\|_{L^2 L^2}
			\\&\hspace{2cm}
		+
		\nu^{-1/2}
		\| e^{\aaa \nu^{1/3} t} f^2_k\|_{L^2L^2}
		\big).
	\end{split}
	\label{EQ125a}
\end{align} 
For $k\neq 0$, using \eqref{nlnrty} we have
\begin{align}
\begin{split}
	\|e^{\aaa \nu^{1/3} t} f^1_k\|_{L^2_t L^2_y} 
	&
	\leq
	\|u^{(1)}_0\|_{L^\infty_tL^\infty_y}
	\|e^{\aaa \nu^{1/3} t}  \omega_k\|_{L^2_tL^2_y}
	+
	\|e^{\aaa \nu^{1/3} t} u^{(1)}_k\|_{L^2_t L^\infty_y}
	\| \omega_0\|_{L^\infty_t L^2_y} 
	\\&\qquad
	+
	\sum_{\ell\in\mathbb Z\backslash\{0,k\}}
	\|e^{\aaa \nu^{1/3} t} u^{(1)}_\ell  \omega_{k-\ell}\|_{L^2_t L^2_y}\\
	&=:T_1+T_2+T_3.
\end{split}
\label{EQ124e}
\end{align}
By the Poincar\'{e} and Agmon inequalities,
\begin{align}
	\begin{split}
	\| u^{(1)}_0\|_{L_t^\infty L^\infty_y}
	\lesssim 
	\|u^{(1)}_0\|_{L^\infty_t L^2_y}^{1/2}
	\|\pa_y u^{(1)}_0 \|_{L^\infty_t L^2_y}^{1/2}
	\lesssim 
	\|\pa_y u^{(1)}_0 \|_{L^\infty_t L^2_y}
	=
	\|  \omega_0\|_{L^\infty_t L^2_y}.
	\end{split}
\end{align}
Applying Poincar\'{e} and Agmon inequalities again gives
\begin{align}
\begin{split}
	&T_1+T_2
	\lesssim 
	\|  \omega_0 \|_{L^\infty_t L_y^2} 
	\lf(\|e^{\aaa \nu^{1/3} t}  \omega_k\|_{L^2_t L_y^2}
	+
	\|\pa_y u^{(1)}_k\|_{L^2_t L_y^2}^{1/2}
	\| e^{\aaa \nu^{1/3} t} u^{(1)}_k\|_{L_t^2 L_y^2}^{1/2}\rg)
	\\&
	\lesssim 
	\|  \omega_0 \|_{L^\infty_t L_y^2} \lf(\|e^{\aaa \nu^{1/3} t}  \omega_k\|_{L_t^2 L_y^2}
	+
	(\| \omega_k\|_{L_t^2 L_y^2}
	+
	\||k| u^{(2)}_k\|_{L^2_t L_y^2})^{1/2}
	\|e^{\aaa \nu^{1/3} t} u^{(1)}_k\|_{L_t^2 L_y^2}^{1/2}\rg)
	\\&
	\lesssim 
	\|   \omega_0 \|_{L^\infty_t L_y^2} \lf(\| e^{\aaa \nu^{1/3} t}  \omega_k\|_{L_t^2 L_y^2}
	+ 
	|k|^{1/2} \|e^{\aaa \nu^{1/3} t} 
	u_k\|_{L^2_t L_y^2}\rg)
	\\&
	\lesssim 
	\nu^{-1/4}
	|k|^{-1/2} E_0E_k.
\end{split}
\label{EQ124d}
\end{align}
%Here, we recall that $(u_1)_k\big|_{y\pm 1}=0$.  
A similar argument to \eqref{EQ124b} yields
\begin{align}
	\lf \|\frac{e^{\aaa \nu^{1/3} t} u^{(1)}_m(t,y)}{(1-|y|)^{1/2}}\rg\|_{L_t^2 L_y^\infty}^2
	\leq 
	\|e^{\aaa \nu^{1/3} t}
	\pa_y u^{(1)}_m\|_{L_t^2L_y^2}^2,
	\quad
	\forall m\in\mathbb{Z}.
	\label{EQ126a}
\end{align}
%A similar estimate is derived in the \cite{CLWZ20}, but the weight is slightly different, we present the proof here
%\begin{align}
%\begin{split}
% 	 & \lf \|\frac{(u_1)_k(t,y)}{(1-|y|^2)^{1/2}}
% 	 \rg \|_{L^2L^\infty}^2=\lf\|\sup_{y\in[-1,1]}\frac{|(u_1)_k|^2}{(1-|y|^2)^{1/2}}\rg\|_{L_t^1}\\
% &\leq \lf\|\max\lf\{\sup_{y\in[0,1]}\frac{|\int_1^y \pa_z (u_1)_k(t,z)dz|^2}{1-|y|},\sup_{y\in[-1,0]}\frac{|\int_{-1}^y \pa_z (u_1)_k(t,z)dz|^2}{1-|y|}\rg\}\rg\|_{L_t^1}\\
% &\leq C\lf\|\max\lf\{\sup_{y\in[0,1]}\frac{(1-|y|)^{1/2\times 2} \|\pa_z (u_1)_k(t,z)\|_{L^2_z}^2}{1-|y|},\sup_{y\in[-1,0]}\frac{(1-|y|)^{1/2\times 2} \|\pa_z (u_1)_k(t,z)\|_{L^2_z}^2}{1-|y|}\rg\}\rg\|_{L_t^1}\\
% &\leq C\|\pa_y(u_1)_k\|_{L_t^2L_y^2}^2.
%\end{split}
%\end{align}
Therefore, 
\begin{align}
\begin{split}
	 T_3
	 %&	 = \sum_{\ell\in\mathbb Z\setminus\{0,k\}} \|e^{\aaa \nu^{1/3} t} u^{(1)}_\ell  \omega_{k-\ell}\|_{L^2_t L_y^2}	 \\
	 &
	\lesssim
	\sum_{\ell\in\mathbb Z\setminus\{0,k\}}
	\lf\|\frac{e^{\aaa \nu^{1/3} t} u^{(1)}_\ell}{(1-|y|)^{1/2}}\rg\|_{L_t^2 L_y^\infty}
	\|(1-|y|)^{1/2} \omega_{k-\ell}\|_{L_t^\infty L_y^2}
	\\&
	\lesssim 
	\sum_{\ell\in\mathbb Z\setminus\{0,k\}}
	\|e^{\aaa \nu^{1/3} t} \pa_y u^{(1)}_\ell\|_{L_t^2L_y^2}\|(1-|y|)^{1/2} \omega_{k-\ell}\|_{L_t^\infty L_y^2}
\end{split}
\end{align}
which leads to
\begin{align}
	\begin{split}
		T_3&
		\lesssim 
		\sum_{\ell\in\mathbb Z\backslash\{0,k\}}
		\lf(\|e^{\aaa \nu^{1/3} t}  \omega_\ell\|_{L_t^2 L_y^2}
		+
		\ell^2 \|e^{\aaa \nu^{1/3} t} \Delta_\ell^{-1} \omega_\ell\|_{L_t^2L_y^2}\rg)\|(1-|y|)^{1/2}
		 \omega_{k-\ell}\|_{L^\infty_t L_y^2}
		\\&
		\lesssim 
		\sum_{\ell\in\mathbb Z\setminus\{0,k\}}
		\|e^{\aaa \nu^{1/3} t}  \omega_\ell\|_{L^2_t L_y^2}\|
%		\br{\nu^{1/3}t}^{-1/2}
		e^{\aaa \nu^{1/3} t}(1-|y|)^{1/2} \omega_{k-\ell}\|_{L^\infty_t L_y^2}
		\\&
		\lesssim \sum_{\ell\in\mathbb Z\backslash\{0,k\}}\nu^{-1/4}
		|k|^{-1/2}E_\ell E_{k-\ell}.
	\end{split}
\label{EQ124c}
\end{align}

Using \eqref{f_2_L22}, \eqref{E_0_est1}, \eqref{EQ125a}, \eqref{EQ124e}, \eqref{EQ124d}, and \eqref{EQ124c}, we obtain
\begin{align}
\begin{split}
	\sum_{k\in \mathbb{Z}}
	E_k&
	\lesssim 
	|\log\nu|\sum_{k\in \mathbb Z}
	|k|^{-1}
	\| \omega_{\mathrm{in};k}\|_{H^2_k}
	+
	\nu^{-1/2}|\log\nu|
	\sum_{k\in\mathbb Z}\sum_{\ell\in \mathbb{Z}}E_{\ell} E_{k-\ell}
	\\&
	\lesssim 
	|\log\nu|\sum_{k\in \mathbb Z}	|k|^{-1}
	\| \omega_{\mathrm{in};k}\|_{H^2_k}
	+
	\nu^{-1/2}|\log\nu|
	\lf(\sum_{k\in\mathbb Z} E_{k }\rg)^2,
	\label{sum_E_k}
\end{split}
\end{align} 
%where we used $\nu k^2 \les 1$.
Noting that the initial data $\Vert  \omega_{\rm in}\Vert_{H^2} \leq \eps \nu^{1/2} |\log \nu|^{-2}$ implies 
%\siming{One needs higher regularity here. }
\begin{align}
	\label{Initial_pt}
	\sum_{k\in \mathbb Z}
		|k|^{-1}
	\| \omega_{\mathrm{in};k}\|_{H^2_k}
	\les
	\epsilon\nu^{1/2} |\log \nu|^{-2}.
\end{align}
A continuity argument yields
\begin{align}
	\sum_{k\in \mathbb{Z}}
	E_k
	\leq
	C\epsilon \nu^{1/2} |\log \nu|^{-1},
	\label{EQ127c}
\end{align}
concluding the proof of the theorem.
\end{proof}

%\newpage
%\section{Velocity Decomposition}

%\newpage
\appendix
%\section{Properties of the Airy Functions}
%\input{Airy.tex}
%\section{Time Laplace-Fourier Relations}
%\input{Laplace-Fourier.tex}
\section{Properties of the Singular Integral Operator $\mathfrak{J}_k$}
This section collects the necessary estimates related to the physical side ghost singular integral operator defined as
\begin{align}\label{J_k}
\mathfrak{J}_k[f](y) := \abs{k}\mh{^{1-\delta}} \text{p.v.} \frac{k}{\abs{k}} \int_{-1}^1 \frac{1}{2i(y-y')} G_k(y, y') f(y') \dee y', 
\end{align}
where $G_k$ is the Green's function for $\Delta_k:=-|k|^2+\pa_{yy}$ with homogeneous Dirichlet boundary conditions and has the explicit form \eqref{green}.
%\begin{align} 
%G_k(y,y')=-\frac{1}{k\sinh(2k)}\left\{\begin{array}{cc}\sinh(k(1-y'))\sinh(k(1+y)),&\quad y\leq y';\\
%\sinh(k(1-y))\sinh(k(1+y')),&\quad y\geq y'.\end{array}\right.\label{G_k}
%\end{align}

\begin{lemma}[\cite{BHIW23I}] \label{lem:BoundI}
The singular integral operator $\mathfrak{J}_k$ extends to a bounded linear operator on $L^2 \to L^2$ and moreover
\begin{align*}
\norm{\mathfrak{J}_k}_{L^2 \to L^2} \lesssim 1. 
\end{align*}
\end{lemma}

The next lemma captures the commutator between $\partial_y$ and $\mathfrak{J}_k$. 
\begin{lemma}[\cite{BHIW23I}] 
\label{lem:CommIpy}
If $f \in H^1$ then $\mathfrak{J}_k[f] \in H^1$. We define 
\begin{align*}
	 \mathfrak{H}_k:=
[\partial_y, \mathfrak{J}_k].
\end{align*}
Then
\begin{align*}
\mathfrak{H}_k[f] = |k| \textup{p.v.} \int_{-1}^1 \frac{H_k(y,y')}{2i(y-y')} f(y') \dee y', 
\end{align*}
where $H_k$ is given by 
\begin{align}
H_k(y,y') := -\frac{\sinh(k(y+y'))}{\sinh(2k)}
\end{align}
\end{lemma}

We have the following commutator estimate on $\mathfrak{H}_k$.
\begin{lemma}[\cite{BHIW23I}]
\label{lem:H}
There holds the estimate
\begin{align} \label{frakH}
\norm{\mathfrak{H}_k}_{L^2 \to L^2} \lesssim |k|. 
\end{align}
\end{lemma}

%Finally we point out the following symmetry properties. 
\begin{lemma}[\cite{BHIW23I}]
\label{lem:AS}
For all $f,g \in L^2$ there holds 
\begin{align*}
\overline{\mathfrak{J}_k[f]} 
= 
- \mathfrak{J}_k[\overline{f}]
\end{align*}
and
\begin{align}
\int_{-1}^1 \bar{f} \mathfrak{J}_k[g] dy = - \int_{-1}^1 \mathfrak{J}_k[\bar{f}] g  dy, 
\end{align}
which in particular, implies $\mathfrak{J}_k = \mathfrak{J}^\ast_k$. 
\end{lemma}

%%%%%%%%%%%%
\begin{lemma}[\cite{BHIW23I}] \label{lem:Trans}
Under the regularity hypotheses, for all $\delta>0$ sufficiently small (depending only on universal constants), the following estimate holds 
\begin{align}
\label{est_T_1}
\Re \lan -\U ik w_k,\mf J_k w_k\ran \leq -\frac{|k|^2}{16}\|\na_k\phi_k\|_{L^2}^2,
\\
|\Re\lan \U'' ik\phi_k, \mf J_k w_k\ran|\leq C\delta \|\na_k\phi_k\|_{L^2}^2,
\label{est_T_2}\\
\Re \lan \nu \de_kw_k,\mf J_k w_k\ran\leq C\nu\|\na_k w_k\|_{L^2}^2,
\label{est_T_3} 
\end{align}
where the implicit constant depends on $\|\U\|_{ C^3}.$
\end{lemma}

%\section{Properties of the Boundary Corrector Functions}
\section{Auxiliary Lemmas}
\label{Applemma}
\begin{proof}[Proof of Proposition \ref{pro_Nav}]
Without loss of generality, we assume that $k>0$. 
First, we develop the estimates for the spectral stream function $\phi_\nv$ and the spectral vorticity gradient $\na  w_\nv$ using the singular integral operator $\mf J_k$. 
Next, we use some standard resolvent estimate techniques to derive the enhanced dissipation estimates for the spectral vorticity $ w_\nv$. 
%In the final step, we fill in the remaining estimates. %employ a perturbation argument to show that the same estimates derived can be extended to the $k \operatorname{Im} \lambda \geq-\delta_0 \nu^{1 / 3}|k|^{2 / 3}$ regime.

We aim to prove the following key estimates 
	\begin{align}
		&|k|\|\na_k\phi_\nv\|_{L^2}+\sqrt{\nu}\|\na_k  w_{\rm Na}\|_{L^2}+\sqrt{|k|(\mathrm{Im}\lambda)_+}\| w_\nv \|_{L^2}\n
		\\&\quad
		\leq C \delta_0^{-\frac{1}{2}}\nu^{-\frac{1}{6}}|k|^{-\frac{1}{3}}\|F_k\|_{L^2}
		+
		C\delta_0^{\frac{1}{2}}\nu^{\frac{1}{6}}|k|^{\frac{1}{3}}\| w_\nv\|_{L^2},
		\quad
		\text{if~}F_k\in L^2,\label{na_phi_FL2}
		\\&
		|k|\|\na_k\phi_{\rm Na}\|_{L^2}+\sqrt{\nu}\|\na_k  w_{\rm Na}\|_{L^2}
		+\sqrt{|k|(\mathrm{Im}\lambda)_+}\| w_\nv \|_{L^2}\n\\
		&\quad
		\leq C
		|k|^{-1}\|\na_k F_k\|_{L^2}
		+
		C\delta_0^{\frac{1}{2}}\nu^{\frac{1}{6}}|k|^{\frac{1}{3}}\| w_{\rm Na}\|_{L^2},
		\quad
		\text{if~}
		F_k\in H_0^1,\label{na_phi_FH1}
		\\&
		|k|\|\na_k\phi_\nv\|_{L^2}+\sqrt{\nu}\|\na_k  w_{\rm Na}\|_{L^2}
		+\sqrt{|k|(\mathrm{Im}\lambda)_+}\| w_\nv \|_{L^2}\n\\
		&\quad
		\leq C\nu^{-\frac{1}{2}}
		\|F_k\|_{H^{-1}}
		+
		C\delta_0^{\frac{1}{2}}\nu^{\frac{1}{6}}|k|^{\frac{1}{3}}\| w_\nv\|_{L^2},\quad
		\text{if~}
		F_k\in H^{-1}.\label{na_phi_FH-1}
	\end{align} 
First, we multiply the equation \eqref{Nav_eqn}$_1$ by $\overline{ w_\nv}$ and take the real part to obtain 
\begin{align}
	&-\nu \int \de_k  w_\nv \overline{ w_\nv} dy+\int (\U-\lambda)ik w_\nv \overline{ w_\nv} dy
	=
	\int F_k\overline{ w_\nv}dy+ik\int\U'' \phi_\nv\overline{\de_k\phi_\nv} dy.
\end{align}
We integrate by parts to get
\begin{align}
	&
	\nu\|\na_k  w_{\rm Na}\|_{L^2}^2 +ik\int (\U-\mathrm{Re}\lambda)| w_\nv|^2dy+k\mathrm{Im}\lambda\| w_\nv\|_{L^2}^2\\ \n
	&\quad
	=
	\int F_k\overline{ w_\nv}dy-ik\int \na_k(\U'' \phi_\nv)\cdot\overline{\na_k\phi_\nv}.
\end{align}
Taking the real part, we are led to
\begin{align} 
\begin{split}
&\nu\|\na_k  w_{\rm Na}\|_{L^2}^2 +k\max\{0,\mathrm{Im}\lambda\}\| w_\nv \|_{L^2}^2
\leq 
\lf|\mathrm{Re}\lan F_k,  w_{\rm Na}\ran\rg|
\\&\quad
+
\|\U'''\|_{L^\infty}\|k\phi_\nv\|_{L^2}\|\pa_y\phi_{\rm Na}\|_{L^2}
-
k\min\{0,\mathrm{Im}\lambda\} \| w_\nv \|_{L^2}^2.
\end{split}
\label{PRO_Navpf1}
\end{align}
Thanks to the spectral constraint $k \mathrm{Im} \lambda \geq -\delta_0 \nu^{1/3} |k|^{2/3}$ and the fact that $ w_\nv\in H^1_0$, we apply interpolations based on the regularity of $F_k$ for the first term on the right side of \eqref{PRO_Navpf1}, yielding
\begin{align}
	\begin{split}
		&
	\nu\|\na_k  w_\nv\|_{L^2}^2 +|k|\max\{0,\mathrm{Im}\lambda\}\| w_\nv \|_{L^2}^2
	\leq 
	\min\big
	\{
	\|F_k\|_{L^2}\| w_\nv\|_{L^2},
	\\&\quad
	\|F_k \|_{H^{-1}}\| w_\nv\|_{H_k^1},
	\|\na_k F_k\|_{L^2}\|\na_k\phi_\nv\|_{L^2}\big\} 
	\\&\quad
	+\|\U'''\|_{L^\infty}\|\na_k\phi_\nv\|_{L^2}^2 +\delta_0\nu^{1/3}|k|^{2/3}\mathbbm{1}_{\mathrm{Im}\lambda\leq 0}\| w_\nv\|_{L^2}^2.
	\end{split} 
\label{Na_est_4}
\end{align}
%Here, we only consider the third case (i.e., $\|\na_k F_k\|_{L^2}\|\na_k\phi_\nv\|_{L^2}$) in the minimum expression if $F_k \in H_0^1$. 
Next, we test the equation \eqref{Nav_eqn}$_1$ by $\mf J_kw_k$ and take the real part to obtain 
\begin{align}
	\mathrm{Re} \lf\lan-\nu \de_k  w_{\rm Na}
	+
	ik(\U-\mathrm{Re}\lambda-i\mathrm{Im}\lambda) w_{\rm Na}-\U''ik\phi_{\rm Na},\mf J_k  w_{\rm Na}\rg\ran 
	=
	\mathrm{Re}\lan  F_k,
	\mf J_k  w_{\rm Na}
	\ran.  \label{Na_est_5}
\end{align} 
Rearranging the terms yields
\begin{align}
	\begin{split}
		&\mathrm{Re}\lan \U ik  w_\nv, \mf J_k w_\nv\ran = \nu\mathrm{Re}\lan \de_k w_\nv,\mf J_k  w_\nv\ran+k\mathrm{Re}\lambda\; \mathrm{Re}\lf( i\lan  w_\nv,\mf J_k w_\nv\ran \rg)
		\\& -k\mathrm{Im}\lambda\;\mathrm{Re}\lan  w_\nv,\mf J_k w_\nv\ran+k\mathrm{Re}\br{\U'' i\phi_\nv, \mathfrak{J}_k  w_\nv}
		+
		\mathrm{Re}\br{F_k,\mf J_k  w_\nv}
		=:\sum_{j=1}^4 T_j.
	\end{split}
	\label{Na_est_5.4}
\end{align}
The left-hand side of  \eqref{Na_est_5.4} and the term $T_1$ can be estimates using Lemma~\ref{lem:Trans} as
\begin{align}
	&
	\mathrm{Re}\lan \U ik  w_\nv, \mf J_k w_\nv\ran\geq \frac{k^2}{16}\|\na_k\phi_\nv\|_{L^2}^2,\label{Na_est_5.41}
	\\&
	T_1
	\leq C\nu\|\na_k  w_\nv\|_{L^2}^2.\label{Na_est_5.42}
\end{align} 
Thanks to Lemma \ref{lem:AS}, we have the identity
\begin{align*}
	\overline{\int_{-1}^1 \bar{f} \mathfrak{J}_k[f] dy} 
	=
	-\int_{-1}^1 f\mathfrak{J}_k[\bar{f}]   dy
	=
	\int_{-1}^1 \bar{f} \mathfrak{J}_k[f] dy,
	\quad\forall f\in L^2. 
\end{align*} 
Hence, we deduce that $\int_{-1}^1 \bar{f} \mathfrak{J}_k[f] dy\in \rr$ and 
in particular, $\int_{-1}^{1} \overline{ w_{\rm Na}} \mathfrak{J}_k  w_{\rm Na} dy \in \rr$. 
Using Lemma~\ref{lem:BoundI}, we obtain
\begin{align}
	\begin{split}
		&T_2=k\mathrm{Re}\lambda\; \mathrm{Re}\lf( i\lan  w_\nv,\mf J_k w_\nv\ran \rg)=0,
		\\&
		|T_3|
		\leq
		|k|\lf(\delta_0\nu^{1/3}|k|^{-1/3}\mathbbm{1}_{\mathrm{Im}\lambda\leq 0}+\max\{0,\mathrm{Im}\lambda\}\rg)\| w_\nv\|_{L^2}^2.  
	\end{split}\label{Na_est_5.43}
\end{align}
Now we estimate the term $T_4$.
Based on the type of regularity, we have following cases: \begin{itemize} 
	\item If $F_k\in L^2$, we invoke Lemma~\ref{lem:BoundI} to obtain
	\begin{align}
		|T_4|&
		=
		|\lan F_k, \mf J_k  w_{\rm Na}\ran|
		\lesssim
		\|F_k\|_{L^2}\| w_\nv\|_{L^2}.
	\end{align}
	\item If $F_k \in H^1_0$, using Lemmas \ref{lem:BoundI} and~\ref{lem:H} we obtain
	\begin{align}
		\begin{split}
			|T_4|&
			=
			|\lan F_k, \mf J_k  w_{\rm Na}\ran|
			\lesssim
			|\lan |k| F_k, \mf J_k |k|\phi_{\rm Na}\ran|
			+
			|\lan F_k, \mf J_k \pa_{yy}
			\phi_{\rm Na}\ran|
			\\&
			\lesssim \|\na_kF_k\|_{L^2}\|\na_k\phi_{\rm Na}
			\|_{L^2}
			+
			|\lan F_k, \pa_y\mf J_k \pa_y \phi_{\rm Na}\ran|
			+
			|\lan F_k,[ \pa_y,\mf J_k] \pa_y \phi_{\rm Na}\ran|
			\\ &
			\lesssim
			\|\na_kF_k\|_{L^2}\|\na_k\phi_{\rm Na}
			\|_{L^2}
			+
			|\lan\pa_y F_k, \mf J_k \pa_y \phi_{\rm Na}\ran|
			\\&\qquad
			+
			\|[\mf J_k,\pa_y]\|_{L^2\rightarrow L^2}\|F_k\|_{L^2}\|\na_k\phi_{\rm Na}\|_{L^2}
			\\&
			\lesssim \|\na_kF_k\|_{L^2}\|\na_k\phi_{\rm Na}
			\|_{L^2}.
			\label{Na_est_7}
		\end{split}
	\end{align}
	\item If $F_k\in H^{-1}$, using the definition of $\mf J_k$ in \eqref{J_k} and $G_k(\pm1,y') = 0$, we get 
	\begin{align}
		\mf J_k w_k(y=\pm 1) =k\; \text{p.v.}\int_{-1}^1\frac{1}{2i}
		\cdot
		\frac{G_k(\pm 1,y') }{(\pm 1-y')}
		{w_k}(y')dy'=0.
	\end{align}
	By Lemmas~\ref{lem:BoundI},~\ref{lem:CommIpy}, and~\ref{lem:H}, we deduce
	\begin{align*}
		\|\pa_y\mf J_k w_\nv\|_{L^2}
		\leq 
		\|\mf J_k \pa_y  w_\nv\|_{L^2}+\|[\pa_y,\mf J_k] w_\nv\|_{L^2}
		\lesssim 
		\|\na_k  w_\nv\|_{L^2}.
	\end{align*}
	Hence, $\mf J_k w_k\in H^1_0$ and
	\begin{align}
		|T_4|
		\lesssim 
		\|F_k\|_{H^{-1}}\|\na_k  w_\nv\|_{L^2}.
	\end{align}
\end{itemize}
To summarize, we have the following bound for the term $T_4$:
\begin{align}
	|T_4|\lesssim \min\Bigl\{\|F_k\|_{L^2}\| w_\nv\|_{L^2},\;\|\na_kF_k\|_{L^2}\|\na_k\phi_{\rm Na} \Vert_{L^2}, \;\|F_k\|_{H^{-1}}\|\na_k  w_\nv\|_{L^2}
	\Bigr\}. \label{Na_est_5.44}
\end{align}
%We emphasize that we only consider the case when $F_k\in H^1_0$. 
% \footnote{One can also use the kernel form of the multiplier to show that $\mf J_k F_k\big|_{\pm 1}$ is bounded by $\|F\|_{L^2}$ because the Green's kernel vanishes at the boundary:
% \begin{align}
% \lf|\mf J_k f_k(y=\pm 1)\rg| =\big| k\text{p.v.}\int_{-1}^1\frac{1}{2i}\underbrace{\frac{G_k(\pm 1,y') }{(\pm 1-y')}}_{\lesssim 1?}f_k(y')dy'\big|\leq C |k|\|f_{k}\|_{L^1}.
% \end{align}} 
%From \eqref{Na_est_5},\eqref{Na_est_7}, \eqref{est_T_1}, \eqref{est_T_2}, and \eqref{est_T_3} 
From \eqref{Na_est_5.4}, \eqref{Na_est_5.41}, \eqref{Na_est_5.42}, \eqref{Na_est_5.43}, and \eqref{Na_est_5.44}, we infer that there exists a constant $C_0\geq 1$ such that 
\begin{align}
	\begin{split}
		&
		|k|^2\|\na_k\phi_{\rm Na}\|_{L^2}^2
		\leq
		{C_0}\Bigl(\nu\|\na_k w_{\rm Na}\|_{L^2}^2+
		\min\big\{\|F_k\|_{L^2}\| w_\nv\|_{L^2},\;\|F_k\|_{H^{-1}}\|\na_k  w_\nv\|_{L^2},\;
		\\&\quad
		\|\na_kF_k\|_{L^2}\|\na_k\phi_{\rm Na}
		\|_{L^2}\big\}
		%	\\&
		+
		|k|\lf(\delta_0\nu^{1/3}|k|^{-1/3}\mathbbm{1}_{\mathrm{Im}\lambda\leq 0}+\max\{0,\mathrm{Im}\lambda\}\rg)\| w_\nv\|_{L^2}^2\Bigr). 
	\end{split}
	\label{Na_est_8}
\end{align} 
We multiply \eqref{Na_est_8} by $1/2C_0$ and add the resulting inequality to \eqref{Na_est_4} to obtain that 
\begin{align}
	&\hspace{-0.25cm}\nu\|\na_k  w_\nv\|_{L^2}^2 + |k|\max\{0,\mathrm{Im}\lambda\}\| w_\nv \|_{L^2}^2+\frac{|k|^2}{2C_0}\|\na_k\phi_\nv\|_{L^2}^2\\
	\leq & 2\min\lf\{\|F_k\|_{L^2}\| w_\nv\|_{L^2},\;\|F_k\|_{H^{-1}}\| w_\nv\|_{H_k^1},\;\|\na_kF_k\|_{L^2}\|\na_k\phi_\nv\|_{L^2}\rg\} \\
	&+\|\U'''\|_{L^\infty}\|\na_k\phi_k\|_{L^2}^2 +2\delta_0\nu^{1/3}|k|^{2/3}\mathbbm{1}_{{\mathrm{Im}\lambda\leq 0}}\| w_\nv\|_{L^2}^2\\
	&+\frac{\nu}{2}\|\na_k  w_\nv\|_{L^2}^2 + \frac{|k|}{2}\max\{0,\mathrm{Im}\lambda\}\| w_\nv \|_{L^2}^2.
\end{align}
Consequently,
\begin{align}\begin{split}
		&\nu\|\na_k  w_{\rm Na}\|_{L^2}^2+ |k|\max\{0,\mathrm{Im}\lambda\}\| w_\nv \|_{L^2}^2+\frac{|k|^2}{4C_0}\|\na_k\phi_{\rm Na}\|_{L^2}^2
		\\
		&\leq 4\min\lf\{\|F_k\|_{L^2}\| w_\nv\|_{L^2},\;\|F_k\|_{H^{-1}}\| w_\nv\|_{H_k^1},\;\|\na_kF_k\|_{L^2}\|\na_k\phi_\nv\|_{L^2}\rg\} \\
		&\quad+4\delta_0\nu^{1/3}|k|^{2/3}\mathbbm{1}_{\mathrm{Im}\lambda\leq 0}\| w_\nv\|_{L^2}^2.\end{split}\label{Na_est_9}
\end{align}
Based on \eqref{Na_est_9}, we conclude that:
\begin{itemize}
	\item If $F_k\in L^2$, \eqref{Na_est_9} implies
	\begin{align}
		&\nu\|\na_k  w_{\rm Na}\|_{L^2}^2+ |k|\max\{0,\mathrm{Im}\lambda\}\| w_\nv \|_{L^2}^2+\frac{|k|^2}{4C_0}\|\na_k\phi_{\rm Na}\|_{L^2}^2
		\\
		&\leq \frac{8}{\delta_0\nu^{1/3}|k|^{2/3}}\|F_k\|_{L^2}^2+8\delta_0\nu^{1/3}|k|^{2/3}\mathbbm{1}_{\mathrm{Im}\lambda
			\leq \delta_0L^{-1}}\| w_\nv\|_{L^2}^2\\
		&\quad+\frac{1}{2}\mathbbm{1}_{\mathrm{Im}\lambda\geq \delta_0L^{-1}}|k|\max\{0,\mathrm{Im}\lambda\}\| w_\nv \|_{L^2}^2,
	\end{align}
	which leads to
	\begin{align}
		&\sqrt{\nu}\|\na_k  w_{\rm Na}\|_{L^2}+ \sqrt{|k|\max\{0,\mathrm{Im}\lambda\}}\| w_\nv \|_{L^2}+\frac{|k|}{2\sqrt{C_0}}\|\na_k\phi_{\rm Na}\|_{L^2}
		\\
		&\leq C\delta_0^{-1/2}\nu^{-1/6}|k|^{-1/3}\|F_k\|_{L^2}
		+
		C
		\mathbbm{1}_{\mathrm{Im}\lambda
			\leq \delta_0L^{-1}}
		\sqrt{\delta_0}\nu^{1/6}|k|^{1/3}\| w_\nv\|_{L^2},
	\end{align}
and \eqref{na_phi_FL2} follows.
	\item If $F\in H^1_0$, \eqref{Na_est_9} implies
	\begin{align}
		&\nu\|\na_k  w_{\rm Na}\|_{L^2}^2+\frac{|k|^2}{4C_0}\|\na_k\phi_{\rm Na}\|_{L^2}^2
		\\
		&\leq 32C_0|k|^{-2}\|\na_kF_k\|_{L^2}^2+\frac{|k|^2}{8C_0}\|\na_k\phi_\nv\|_{L^2}^2
		+
		4\delta_0\nu^{1/3}|k|^{2/3}
		\mathbbm{1}_{\mathrm{Im}\lambda
			\leq 0}
		\| w_\nv\|_{L^2}^2,
	\end{align}%33
and \eqref{na_phi_FH1} follows.
	
	\item If $F_k\in H^{-1}$, from \eqref{Na_est_9} it follows that
	\begin{align}
		&\nu\|\na_k  w_{\rm Na}\|_{L^2}^2+\frac{|k|^2}{4C_0}\|\na_k\phi_{\rm Na}\|_{L^2}^2
		\\
		&\leq 8\nu^{-1}\|F_k\|_{H^{-1}}^2+\frac{\nu}{2}\|\na_k w_\nv\|_{L^2}^2
		+
		4\delta_0\nu^{1/3}|k|^{2/3}
		\mathbbm{1}_{\mathrm{Im}\lambda
			\leq 0}
		\| w_\nv\|_{L^2}^2,
	\end{align}
	which yields \eqref{na_phi_FH-1}.
\end{itemize}

Next we estimate the spectral vorticity $ w_\nv$ following the strategy in  \cite{CWZ23,ChenWeiZhang20}.
%\myb{HS: If necessary, we can also be hand-wavy here.} %To control the nonlocal term $\U''ik\phi_\nv$ in the expression, one needs to obtain a sharp inviscid damping estimate, i.e., an estimate on the spectral stream function that is independent of $\nu^{-1}$. Among all the estimates presented in the proposition, the only candidate is the $H^1_0$-source estimate \eqref{F_in_H1_est}. As a result, we start from there. 
We introduce the following test function $\chi_1\in C^\infty([-1,1])$:
\begin{align}
	\chi_1 (y)
	=
	\frac{1}{(\U(y)-\mathrm{Re} \lambda)+iL^{-1}}.
\end{align}
%\myb{HS: I found that choosing the $+$-sign is better. It is different from the setup in the \cite{ChenWeiZhang20}. However, since it is the conjugate of the $\chi_1$ in the paper, the estimates developed hold.} 
It is readily checked that
\begin{align}\begin{split}
		&L^{-1}\|\chi_1\|_{L^\infty}+L^{-\frac{1}{2}}\|\chi_1\|_{L^2}+L^{-2}\|\partial_y\chi_1\|_{L^\infty}\leq C,
		\\
		&\|\chi_1 f\|_{L^2}\leq \|\chi_1\|_{L^2}\|f\|_{L^\infty}\leq CL^{\frac{1}{2}}\|f\|_{L^2}^{\frac{1}{2}}\|\na_k f\|_{L^2}^{\frac{1}{2}},\quad \forall f\in H_0^1([-1,1]). \end{split}\label{pro_chi_1}
\end{align} 
%\myb{For checking purpose only: Since we assume that $\U$ is strictly monotone on the interval $[-1,1]$, we have that 
%	\begin{align}
%		\|\chi_1\|_{L^2}^2
%		&=\lf|\int_{-1}^1\frac{\U'}{(\U(y)-\mathrm{Re}\lambda)^2+L^{-2}} \frac{1}{\U'}dy\rg|
%		\\&
%		\leq \|(\U')^{-1}\|_{L^\infty_y}\lf|L\int_{-1}^{1}
%		\pa_y\lf(\arctan\lf({L}(\U(y)-\mathrm{Re}\lambda)\rg)\rg)dy\rg| 
%		\\&
%		\leq  \pi\|(\U')^{-1}\|_{L^\infty}L.
%\end{align}}
We test the equation \eqref{Nav_eqn}$_1$ by $\chi_1  w_\nv$ and taking the imaginary part, obtaining
\begin{align}
\begin{split}
		&\mathrm{Im}\Big(ik\int (\U(y)-\lambda)\overline{\chi_1}| w_\nv|^2 dy\Big)
	\\ 
	&=-\nu\mathrm{Im} \int \na_k  w_\nv\cdot \overline{\na_k(\chi_1  w_\nv)} dy+\mathrm{Im}\int F_k\overline{\chi_1  w_\nv}dy
	\\&\qquad
	+\mathrm{Im}\Big(ik\int\U'' \phi_\nv\overline{\chi_1  w_\nv} dy\Big)
	\\&
	=:T_5+T_6+T_7.\label{Nav_est_15}
\end{split}
\end{align}
% If we choose the $\chi_1$ with $-iL$: Assume that $|\min\{0,\mathrm{Im}\lambda\}|\nu^{-\frac{1}{3}}|k|^{\frac{1}{3}}\leq\delta_0\leq \frac{1}{2}$.
% We estimate each term in the relation \eqref{Nav_est_15}. For the left-hand side of the relation \eqref{Nav_est_15}, we implement direct computation together with the estimates \eqref{Na_est_1}, \eqref{pro_chi_1} to obtain that for $\mathrm{Im}\lambda\leq0$,
For the left-hand side of the relation \eqref{Nav_est_15}, direct computation and \eqref{pro_chi_1}$_2$ give
\begin{align}\begin{split}
		\eqref{Nav_est_15}_{\rm L.H.S}=&k\mathrm{Re}\int | w_\nv|^2\frac{(\U-\mathrm{Re}\lambda)^2+L^{-1}\mathrm{Im}\lambda}{(\U-\mathrm{Re}\lambda)^2+L^{-2}}dy\\
		=&k\| w_\nv\|_{L^2}^2+kL^{-1}(\mathrm{Im}\lambda-L^{-1})\| 
		\chi_1  w_\nv\|_{L^2}^2\\
		\geq&k\| w_\nv\|_{L^2}^2-C\mathbbm{1}_{\mathrm{Im}\lambda \leq L^{-1}}kL^{-1}\| w_\nv\|_{L^2}\|\na_k w_\nv\|_{L^2}\\
		\geq& \frac{k}{2}\| w_\nv\|_{L^2}^2 -C\mathbbm{1}_{\mathrm{Im}\lambda \leq L^{-1}}|k|^{\frac{1}{3}}\nu^{\frac23}\|\na_k w_\nv\|_{L^2}^2.\end{split}\label{Nav_est_16}
\end{align}
%Here, in the last line, we used that $L^{-1}=\nu^{1/3}|k|^{-1/3}.$
% \begin{align}
%     \eqref{Nav_est_15}_{\rm L.H.S}=&k\mathrm{Re}\int | w_\nv|^2\frac{(\U-\mathrm{Re}\lambda)^2-L^{-1}\mathrm{Im}\lambda}{(\U-\mathrm{Re}\lambda)^2+L^{-2}}dy\\
%     =&k\| w_\nv\|_{L^2}^2-kL^{-1}(L^{-1}+\mathrm{Im}\lambda)\|\chi_1  w_\nv\|_{L^2}^2\\
%     \geq&k\| w_\nv\|_{L^2}^2-CkL^{-1}(1-\delta_0)\| w_\nv\|_{L^2}\|\na_k w_\nv\|_{L^2}\\
%     \geq& \frac{k}{2}\| w_\nv\|_{L^2}^2 -C|k|^{\frac{1}{3}}\nu^{\frac23}\|\na_k w_\nv\|_{L^2}^2.\label{Nav_est_16}
% \end{align}
For the term $T_5$, we use \eqref{pro_chi_1}$_1$ to get
\begin{align}
	\begin{split}
		|T_5|\leq& \nu\|\na_k  w_\nv\|_{L^2}(\|\chi_1'\|_{L^\infty}\| w_\nv\|_{L^2}+\|\chi_1\|_{L^\infty}\|\na_k  w_\nv\|_{L^2})\\
		\leq&C\nu^{\frac13}|k|^{\frac{2}{3}}\|\na_k  w_\nv\|_{L^2}\| w_\nv\|_{L^2}+C\nu^{\frac{2}{3}}|k|^{\frac{1}{3}}\|\na_k  w_\nv\|_{L^2}^2\\
		\leq& C\delta_0^{-\frac{1}{4}}\nu^{\frac{2}{3}}|k|^{\frac{1}{3}}\|\na_k  w_\nv\|_{L^2}^2+ {\delta_0^{\frac{1}{4}}}|k|\| w_\nv\|_{L^2}^2. \end{split}\label{Nav_est_17}
\end{align}
Similarly, the term $T_6$ is estimated using  \eqref{pro_chi_1} as follows:
\begin{align*}
	|T_6|\leq 
	\left\{
	\begin{aligned} &\| w_\nv\|_{L^2}\|\chi_1 F_k\|_{L^2}
		\leq C\nu^{-\frac{1}{3}}|k|^{\frac{1}{3}}
		\| w_\nv\|_{L^2}
		\|F_k\|_{L^2},\quad 
		\text{if~}
		 F_k\in L^2;\\ & \| w_\nv\|_{L^2}\|\chi_1 F_k\|_{L^2}\leq C\nu^{-\frac{1}{6}}|k|^{-\frac{1}{3}}\| w_\nv\|_{L^2}
		\|\na_k F_k\|_{L^2},\quad
		\text{if~} F_k\in H_0^1;\\
		&\|\na_k(\chi_1  w_\nv)\|_{L^2}
		\|F_k\|_{H^{-1}}
		\leq C
		\lf(\nu^{-\frac{2}{3}}|k|^{\frac{2}{3}}\| w_\nv\|_{L^2}
		+
		\nu^{-\frac{1}{3}} |k|^{\frac{1}{3}}
		\|\na_k w_\nv\|_{L^2}\rg)
		\|F_k\|_{H^{-1}},\\&\quad \text{if~}
		F_k\in H^{-1}.
	\end{aligned}
	\right.
\end{align*}
Finally, we estimate the term $T_7$ using \eqref{na_phi_FL2}, \eqref{na_phi_FH1}, \eqref{na_phi_FH-1}, and the bound $\|\U''\|_{L^\infty}\leq \kappa$ as
\begin{align}
	\begin{split}
		|T_7|\leq& C\|\U''\|_{L^\infty}|k|\| w_\nv\|_{L^2}\|\phi_\nv \chi_1\|_{L^2}
		\leq C\|\U''\|_{L^\infty}|k|\| w_\nv\|_{L^2}\nu^{-\frac{1}{6}}|k|^{-\frac{1}{3}}\|\na_k\phi_\nv\|_{L^2}\\
		\leq&
		C\kappa \nu^{-\frac{1}{6}}|k|^{-\frac{1}{3}}\| w_\nv\|_{L^2}\min\{ \delta_0^{-\frac{1}{2}} \nu^{-\frac{1}{6}}|k|^{-\frac{1}{3}}\|F_k\|_{L^2},|k|^{-1}\|\na_kF_k\|_{L^2},
		\\&\qquad
		\nu^{-\frac{1}{2} }\|F_k\|_{H^{-1}}\}
%		\\&
		+C\kappa\delta_0^{\frac{1}{2}}\| w_\nv\|_{L^2}^2.
	\end{split}\label{Nav_est_18}
\end{align}%\min\Bigl\{ \; |k|^{-1}\|\na_k F_k\|_{L^2},\;\Bigr\}
%Finally, we combine \eqref{Nav_est_15}, \eqref{Nav_est_16}, \eqref{Nav_est_17}, \eqref{Nav_est_18}, and the estimate \eqref{Na_est_1} to obtain that 
Combining \eqref{na_phi_FL2}, \eqref{na_phi_FH1}, \eqref{na_phi_FH-1}, \eqref{Nav_est_15}, \eqref{Nav_est_16}, \eqref{Nav_est_17}, and \eqref{Nav_est_18}, we obtain
\begin{itemize}
		\item If $F_k\in L^2$, 
	\begin{align*}
		\| w_\nv\|_{L^2}^2 &\lesssim \delta_0^{-\frac{1}{4}}\nu^{\frac{2}{3}}|k|^{-\frac{2}{3}}\|\na_k  w_\nv\|_{L^2}^2+\delta_0^{\frac{1}{4}}\| w_\nv\|_{L^2}^2
		+\nu^{-\frac{1}{3}}|k|^{-\frac{2}{3}}\| w_\nv\|_{L^2}\|F_k\|_{L^2}		\\
		&\quad
		+\kappa \nu^{-\frac{1}{3}}|k|^{-\frac{5}{3}}\| w_\nv\|_{L^2}\|F_k\|_{L^2}+\kappa\delta_0^{\frac{1}{2}}|k|^{-1}\| w_\nv\|_{L^2}^2\\
		&   \lesssim \delta_0^{-\frac{5}{4}}\nu^{-\frac{2}{3}}|k|^{-\frac{4}{3}}\|F_k\|_{L^2}^2+\delta_0^{\frac{1}{4}}\| w_\nv\|_{L^2}^2.
	\end{align*}
	We take $\delta_0>0$ sufficiently small to obtain
	\begin{align}
		\nu^{\frac{1}{3}}|k|^{\frac{2}{3}}\| w_\nv\|_{L^2}\lesssim \|F_k\|_{L^2}.
	\end{align}  
	The estimates of $\mathbb{E} [ w_{\rm Na}, \phi_{{\rm Na}}]$ and $\sqrt{|k| ({\rm Im} \lambda)_+} \Vert  w_{\rm Na}\Vert_{L^2}$ in  \eqref{F_in_L2_est} follows by appealing to \eqref{na_phi_FL2}.
	
	\item If $F_k\in H^1_0$, 
	\begin{align*}
		\| w_\nv\|_{L^2}^2 &\lesssim 
		\delta_0^{-\frac{1}{4}}\nu^{\frac{2}{3}}|k|^{-\frac{2}{3}}\|\na_k  w_\nv\|_{L^2}^2+\delta_0^{\frac{1}{4}}\| w_\nv\|_{L^2}^2
%		\\&\quad
		+
		\nu^{-\frac{1}{6}}|k|^{-\frac{4}{3}}
		\| w_\nv\|_{L^2}\|\na_kF_k\|_{L^2}
		\\&\quad
		+\kappa \nu^{-\frac{1}{6}}|k|^{-\frac{7}{3}}\| w_\nv\|_{L^2}\|\na_kF_k\|_{L^2}+\kappa\delta_0^{\frac{1}{2}}|k|^{-1}\| w_\nv\|_{L^2}^2\\
		&   \lesssim \delta_0^{-\frac{1}{4}}\nu^{-\frac{1}{3}}|k|^{-\frac{8}{3}}\|\na_kF_k\|_{L^2}^2+\delta_0^{\frac{1}{4}}\| w_\nv\|_{L^2}^2.
	\end{align*}
	We take $\delta_0>0$ sufficiently small to get
	\begin{align}
		\label{Nav_est_19}
		\nu^{\frac{1}{6}}|k|^{\frac{4}{3}}\| w_\nv\|_{L^2}
		\les
		\|\na_kF_k\|_{L^2}. 
	\end{align}
	From \eqref{Nav_est_19} and \eqref{na_phi_FH1} it follows that
	\begin{align*}
		|k|\|\na_k\phi_\nv\|_{L^2}+\sqrt{\nu}\|\na_k w_\nv\|_{L^2}+\nu^{\frac{1}{6}}|k|^{\frac{1}{3}}\| w_\nv\|_{L^2}
		\les
		|k|^{-1}\|\nabla_kF_k\|_{L^2},
	\end{align*}
	 and the estimates of $\mathbb{E} [ w_{\rm Na}, \phi_{{\rm Na}}]$ and $\sqrt{|k| ({\rm Im} \lambda)_+} \Vert  w_{\rm Na}\Vert_{L^2}$ in \eqref{F_in_H1_est} follows.

	\item If $F_k\in H^{-1}$,
	\begin{align*}
		\| w_\nv\|_{L^2}^2 
		&\lesssim \delta_0^{-\frac{1}{4}}\nu^{\frac{2}{3}}|k|^{-\frac{2}{3}}\|\na_k  w_\nv\|_{L^2}^2+\delta_0^{\frac{1}{4}}\| w_\nv\|_{L^2}^2
		+\nu^{-\frac{2}{3}}|k|^{-\frac{1}{3}}\| w_\nv\|_{L^2}\|F_k\|_{H^{-1}}		\\&\quad
		+\nu^{-\frac13}|k|^{-\frac{2}{3}}\|F_k\|_{H^{-1}} \Vert \nabla_k  w_{\rm Na}\Vert_{L^2}
%		\\&\quad
		+
		\kappa\delta_0^{\frac{1}{2}}|k|^{-1}\| w_\nv\|_{L^2}^2\\
		&   \lesssim \delta_0^{-\frac{1}{4}}\nu^{-\frac{4}{3}}|k|^{-\frac{2}{3}}\|F_k\|_{H^{-1}}^2+\delta_0^{\frac{1}{4}}\| w_\nv\|_{L^2}^2.
	\end{align*}
Similar arguments as in other two cases yield the estimates of $\mathbb{E}[ w_{\rm Na}, \phi_{{\rm Na}}]$ and the term $\sqrt{|k| ({\rm Im} \lambda)_+} \Vert  w_{\rm Na}\Vert_{L^2}$ in \eqref{F_in_H-1_est}.
\end{itemize}

To conclude the proof of the proposition, it remains to prove
\begin{align}
	\nu\|\de_k  w_\nv\|_{L^2}\lesssim \min\{\|F_k\|_{L^2},\nu^{1/6}|k|^{-2/3}\|\na_kF_k\|_{L^2}\}.\label{Nav_est_20}
\end{align}
We test the equation \eqref{Nav_eqn}$_1$ by $\nu \de_k  w_\nv$ and take the real part to obtain
%{\color{cyan}For Checking purpose only: Denote $\na_k:=(ik,\pa_y)$, then 
%	\begin{align*}
%		&-\nu \Delta_k  w_\nv+i k\left(\U-\mathrm{Re}\lambda-i \mathrm{Im}\lambda\right)  w_\nv-i \U^{\prime \prime} k \phi_\nv=F_k\\
%		\Longrightarrow&-\int\nu \Delta_k  w_\nv \overline{\Delta_k  w_\nv}+ i k\int\left(\U-\mathrm{Re}\lambda\right)  w_\nv \overline{\Delta_k  w_\nv}\\
%		&\hspace{2.5cm} +k \mathrm{Im}\lambda\int  w_\nv \overline{\Delta_k  w_\nv}-i\int\U'' k\phi_\nv\overline{\de_k  w_\nv}=\int F_k\overline{\de_k w_\nv}\\ 
%		\Longrightarrow& \nu\|\de_k  w_\nv\|_{L^2}^2+k\mathrm{Im}\lambda \|\na_k  w_\nv\|_{L^2}^2+\mathrm{Re}\ \int ik \U'  w_\nv\overline{\pa_y  w_\nv}\\
%		&-\mathrm{ Re}\ i \int \U'''k\phi_\nv \overline{\pa_y w_\nv}-\mathrm{Re}\int i\U'' k\na_k\phi_\nv\cdot\overline{\na_k  w_\nv}=\begin{cases}\mathrm{Re}\int \na_kF_k\cdot \overline{\na_k  w_\nv} ,\quad F_k\in H^1_0,\\ -\mathrm{Re}\int F_k \overline{\de_k w_\nv},\quad F_k\in L^2.\end{cases}
%	\end{align*}
%	Now, we multiply the relation by $\nu$, recall that $\mathrm{Im}\lambda\geq -\delta_0\nu^{1/3}|k|^{-1/3}$, and apply the established estimate (from Step \# 1, \# 2)
%	\begin{align*}
%		\mathbb{E}[ w_\nv,\phi_\nv]=&\nu^{1/6}|k|^{4/3}\|\na_k\phi_\nv\|_{L^2}+\nu^{2/3}|k|^{1/3}\|\na_k  w_{\rm Na}\|_{L^2}+\nu^{1/3}|k|^{2/3}\| w_\nv\|_{L^2}\n\\
%		\hspace{2.5cm}\leq& C_{\delta_0} \min\lf\{\|F_k\|_{L^2},
%		\nu^{1/6}|k|^{-2/3}\|\na_k F_k\|_{L^2}\rg\},
%	\end{align*}
%}
\begin{align}
\begin{split}
		& \nu^2\|\de_k w_\nv\|_{L^2}^2+\nu|k|(\mathrm{Im} \lambda)_+\|\na_k  w_\nv\|_{L^2}^2\\
	&\leq |(\mathrm{Im}\lambda)_-|
\nu|k|\|\na_k  w_\nv\|_{L^2}^2
	\\&
	+
	C\nu|k|\|\U'\|_{L^{\infty}}\| w_\nv\|_{L^2}\|\na_k  w_\nv\|_{L^2}+C\nu|k|\|\U''\|_{W^{1,\infty}}\|\na_k\phi_\nv\|_{L^2}\|\na_k  w_\nv\|_{L^2}\\
	&\quad+C\min\Big\{\nu\|F_k\|_{L^2}\|\de_k  w_\nv\|_{L^2}, \nu \|\na_k F_k\|_{L^2}
	\|\na_k  w_\nv\|_{L^2}\Big\}\\
	&\leq C \mathbb{E}[ w_\nv,\phi_\nv]^2+C\min\{\|F_k\|_{L^2}^2,\nu^{1/3}|k|^{-4/3}\|\na_kF_k\|_{L^2}^2\} +\frac{1}{2}\nu^2\|\de_k w_\nv\|_{L^2}^2,
	\label{EQ200a}
\end{split}
\end{align}
where we used the relation $\nu^{1/3}|k|^{-1/3} \leq \nu^{1/6}|k|^{-2/3}\|\U'\|_{L^\infty}^{1/6}$ in the last step, which follows from $\nu|k|^2\leq \|\U'\|_{L^\infty}$. 
From the already-established bound
\begin{align}
	\mathbb{E}[ w_\nv,\phi_\nv]^2
	\lesssim \min\{\|F_k\|_{L^2}^2,\nu^{1/3}|k|^{-4/3}\|\na_kF_k\|_{L^2}^2\},
\end{align} 
we deduce \eqref{Nav_est_20} from \eqref{EQ200a}.
\end{proof}

\begin{lemma}
For $\mathcal{E}$ defined in \eqref{Frz_E}, we have
\label{LemmaE}
\begin{align}
	\|e^{\delta_\ast\nu^{1/3}t} \omega\|_{L^2([0,T];\mathcal{Z})}^2
	\lesssim 
	\mathcal{E},
	\label{EQ123c}
\end{align}
where $\delta_\ast>0$ is the constant from Proposition~\ref{pro:lin_dy}.
\end{lemma}

\begin{proof}
Using \eqref{Frz_om_dmp}, we have
\begin{align}
\begin{split}
	&
		\|e^{\delta_\ast\nu^{1/3}t} \omega\|_{L^2([0,T];\mathcal{Z})}^2
	=
	\sum_{j=0}^N
	\int_{\mathcal{I}_{[j]}}
	e^{2\delta_\ast \nu^{1/3} t}
	\Big\| \sum_{j'=0}^j  \omega_{[j']}
	\Big\|^2_{\mathcal{Z}}
%	\\&\quad
	\leq
	\sum_{j=0}^N
	e^{2\delta_\ast (j+1)}
	\int_{\mathcal{I}_{[j]}}
	\lf(\sum_{j'=0}^{j} 
	\Vert  \omega_{[j']}\Vert_{\mathcal{Z}}\rg)^2.
%		\les
%	\sum_{j=0}^N
%	e^{2\delta_\ast (j+1)}
%	\lf(\sum_{j'=0}^{j} 
%	\Vert  \omega_{[j']}\Vert_{L^2 (\mathcal{I}_{[j]};\mathcal{Z})}\rg)^2,
%		\\&\quad
%		\les
%		\sum_{j=0}^N
%		e^{6\delta_\ast (j+1)}
%		\int_{\mathcal{I}_{[j]}}
\end{split}
\end{align}
An application of the Cauchy-Schwarz inequality leads to
\begin{align}
	\begin{split}
		\int_{\mathcal{I}_{[j]}}
	\lf(\sum_{j'=0}^{j} 
	\Vert  \omega_{[j']}\Vert_{\mathcal{Z}}\rg)^2
	&\leq
	\int_{\mathcal{I}_{[j]}}
	\lf(\sum_{j'=0}^{j}
	(j-j'+1)^2
	\Vert  \omega_{[j']}\Vert_{\mathcal{Z}}^2\rg)
	\lf(\sum_{j'=0}^{j}
	\frac{1}{(j-j'+1)^2}\rg)
	\\&
	\leq
	C\lf(\sum_{j'=0}^{j}
	(j-j'+1)
	\Vert  \omega_{[j']}\Vert_{L^2 (\mathcal{I}_{[j]};\mathcal{Z})}\rg)^2
	=C Y_{[j]}^2.
	\end{split}
\end{align}
Therefore, we conclude the proof of the lemma.
\end{proof}

%Below, we present the technical lemmas. 
\begin{lemma}
\label{lem:Frz_cf}
Assume that $U_{\rm in}\in H^4 (-1,1)$ satisfies the compatibility condition \eqref{eq:bc:Uin}. 
Then for $U(t,y)$ defined in \eqref{def:barU}, we have
\begin{align}
	\Vert U (t,\cdot)-y\Vert_{H^4}
	&\les \Vert U_{\rm in} -y\Vert_{H^4},
	\quad
	\forall t\geq 0,
	\label{Prfl_est0}
	\\
	\sup_{y\in [-1,1]}
	\lf|\frac{U(t,y)-U(s,y)}{1-|y|}\rg|
	&\les\nu \|U_{\rm in}\|_{H^4} |t-s|,
		\label{Prfl_est}
	\\
	\Vert U''(t,\cdot) - U''(s,\cdot)\Vert_{L^2_y}
	&\les
	 \nu \|U_{\rm in}\|_{H^4}|t-s|,
	\label{Prfl_est2}
\end{align}
for any $t\geq s\geq 0$.
\end{lemma}
\begin{proof}
%	The proof is based on the property of $U$ and will be decomposed into two steps. {\color{red} HS: Here, the proof requires $C^4$ regularity of $U_{\rm in}$. But I guess that $C^3$ is enough. If one really want to reduce the regularity of the $U_{\rm in}$, then we might have
%\begin{align} \lf|\frac{U(t,y)-U(s,y)}{(1-|y|)^{1/2}}\rg|\leq C(\|U_{\rm in}\|_{\color{red}C^{2,1/2}})\nu |t-s|. 
%\end{align}}
From \eqref{def:barU} it follows that $(U-y)$ solves heat equation
\begin{align}
	\left\{
	\begin{aligned}
		&\pa_t (U-y)=\nu\pa_{y}^2(U-y),
		\\ &
		(U-y)\big|_{y=\pm 1}=0,
	\quad
		(U-y) (t=0)= U_{\rm in}  -y.
	\end{aligned}
	\right.
	\label{heateqn}
\end{align}
A standard basic energy estimate yields \eqref{Prfl_est0}.
Hence, $\pa_y^2(U-y)|_{y=\pm 1}=0$ and
\begin{align}
	\int_{-1}^1 \pa_y (U-y) dy=(U-y)\Big|_{y=-1}^{y=1}
	=0.
\end{align}
By the Agmon inequality and Poincar\'e inequality, we have
\begin{align}
	\begin{split}
		\|\pa_y (U-y)\|_{L^\infty}
		&
		\lesssim
		\|\pa_y (U-y)\|_{L^2}^{1/2}\|\pa_y(U-y)\|_{H^1}^{1/2}
		%	\\&
		\lesssim   
		\| \partial_y^2 (U-y)\|_{L^2},
		\label{EQ122a}
	\end{split}
\end{align}
for any $t\geq 0$.
Differentiating \eqref{heateqn}
gives
\begin{align}
	\left\{
	\begin{aligned}
		&
		\pa_t \pa_y^2(U-y)=\nu \pa_{y}^2 (\pa_y^2(U-y)),
		\\& \pa_y^2(U-y)\Big|_{y=\pm 1}=0, \quad
		\partial_y^2 (U-y) (t=0)
		=\partial_y^2 U_{\rm in}.	
	\end{aligned}
	\right.
	\label{EQ200b}
\end{align}
Using \eqref{EQ122a} and the basic energy inequality of the heat equation \eqref{EQ200b}, we infer
\begin{align}
	\Vert \partial_y (U-y)\Vert_{L^\infty_{t} L^\infty_y}
	\les
	\Vert \partial_y^2 (U-y)\Vert_{L^\infty_{t} L^2_y}
	\les
	\Vert \partial_{y}^2 U_{\rm in}\Vert_{L^2_y}
	\les
	\Vert U_{\rm in}\Vert_{H^2}.
\end{align}
Similar argument yields 
\begin{align}
	\|\pa_y^3(U-y)\|_{L_{t}^\infty L^\infty_y}
	\les
	\Vert \pa_y^4 (U-y)\Vert_{L_{t}^\infty L^2_y}
	\les
	\|\partial_{y}^4 U_{\rm in}\|_{L^2_y}
	\les
	\Vert U_{\rm in}\Vert_{H^4}.
	\label{EQ122b}
\end{align}
Let $t\geq s\geq 0$.
From the fundamental theorem of calculus it follows that
\begin{align}
	U(t,y)- U(s,y)
	=
	\int_{s}^t
	\partial_{\tau} 
	U(\tau, y) d\tau
	=
	\nu
	\int_{s}^t
	\partial_{y}^2
	U(\tau,y) d\tau,
	\quad
	y\in [-1,1].
	\label{EQ200c}
\end{align}
%
%\begin{align}
%	\lf|\frac{U(t,y)-U(s,y)}{1-|y|}\rg|
%	=
%	\lf|\frac{1}{1-|y|}\int_s^t \pa_\tau U(\tau,y)d\tau\rg|
%	=
%	\lf|\nu\int_s^t  \frac{\pa_{y}^2U(\tau,y)-0}{1-|y|}d\tau\rg|.
%\end{align}
Using \eqref{EQ122b}, \eqref{EQ200c}, and the mean value theorem, we arrive at
\begin{align}
\begin{split}
		&\lf|\frac{U(t,y)-U(s,y)}{1-y}\rg|
	\leq
	\nu 
	\int_s^t 
	\lf|
	\frac{\partial_y^2 U(\tau, y)  - \partial_y^2 U(\tau, 1)}{y-1}
	\rg|
	d\tau
	\\&\quad
	\leq
	\nu 
	\int_s^t
	\sup_{y\in [-1,1]}
	|\partial_y^3 U(\tau, y) |
	d\tau
	\leq C\nu \|U_{\rm in}\|_{H^4}|t-s|,
	\label{EQ200d}
\end{split}
\end{align}
for any $y\in[-1,1]$.
A similar argument shows that
\begin{align}
	\lf|\frac{U(t,-y)-U(s,-y)}{1-y}\rg|
%	\leq
%	\nu 
%	\int_s^t 
%	\lf|
%	\frac{\partial_y^2 U(\tau, -y)  - \partial_y^2 U(\tau, -1)}{-y-(-1)}
%	\rg|
%	d\tau
	\leq C\nu \|U_{\rm in}\|_{H^4}|t-s|,
\end{align}
for any $y\in [-1,1]$.
Hence, the proof of \eqref{Prfl_est} is completed.

To prove \eqref{Prfl_est2}, we proceed similarly as in \eqref{EQ200d} to obtain
\begin{align}
	\begin{split}
		&\Vert
		U''(t,\cdot)- U''(s,\cdot) \Vert_{L^2_y}
		\leq
		\int_s^t 
		\Vert
		\partial_t U''(\tau, \cdot)  
		\Vert_{L^2_y}
		d\tau
		\leq
		\nu
		\int_s^t 
		\Vert
		\partial_y^4 U(\tau, \cdot)  
		\Vert_{L^2_y}
		d\tau
		\\&\quad
		\leq
		\nu 
		\int_s^t
	\Vert \pa_y^4 U_{\rm in} \Vert_{L^2}
		d\tau
		\leq \nu \|U_{\rm in}\|_{H^4}|t-s|,
	\end{split}
\end{align}
where we also used the basic energy estimate for the one-dimensional heat equation.
\end{proof}

\begin{lemma}\label{lem:Frzfrc}
Suppose that $\nu k^2 \les 1$.
Then we have
\begin{align}
	\label{Frzfrc_est}
	&\|\mathbbm{f}_{{\rm Frozen} [j]}\|_{L^2(\mathcal{I}_\jj;L^2)}
%		\lesssim \myr{\nu^{1/2}|k|^{-1/3}}\sum_{j'=0}^{j-1}(|j-j'|+1)\| \omega_{\jjj}\|_{L^2(\mathcal{I}_{\jj};\mathcal{Z})}
	\lesssim \nu^{1/2}|k|^{-1/3}Y_{\jj},
	\quad j\geq 1,
	\\&
	\sum_{j=0}^N
|\log \nu|^2
\nu^{-1/3} |k|^{4/3}
\Vert e^{\delta_\ast \nu^{1/3} t} \mathbbm{f}_{{\rm Disc} [j]}\Vert_{L^2 (\mathcal{I}_{[j]}; L^2)}^2
\les
|\log \nu|^2
\nu^{1/3}
\EE.
	\label{EQ123d}
\end{align}
%\begin{align}
%	\sum_{j=0}^N 
%	e^{2\delta_\ast j}\|\mathbbm{f}_{{\rm Frozen}}\|_{L^2(\mathcal{I}_\jj;L^2)}^2\lesssim \nu|k|^{-2/3}E.
%\end{align}
\end{lemma}
\begin{proof}
Using \eqref{frozenj}, we write
\begin{align}
\begin{split}
		\mathbbm{f}_{{\rm Frozen} [j]}
	&
	=-
	\sum_{j'=0}^{j-1}(U_\jj-U_\jjj) \omega_\jjj
	+
	\sum_{j'=0}^{j-1}(U_\jj''-U_\jjj'')\psi_\jjj
	\\&
	=:
	\sum_{j'=0}^{j-1}\mathbbm{f}_{{{\rm Frozen}}\jjj}^{(1)}
	+
	\sum_{j'=0}^{j-1}\mathbbm{f}_{{\rm Frozen}\jjj}^{(2)}.
\end{split}
	\label{Frzfrcest1}
\end{align}
%Now we estimate each term in \eqref{Frzfrcest1}. 
For $0\leq j'\leq j-1$,
using Lemma~\ref{lem:Frz_cf} we obtain
\begin{align}
	\label{Frzfrcest2}
\begin{split}
		&\|\mathbbm{f}_{{\rm Frozen}\jjj}^{(1)}\|_{L^2(\mathcal{I}_\jj;L^2)}
	= \lf\|\frac{U_\jj-U_{\jjj}}{(1-|y|)^{1/2}}(1-|y|)^{1/2} \omega_{\jjj}\rg\|_{L^2_t(\mathcal{I}_{\jj};L^2)}
	\\&
	\lesssim \lf\|\frac{U_\jj-U_{\jjj}}{(1-|y|)^{1/2}}\rg\|_{L^\infty_y}
	\| \rho_k^{1/2}
	 \omega_{\jjj} \|_{L^2(\mathcal{I}_{\jj};L^2)}
%	\\&
%		\lesssim \myr{\nu^{2/3}|j-j'|\nu^{-1/6}|k|^{-1/3}\| \omega_{\jjj}\|_{L^2(\mathcal{I}_{\jj};\mathcal{Z})}
	\lesssim
	\nu^{1/2} |k|^{-1/3}
|j-j'| \| \omega_{\jjj}\|_{L^2(\mathcal{I}_{\jj};\mathcal{Z})}.
\end{split}
\end{align}
From Lemma~\ref{lem:Frz_cf} it follows that
\begin{align}
	\label{Frzfrcest3}
\begin{split}
		&\|\mathbbm{f}_{{\rm Frozen}\jjj}^{(2)}\|_{L^2(\mathcal{I}_\jj;L^2)}
	= \lf\|(U_\jj''-U_{\jjj}'')\psi_{\jjj}\rg\|_{L^2(\mathcal{I}_\jj;L^2)}
	\\
	&\lesssim \lf\|U_\jj''-U_{\jjj}''\rg\|_{L^\infty_y}|k|^{-1}
	\|\na_k\de_k^{-1} \omega_{\jjj}\|_{L^2(\mathcal{I}_{\jj};L^2)}\\
	&\lesssim 
	\nu^{2/3}
	|k|^{-1}
	|j-j'|\| \omega_{\jjj}\|_{L^2(\mathcal{I}_{\jj};\mathcal{Z})}.
\end{split}
\end{align}
Combining \eqref{Frzfrcest1}, \eqref{Frzfrcest2}, and \eqref{Frzfrcest3}, we arrive at
\begin{align}
	\begin{split}
	\Vert \mathbbm{f}_{{\rm Frozen} [j]}\Vert_{L^2(\mathcal{I}_\jj;L^2)}
	&
	\les
	\sum_{j'=0}^{j-1}	
	(	\nu^{1/2} |k|^{-1/3}
		+
		\nu^{2/3} |k|^{-1})
	|j-j'| \| \omega_{\jjj}\|_{L^2(\mathcal{I}_{\jj};\mathcal{Z})}
%	\\&
	\les
	\nu^{1/2} |k|^{-1/3} 
	Y_{[j]}.
	\end{split}
\end{align}
completing the proof of \eqref{Frzfrc_est}.

%	\begin{align}
%		\begin{split}
%		\sum_{j=0}^N
%	|\log \nu|^2
%	\nu^{-1/3} |k|^{4/3}
%	\Vert e^{3\delta_\ast \nu^{1/3} t} \mathbbm{f}_{{\rm Disc} [j]}\Vert_{L^2 (\mathcal{I}_{[j]}; L^2)}^2,	
%		\end{split}
%	\end{align}

For $0\leq j\leq N$, we use \eqref{discj2} and \eqref{discj3} to write
\begin{align}
\begin{split}
		\mathbbm{f}_{{{\rm Disc}}[j]}
	&=
			\sum_{j'=0}^j
	(U_\jj (y)-U(t,y))
	 \omega_{[j']}
	-
			\sum_{j'=0}^j
	(U''_\jj (y)-U''(t,y))
	\psi_{[j']}.
\end{split}
	\label{discj4}
\end{align}
An application of the Cauchy-Schwarz inequality leads to
\begin{align}
	\begin{split}
		&
\lf	\Vert e^{\delta_\ast \nu^{1/3} t}
		\sum_{j'=0}^j
	(U_\jj (y)-U(t,y))
	 \omega_{[j']}
	\rg\Vert_{L^2 (\mathcal{I}_{[j]}; L^2)}^2
	\\&
	\les
	\int_{\mathcal{I}_{\jj}}
		\lf(	\sum_{j'=0}^j	\Vert e^{\delta_\ast \nu^{1/3} t}
	(U_\jj (y)-U(t,y))
	 \omega_{[j']}
	\Vert_{L^2}\rg)^2
	\\&
	\les
		\int_{\mathcal{I}_{\jj}}
	\lf(	\sum_{j'=0}^j
	(j-j'+1)^2
		\Vert e^{\delta_\ast \nu^{1/3} t}
	(U_\jj (y)-U(t,y))
	 \omega_{[j']}
	\Vert_{L^2}^2
	\rg)
	\lf(
	\sum_{j'=0}^{j}
	\frac{1}{(j-j'+1)^2}
	\rg)
	\\&
	\les
	\sum_{j'=0}^j
	e^{2\delta_\ast j}
	(j-j'+1)^2
	\int_{\mathcal{I}_{\jj}}
	\Vert 
	(U_\jj (y)-U(t,y))
	 \omega_{[j']}
	\Vert_{L^2}^2.
	\end{split}
\label{EQ201a}
\end{align}
Similar arguments as in \eqref{Frzfrcest2} give
\begin{align}
	\begin{split}
		&
	\Vert 
	(U_\jj (y)-U(t,y))
	 \omega_{[j']}
	\Vert_{L^2 (\mathcal{I}_{\jj}; L^2 )}^2
	\les
	\nu |k|^{-2/3}
	\Vert 
	 \omega_{[j']}
	\Vert_{L^2 (\mathcal{I}_{\jj}; \mathcal{Z})}^2.
	\end{split}
\label{EQ201b}
\end{align}
From \eqref{EQ201a} and \eqref{EQ201b} it follows that
\begin{align}
	\begin{split}
		&
		\sum_{j=0}^N
	|\log \nu|^2
	\nu^{-1/3} |k|^{4/3}
	\Vert e^{\delta_\ast \nu^{1/3} t} 
		\sum_{j'=0}^j
	(U_\jj (y)-U(t,y))
	 \omega_{[j']}
	\Vert_{L^2 (\mathcal{I}_{[j]}; L^2)}^2
	\\&
	\les
	|\log \nu|^2
	\nu^{2/3} |k|^{2/3}
	\sum_{j=0}^N
		\sum_{j'=0}^j
	e^{2\delta_\ast j}
	(j-j'+1)^2
	\Vert  \omega_{\jjj}\Vert^2_{L^2 (\mathcal{I}_{\jj};\mathcal{Z})}
	\les	
		|\log \nu|^2
	\nu^{1/3}
	\EE.
\end{split}
\label{discj5}
\end{align}
A similar argument as in \eqref{Frzfrcest3} shows that
\begin{align}
		&\hspace{-0.5cm}
		\sum_{j=0}^N
		|\log \nu|^2
		\nu^{-1/3} |k|^{4/3}
		\lf\Vert e^{\delta_\ast \nu^{1/3} t} 
		\sum_{j'=0}^j
	(U''_\jj (y)-U''(t,y))
	\psi_{[j']}
		\rg\Vert_{L^2 (\mathcal{I}_{[j]}; L^2)}^2
		\les	
		|\log \nu|^2
		\nu^{1/3}
		\EE.
\label{discj6}
\end{align}
Combining \eqref{discj4}, \eqref{discj5}, and \eqref{discj6}, we conclude the proof of \eqref{EQ123d}.
\end{proof}

%\begin{proof}
%	Since $\nabla \cdot u=0$, we have
%	\begin{align}
%		\int_{-\pi}^{\pi}
%		\pa_x u^{(1)}
%		+\pa_y u^{(2)}
%		dx
%		=0
%		\implies
%			\int_{-\pi}^{\pi}
%			\pa_y u^{(2)}
%		dx
%		=0.
%	\end{align}
%	Since $u^{(2)}|_{y=\pm 1}=0$, we have
%	\begin{align}
%		\int_{-\pi}^{\pi} u^{(2)} (y=1) dx=0.
%	\end{align}
%Combining the above identity, we get
%	\begin{align}
%	\int_{-\pi}^{\pi} u^{(2)}  dx=0
%\end{align}
%for any $y\in [-1,1]$.
%\end{proof}

\vspace{4 mm}

\vspace{4 mm}

\noindent \textbf{Acknowledgements:}  JB gratefully  acknowledges support from NSF DMS-2108633. SH acknowledges support from NSF DMS-2406293. SI gratefully acknowledges support from NSF DMS2306528 and NSF CAREER DMS2442781. FW gratefully  acknowledges support from the National Natural Science Foundation of China (No. 12471223, 12101396, and 12331008).

\addcontentsline{toc}{section}{References}
\bibliographystyle{abbrv}
\bibliography{bibliography}

\end{document}